\documentclass{amsart}
\usepackage{amsmath,amssymb,amsxtra,amsthm,amscd}
\usepackage{setspace}

\newcommand{\dotr}{\mbox{$\boldsymbol{.}$}}
\usepackage{multicol}
\usepackage{bbm}
\usepackage{graphicx}
\usepackage[all,cmtip]{xy}
\usepackage{lmodern}
\usepackage{enumitem}
\usepackage{color}
\usepackage[color=white,bordercolor=red,linecolor=red]{todonotes}

\usepackage{soul}

\usepackage{comment}

\usepackage{pb-diagram, pb-xy}

\usepackage{mathtools}  
\usepackage{tabulary}
\usepackage{booktabs}

\usepackage{amssymb}
\usepackage{latexsym}
\usepackage{pb-diagram}
\usepackage[mathscr]{euscript}
\usepackage{graphicx}
\usepackage{xcolor}
\usepackage[new]{old-arrows}

\newcommand\im{\mathrm{Im}}

\renewcommand\ker{\mathrm{Ker}}

\newcommand\cat{\underline}

\DeclareMathOperator{\Sdd}{\mathit{S}}
\DeclareMathOperator{\diam}{\mathrm{diam}}

\newcommand\commentout[1]{\marginpar{\tiny $\backslash$commentout}}

\renewcommand\qed{\hfill$\square$}

\def\compcirc {\mbox{\hspace{.05cm}}\raisebox{.04cm}{\tiny  {$\circ$ }}}

\makeatletter
\newcommand{\colim@}[2]{%
  \vtop{\m@th\ialign{##\cr
    \hfil$#1\operator@font colim$\hfil\cr
    \noalign{\nointerlineskip\kern1.5\ex@}#2\cr
    \noalign{\nointerlineskip\kern-\ex@}\cr}}%
}
\newcommand{\colim}{%
  \mathop{\mathpalette\colim@{\rightarrowfill@\textstyle}}\nmlimits@
}
\makeatother

\theoremstyle{plain}
\newtheorem{Theorem}{Theorem}[section]

\newtheorem{Corollary}[Theorem]{Corollary}

\newtheorem{Lemma}[Theorem]{Lemma}
\newtheorem{Proposition}[Theorem]{Proposition}

\theoremstyle{definition}
\newtheorem{Definition}[Theorem]{Definition}

\newtheorem{Remark}[Theorem]{Remark}
\newtheorem{Example}[Theorem]{Example}

\theoremstyle{remark}
\newtheorem{Notation}[Theorem]{Notation}

\newenvironment{Proof}{\par\noindent\textbf{Proof:}}
{\qed}

\usepackage{graphics}

\DeclareMathOperator*{\hocolimprep}{hocolim}                   
\newcommand{\hocolim}[1]%
{\displaystyle\hocolimprep_{\substack{-\!\rightarrow \\ #1}} \, }

\DeclareMathOperator*{\colimprep}{colim}                   
\newcommand{\truecolim}[1]%
{\displaystyle\colimprep_{\substack{-\!\rightarrow \\ #1}} \, }

\DeclareMathOperator{\Or}{Or}

\title{Algebraic \textit{K}-theory of Geometric Groups}
\author{Gunnar Carlsson and Boris Goldfarb}
\begin{document}

\begin{abstract}
In this paper we introduce a homotopy theoretic technique for proving that the $K$-theoretic assembly map is an equivalence.  It is an extension of the methods used to prove split injectivity of the assembly and applies to any geometrically finite group.   Our result is that there are two requirements which need to hold.  The first is that the assembly map for the group regarded as a metric space is an equivalence.  This is a non-equivariant condition and depends only on the coarse type of the word metric on the group.  The second is that the group ring satisfies an algebraic coherence condition, which currently can be verified for all known groups for which the split injectivity statement for the assembly holds.  The two conditions extend very broadly. In particular, both conditions hold for groups of finite asymptotic dimension. 

To state the main theorem precisely, given a regular Noetherian ring $A$ of finite global dimension and a group $\Gamma$ with finite $K(\Gamma,1)$ and finite asymptotic dimension, we prove that the $K$-theoretic assembly map is an equivalence.  Therefore, in all dimensions the $K$-theory of $A[\Gamma]$ is the group homology of $\Gamma$ with coefficients in the $K$-theory spectrum of $A$.  One of the many geometric consequences of this theorem is vanishing of the Whitehead group of $\Gamma$.
\end{abstract}

\maketitle
\tableofcontents
\section{Introduction}

This work is motivated by the Borel conjecture, also known as the topological rigidity conjecture for aspherical manifolds, first stated as a problem by Borel in a letter to Serre in 1953 \cite{aR:16}.  
The conjecture states that whenever two compact aspherical manifolds are homotopy equivalent, any homotopy equivalence is homotopic to a homeomorphism.  In particular, the two manifolds are in fact homeomorphic.
This conjecture has been the focus of much research in manifold topology and beyond after a reformulation in terms of $K$-theory and surgery, and it has had an impact on classification of manifolds more general than aspherical.  The current state of affairs can be learned from Weinberger's book \cite{sW:22}.  It is possible to prove the Borel conjecture for aspherical manifolds in dimensions 5 and higher with a given fundamental group by studying the assembly maps in both $K$-theory and $L$-theory associated to that fundamental group.  What we do in this paper is prove the isomorphism conjecture in $K$-theory for a very large class of fundamental groups specified by their metric properties.  The $L$-theory version of the conjecture can be proved using the same method.

The isomorphism conjectures were first explicitly stated by Farrell and Jones in \cite{tFlJ:93}.  Since the fundamental groups of compact aspherical manifolds have no torsion, and for the application to topological rigidity one is really interested in the integral conjectures, we formulate the simplest and most important version of the conjectures.

The Farrell-Jones conjectures are stated in terms of assembly maps such as $a_K$, which is the focus of this paper.  The definition of the $K$-theory assembly maps goes back to Loday.
The domain is the group homology spectrum $B\Gamma_{+} \wedge K (A)$ of a discrete group $\Gamma$ with coefficients in the nonconnective algebraic $K$-theory spectrum $K (A)$ of a ring $A$.
The stable homotopy groups of $K (A)$ coincide with
the lower $K$-theory of $A$ in negative dimensions and Quillen's
$K$-theory in nonnegative dimensions.

We use $\Gamma$ as the notation for a finitely generated group throughout the paper.  This is justified by the analogy with the traditional notation for discrete lattices in Lie groups but also by the fact we extensively use the letter G for the important \textit{G}-theory constructions.

Given a ring $A$, one can view an
element $\gamma$ of $\Gamma$ as an isomorphism of the free rank one
$A[\Gamma]$-module with the inverse $\gamma^{-1}$. Following
Loday, to each isomorphism $f$ of free finitely generated
$A$-modules there corresponds an $A[\Gamma]$-isomorphism $\gamma
\otimes f$ of finitely generated $A[\Gamma]$-modules. This ``assembly'' on the matrix level induces the assembly map 
\[
a_K (\Gamma,A) \colon B\Gamma_{+} \wedge K (A) \longrightarrow K (A[\Gamma]). 
\] 
If, in addition, the ring $A$ has an involution, there is an analogous map in $L$-theory,
\[
a_L (\Gamma,A) \colon B\Gamma_{+} \wedge L^{\langle-\infty \rangle} (A)
\longrightarrow L^{\langle-\infty \rangle} (A[\Gamma]).\]

Recall that a Noetherian ring $A$ is called \textit{regular} if every finitely generated $A$-module has a finite resolution by finitely generated projective $A$-modules.

{\bf The \textit{K}- and \textit{L}-theory Farrell-Jones isomorphism conjectures.}
For a regular Noetherian ring $A$ and a torsion-free group
$\Gamma$, the algebraic $K$-theory of the group ring $A[\Gamma]$ is isomorphic
to the homology of the group $\Gamma$ with
coefficients in the nonconnective $K$-theory of $A$.  If in addition the ring $A$ has an involution, there is an analogous isomorphism in algebraic $L$-theory.  More explicitly, the assembly maps $a_K (\Gamma,A)$ and
$a_L (\Gamma,A)$ are weak homotopy equivalences.

For surveys on the Farrell-Jones conjectures, the motivations, applications, and the recent advances, we refer to  \cite{aB:19,wL:22,hRmV:19}.

The most important regular Noetherian ring for geometric applications  is the ring of integers $\mathbb{Z}$.  The results of this paper will apply to all regular Noetherian rings of finite global dimension, including $\mathbb{Z}$.  

While closed aspherical manifolds are the main focus of the applications, the results of this paper will apply to all aspherical finite complexes where the fundamental group comes from a specific large class of groups.

Asymptotic dimension is a well-known global characteristic of a finitely generated group as a metric space with a word-length metric.
The groups that have \textit{finite asymptotic dimension} is a vast class, including Gromov's word hyperbolic groups, the mapping class groups, and other important families in geometric group theory.  

First, we state a concise version of the theorem.

{\bf The Main Theorem.}
{\it Let $\Gamma$ be a group with finite $K(\Gamma,1)$, and let $A$ be a regular Noetherian ring.  We endow $\Gamma$ with the word metric with respect to an arbitrary finite generating set.  
If the group $\Gamma$ with this word metric has finite asymptotic dimension, and the ring $A$ has finite global dimension, then the $K$-theoretic assembly map $a_K (\Gamma,A)$ is an equivalence.}

This theorem follows from a more technical expanded version which emphasizes the generality of the result. 

We gave a definition of $G$-theory of the group ring $A[\Gamma]$, as part of controlled $G$-theory $G (M, A)$ for an extended metric space $M$ or more generally a coarse space in the sense of Roe, reviewed later in the paper. 
There is an equivariant assembly map 
\[ 
A_G \colon h^{\mathit{lf}} (\Gamma; G(A)) \longrightarrow G (\Gamma, A).
\] 
The fixed points of $A_G$ give exactly the $G$-theory assembly map
\[
a_G \colon B\Gamma_{+} \wedge G (A) \longrightarrow G (A[\Gamma]). 
\] 
We can further consider the assembly map
\[ 
A_{G,D} \colon h^{\mathit{lf}} (\Gamma \times D; G(A)) \longrightarrow G (\Gamma \times D, A).
\] 
where $D$ is an arbitrary discrete set.

We obtain the Main Theorem from the following general $G$-theoretic theorem and its corollary.

\textbf{The \textit{G}-theory Isomorphism Theorem.}
{\em Suppose the assembly map $A_{G,D}$ is a non-equivariant equivalence for the given group $\Gamma$ with finite $K(\Gamma,1)$ and a regular Noetherian ring $A$, then $a_G$ is an equivalence.}

The assembly maps $a_G$ and $a_K$ are related via Cartan maps that form a commutative square
\[
\xymatrix{
B\Gamma_{+} \wedge G (A)
 \ar^-{a_G}[r]  
&G (A[\Gamma])  \\
B\Gamma_{+} \wedge K (A) \ar_{\mathrm{id} \wedge \kappa_A}[u]
 \ar^-{a_K}[r]
&K (A[\Gamma]) \ar_{\kappa_{A[\Gamma]}}[u]
}
\]
The Cartan map $\kappa_A \colon K(A) \to G(A)$ is an equivalence for all regular Noetherian rings.  We'll say $\Gamma$ is {\em coarsely coherent} if $\kappa_{A[\Gamma]}$ is an equivalence whenever the coefficient ring $A$ has finite global dimension.

\textbf{The \textit{K}-theory Isomorphism Theorem.}  
{\em We start with a regular Noetherian ring $A$ that has finite global dimension.  Suppose the assembly map $A_{G,D}$ is a non-equivariant equivalence for the given group $\Gamma$ and the ring $A$, and $\Gamma$ is coarsely coherent, then the $K$-theoretic assembly map $a_K (\Gamma,A)$
is an equivalence.  All of these assumptions are satisfied in the case the group $\Gamma$ has finite asymptotic dimension. }

To see that the $K$-theory theorem is a corollary to the $G$-theory theorem, we only need the fact that finite asymptotic dimensional groups are coarsely coherent.  That was the main result of \cite{gCbG:16}.

{\bf Comment:} This proof strategy allows to separate checking the $K$-theory Isomorphism Conjecture into two steps: one that checks the $G$-theory assembly $A_{G,D}$ is an equivalence for the given $\Gamma$ and one that checks $\Gamma$ is coarsely coherent.  For example, it's already known that groups $\Gamma$ with finite decomposition complexity are coarsely coherent from \cite{bG:13}, so the $K$-theory Isomorphism Conjecture for this class of groups would follow from the conjecture about $A_{G,D}$. 

We feel  that the two conditions together are a good candidate for an ``if and only if'' statement concerning equivalence of the $K$-theoretic assembly.  Both statements hold for all known examples, and our intuition suggests that failure of either one would force failure of the equivalence statement. 

{\bf Guide to  the proof:}
The method of proof of the split injectivity of the assembly map $a_K (\Gamma, A)$ in \cite{gC:95} and \cite{gCbG:04} is to observe that the assembly map is split by the inclusion $i_{\Gamma} \colon K(X,A)^{\Gamma} \hookrightarrow K(X,A)^{h \Gamma} $, where $K(X,A)$ denotes an equivariant version of the Pedersen-Weibel bounded $K$-theory (see \cite{ePcW:85}) of the universal cover $X$ of a manifold $K(\Gamma,1)$ space.  The method of proof that the assembly is an equivalence is to show that the inclusion $i_{\Gamma}$ is itself  split.  This will guarantee the result, since it assures that the existing splitting (which is surjective, since it is projection on a wedge factor) is also injective. 
So we set out to find an equivariant spectrum $M$ that fits in a commutative triangle
\[
\xymatrix{
 {G(X,A)^{\Gamma}}   \ar[dr]^-{s} \ar[rr]^-{i_{\Gamma}}
&&{G(X,A)^{h \Gamma}}  \ar[dl]_-{t} \\
& {M}^{\Gamma}
}
\]
where $s$ is an equivalence. It turns out, though, that there is a similar result (our Proposition \ref{SFP}) which suffices and doesn't require a  precise analysis which identifies the homotopy fixed point set.  We say a spectrum with $\Gamma$-action $E$ is {\em homotopy computable} if the natural map $E^{\Gamma} \rightarrow E^{h \Gamma}$ is an equivalence.  
Proposition \ref{SFP} asserts that for a $\Gamma$-spectrum $X$, if there is an equivariant map $j \colon X \rightarrow E$, where $E$ is homotopy computable, and so that there is a splitting of $j^{\Gamma}$, then the map $X^{\Gamma} \hookrightarrow X^{h\Gamma} $ also admits a right inverse, and so is a split injection. 

To carry this out, we need to find the right context in which to work.  The proof of the splitting of the assembly map in various situations which was carried out in \cite{gC:95}, \cite{gCeP:95}, and  \cite{gCbG:04} involved the use of {\em coarse geometry} (see  \cite{jR:03}).  Coarse geometry studies non-compact metric spaces (as well as some more general objects) under notions of maps which are  not restricted to being continuous maps but instead work with asymptotic notions of continuity.    In particular,  no local notion of continuity is assumed, and instead they encode a notion of continuity at  infinity.  The maps are required to satisfy a notion of properness, and so the theory is well suited to the study of non-compact manifolds.  As an example, in the case of the universal cover $\tilde{X}$  of a finite $K(\Gamma , 1)$-space, the inclusion of  a single orbit $\Gamma \cdot x \hookrightarrow  \tilde{X}$  is a coarse equivalence. The Pedersen-Weibel bounded $K$-theory $K(M)$  of a proper metric space $M$  (see \cite{ePcW:85}) is an example of a spectrum valued functor on the category of proper metric spaces and coarse maps. It was the key ingredient in the papers \cite{gC:95}, \cite{gCeP:95}, and \cite{gCbG:04}. A key property of the construction (made equivariant) was that $K(\tilde{X},A)^{\Gamma}  \simeq K(A[\Gamma])$ under the action of the group $\Gamma$ by deck transformations, and the splitting was realized by a passage to homotopy fixed point sets.     Our goal is to produce a splitting of the inclusion $K(A[\Gamma]) \hookrightarrow K(\tilde{X},A)^{h \Gamma}$, and we will do this in a way which is entirely parallel to the method of proof of the injectivity of the assembly map.  

Our strategy, much like the strategy adopted in \cite{gC:95}, \cite{gCeP:95}, and \cite{gCbG:04}  is the exploitation of  Pedersen-Weibel  bounded $K$-theory, or rather of a generalization of that theory.  The original proof of injectivity is best cast in the category of {\em extended metric spaces}, i.e. metric spaces which permit the value $+ \infty$ for the metric, and proper coarse maps. We'll denote this category  by $\mathfrak{E}^m$. Where the proof of injectivity  worked in the  ``absolute" category of extend metric spaces, in our case we will  instead need to work in a relative category, the category of objects of $\mathfrak{E}^m$  with $\Gamma$-action equipped with equivariant  coarse reference maps to   $\Gamma$, where $\Gamma$ is  regarded as a metric space  using a word length metric derived from a finite generating set for $\Gamma$, and is equipped with the left multiplication action.    The required homotopy computable $\Gamma$-spectrum will be modeled by a version of bounded $K$-theory applied to pair of simplicial objects  over $\Gamma$,  the inclusion from the fixed point set into the homotopy fixed set is naturally modeled using this pair. The splitting of this  inclusion is then achieved by using three  aspects of the coarse category.  The first is the existence of a cone construction $C(X)$ in which local properties of spaces are expanded to be reflected at infinity, and therefore in versions of bounded $K$-theory. This construction is a modification of the cone construction of Pedersen-Weibel \cite{ePcW:89}, which is shown there to have the desired effect of pushing local topological behavior to infinity, and making it detectable by their bounded $K$-theory.      The second  is a construction which assigns to each metric space $X$  with left $\Gamma$ action a new metric space $X^{bdd}$, also with $\Gamma$-action,  in which the $\Gamma$-action is bounded, i.e for every $\gamma\in \Gamma$, the function $x \rightarrow d(x, \gamma x)$ is a bounded function on $X$.  This construction has the role of forcing the action "at infinity" to be trivial, and is a kind of homotopy orbit space construction for the category of actions of groups on coarse spaces.   The third aspect is the existence of suitable excision properties for a $G$-theoretic version of the Pedersen-Weibel construction, which is what finally permits us to construct the splitting.  

The technical ingredients required to carry out this program are fairly extensive, and we summarize them here.  
\begin{enumerate}
    \item{{\bf $G$-theory of group rings:} The general statement that the $K$-theory assembly map is an equivalence is known to be false, even for the group $\Gamma = \mathbb{Z}$, when the coefficient ring contains nilpotent elements. For this reason, our decision to work around this issue by passing to $G$-theory, where no such counterexamples are known, is a natural one.  One obstacle is that $G$-theory is classically defined only for Noetherian rings, so it is
    necessary to construct a version of $G$-theory which applies broadly to group rings of finitely generated groups.  This construction is carried out in \cite{gCbG:16}.   
     } 
    \item{{\bf Extension of the Pedersen-Weibel bounded $K$-theory to $G$-theory:} The proofs of split injectivity of the assembly in \cite{gC:95} and \cite{gCbG:04} utilize the Pedersen-Weibel $K$-theory of the universal cover of a $K(\Gamma, 1)$-manifold.  We require  an extension of this construction to the $G$-theoretic setting. This is carried out in \cite{gCbG:11}; the reader will find an exposition in section \ref{PWtheories} below. }
    \item{\label{equi}{\bf Equivariant versions of bounded $G$-theory:} Since the universal cover comes equipped with a group action of  $\Gamma$, we will require the construction of bounded $G$ theory to admit an equivariant version.  This is constructed in \cite{gCbG:19}. One essential property is that the fixed points of this construction give precisely the $G$-theory of group rings in (1). }
    \item{\label{fibred} {\bf Fibred versions of bounded $G$-theory:} We will need to develop theories whose ''cells" are copies of $G(A[\Gamma])$, and this will amount to a version of parametrized homotopy theory in the coarse category, where the base is $\Gamma$ regarded as a coarse space using its word length metric. }
\end{enumerate}

From items (\ref{equi}) and (\ref{fibred}) above, given an extended metric  space $X$ with an action of a finitely generated group $\Gamma$, we are able to construct a $\Gamma$-equivariant  spectrum $\mathbb{G}_{\Gamma}(X, A)$, where $X$ is a coarse space and $A$ is a commutative coefficient ring. It is to be regarded as a parametrized cohomology theory over the base $\Gamma$, applied to the space $\Gamma \times X \rightarrow \Gamma$ over $\Gamma$. To orient the reader, we  give some examples of the values this functor takes, as well as some properties it possesses.  

\begin{itemize}
    \item{$\mathbb{G}_{\Gamma} (pt,A) \cong \mathbb{G}_{\{ e \}}(\Gamma,A) $  and $\mathbb{G}_{\Gamma} (pt,A)^{\Gamma} \cong G(A[\Gamma]) $, where $G$ denotes the $G$-theory of group rings constructed in \cite{gCbG:16}. 
    $G_{\{ e \}}(\Gamma, A) $ is equivalent to the bounded $G$-theory of the universal covering of a finite classifying space of $\Gamma$, which is often proved to be equivalent to the locally finite (Borel-Moore) homology of that universal covering.   }
    \item{For any finite set $X$, regarded as an extended metric space where the distances between any two distinct points are $+ \infty$ and equipped with the trivial $\Gamma$-action, we have $\mathbb{G}_{\Gamma}(X,A)^{\Gamma}  \cong h(X_+, G(A)) $.  We also have $\mathbb{G}_{\Gamma}(X,A) \cong h(X_+, G_{\{ e \}}(\Gamma, A)) $. For infinite sets $X$ with the trivial $\Gamma$-action, we similarly find that $\mathbb{G}_{\Gamma}(X,A)^{\Gamma}  \cong h^{\mathit{lf}} (X, G(A)) $, where $h^{lf}$ is the locally finite homology used in \cite{gC:95}.}
    
    \item{For any finite simplicial set  $X_{\dotr}$ with trivial $\Gamma$-action,  $\mathbb{G}_{\Gamma} (C(X)_{\dotr}, A) \cong \Sigma X_{\dotr} \wedge G(A) $ and $\mathbb{G}_{\Gamma} (C(X)_{\dotr}, A)^{\Gamma}  \cong \Sigma X_{\dotr} \wedge G(A[\Gamma])$ where the cone $C(X)_{\dotr} $ is regarded as an object in $\mathfrak{E}^m$. }
    
   \item{The construction $\mathbb{G}_{\Gamma} (X,A) $ actually admits an action by $\Gamma \times \Gamma$, whose restriction to the diagonal $\Gamma \subseteq \Gamma \times \Gamma$ is the action we use unless otherwise indicated.  }
   
    \item{For any finite complex $X$ with $\pi _1(X,x_0) \cong \Gamma$, there is an equivalence 
    $$ \mathbb{G}_{\Gamma}(C(X),A) \cong \mathbb{G}_{\Gamma} (C(\tilde{X}), A) ^{\{ e \} \times\Gamma}$$ where the $\Gamma$ action in question is the one coming from the action of $\Gamma$ on $\tilde{X}$ by deck transformations, that is the action of the group $\{ e \} \times \Gamma \subseteq \Gamma \times \Gamma$. The inclusion $\mathbb{G}_{\Gamma}(C(X),A) \hookrightarrow \mathbb{G}_{\Gamma} (C(\tilde{X}),A) $ of spectra with $\Gamma$-action will be referred to as the ``parametrized asymptotic transfer", and is introduced and discussed in section \ref{transfer}. 
    }
 \item{For finite coverings of coarse spaces with bounded (see above for the definition of this concept) $\Gamma$-actions, there are excision properties which permit the computation of the equivariant theory $\mathbb{G}_{\Gamma} $ in a local to global manner.  This package of theorems is contained in section \ref{ExcThms}.  Some of them are direct analogues of bounded excision in $K$-theory.  However, the crucial excision result, Theorem \ref{ExcGalltwo}, has no analogue in $K$-theory.  This is explained in Appendix \ref{karoubi} and serves to explain the need for coherence assumptions in the $K$-theory isomorphism theorem.}
   \end{itemize}
Our goal is to produce a splitting of an  inclusion $\mathbb{G}_{\Gamma} (pt, A) ^{\Gamma}  \hookrightarrow E^{h \Gamma} $ for a homotopy computable spectrum $E$, and the construction proceeds through several steps. 
\begin{enumerate}
    \item{Suppose we have a triangulation of a finite $X = B \Gamma$, and that it is embedded as a subcomplex in the interior of a triangulation of an $n$-disc $D^n$. After subdivision, we can convert the simplicial complex  into a simplicial set, which we denote by $D^n _{\dotr}$, and equip with the trivial $\Gamma$-action.   The starting point of our section will now be 
    $\mathbb{G}_{\Gamma} (C(D^n)_{\dotr}, C( \partial D^n)_{\dotr}, A)$, whose fixed point spectrum is $\Sigma h(D^n_{\dotr}, \partial  D^n_{\dotr}, G(A[\Gamma]))$. This is a suitable choice for the domain of our splitting, since it is an $(n+1)$-fold suspension of $G(A[\Gamma])$.  }
    
    \item{The next stage is the construction of the subcomplex $W$ of $D^n$, which is restriction of the triangulation of $D^n$ to the closure of the complement of $X$. We denote the corresponding simplicial set by $W_{\dotr}$.  This complex is given the trivial $\Gamma$-action. We construct the spectrum $\mathbb{G}_{\Gamma}(C(D^n)_{\dotr}, C(W)_{\dotr}, A) $. There is an induced map 
    $$ \mathbb{G}_{\Gamma} (C(D^n)_{\dotr}, C( \partial D^n)_{\dotr}, A) \rightarrow \mathbb{G}_{\Gamma}(C(D^n)_{\dotr}, C(W)_{\dotr}, A)
    $$ as well as an excision equivalence
    $$ \mathbb{G}_{\Gamma}(C(D^n)_{\dotr}, C(W)_{\dotr}, A) \overset{\simeq}{\longrightarrow}  \mathbb{G}_{\Gamma}(C(X)_{\dotr}, C(\partial X)_{\dotr},A ).
    $$
    It is an equivariant equivalence of spectra.  }
    
    \item{Form the universal covering $\widetilde{X}$ of $X$, and the pullback of $\widetilde{X}$ along the inclusion $\partial X \hookrightarrow X$, which we call $\partial \widetilde{X}$. Both of these constructions have evident  simplicial set versions, and we apply the parametrized asymptotic transfer (Definition \ref{HCGCXZ}) to them to obtain a map 
   $$  \mathbb{G}_{\Gamma}(C(X)_{\dotr}, C(\partial X)_{\dotr},A ) \overset{\tau}{\rightarrow}
    \mathbb{G}_{\Gamma}
    (C(\widetilde{X}))_{\dotr}, C(\partial  \widetilde{X})_{\dotr} ,A).$$ The realization of the simplicial $\Gamma$-spectrum $\mathbb{G}_{\Gamma}
    (C(\widetilde{X}))_{\dotr}, C(\partial  \widetilde{X})_{\dotr} ,A)$ will be homotopy computable. 
    }
    \item{Next, we map to metric coarse theory and apply the construction $(-)^{bdd} $, and obtain a map 
     $$  
    \mathbb{G}_{\Gamma}(C(\widetilde{X}))_{\dotr}, C(\partial  \widetilde{X})_{\dotr} ,A) \rightarrow  \mathbb{G}_{\Gamma}(C(\widetilde{X})^{bdd}, C(\partial \widetilde{X})^{bdd} ,A ) 
    $$
    All of these are equivariant maps.  After applying the last construction, the action on $C(\widetilde{X})$ becomes bounded.  We describe an equivalence of the last spectrum to the one associated to the trivial action in \cite{gCbG:19}.  So now we have 
     $$  
    \mathbb{G}_{\Gamma}(C(\widetilde{X}))_{\dotr}, C(\partial  \widetilde{X})_{\dotr} ,A) \rightarrow  \mathbb{G}_{\Gamma}(C(\widetilde{X})_0^{bdd}, C(\partial \widetilde{X})_0^{bdd} ,A ). 
    $$
    } 
    \item{Select a small metric ball $B$ in $D$ that lifts isometrically to a ball $\overline{B}$ in $\widetilde{X}$.  Just as we did before in (2) with the complement of $X$ in $D$, we can apply excision to the complement of $\overline{B}$ in $\widetilde{X}$.  We get an induced map 
    $$  
    \mathbb{G}_{\Gamma}(C(\widetilde{X})_0^{bdd}, C(\partial \widetilde{X})_0^{bdd} ,A ) \rightarrow  \mathbb{G}_{\Gamma}(C(\overline{B}), C(\partial \overline{B}), A). 
    $$
    } 
\end{enumerate}

It is this last spectrum that we are using as $M$ in the strategy of proving that $i_{\Gamma} \colon \mathbb{G}_{\Gamma} (C(D^n)_{\dotr}, C( \partial D^n)_{\dotr}, A) \simeq \Sigma^{n+1}  \mathbb{G}_{\Gamma} (pt,A)^{\Gamma} \hookrightarrow \Sigma^{n+1}  \mathbb{G}_{\Gamma} (pt,A)^{h \Gamma} $ is split injective.  We compute in section \ref{proofend} that on fixed points the induced map
$$
s \colon \mathbb{G}_{\Gamma} (C(D^n)_{\dotr}, C( \partial D^n)_{\dotr}, A)^{\Gamma} \to \mathbb{G}_{\Gamma}(C(\overline{B}), C(\partial \overline{B}), A)^{\Gamma}
$$
is an equivalence. This gives the required result that our inclusion is split.  

\section{Coarse Spaces}

The key ingredient in our proofs of the split injectivity of the assembly map for algebraic $K$-theory (the $K$-theory Novikov conjecture  formulation) has been the use of so-called coarse homology theories on the category of metric spaces.  Coarse homology theories are theories that are insensitive to phenomena occurring in bounded parts of the metric, but instead depend only on the structure of the ``space at infinity".  Our proof that the $K$-theoretic assembly map is in fact an equivalence will also depend on such theories, but it will require a more extensive development of their behavior  under homotopy theoretic constructions. We now define the objects of the category we will study.   
\subsection{Definitions}
\begin{Definition} \label{coarsedef} \cite{jR:03}
    A {\em coarse structure} on a set $X$ is a collection $\mathcal{E}$ of subsets of $X \times X$ satisfying the following properties. 
    \begin{enumerate}
        \item{The diagonal $X \subseteq X \times X$ is in $\mathcal{E}$.}
        \item{If $U \in \mathcal{E}$, and $V \subseteq U$, then $V \in \mathcal{E}$.}
        \item{$\mathcal{E}$ is closed under finite unions. }
        \item{ If $U \in \mathcal{E}$, then so is $U^{-1} = \{(x,x^{\prime})| (x^{\prime},x) \in \mathcal{E}$\}}
        \item{If $U$ and $V$ are in $\mathcal{E}$, then so is $U \circ V = \{(x,x^{\prime}) | \exists \xi \in X $ so that $(x,\xi)$ and $(\xi, x^{\prime})$ are both in $\mathcal{E}$\}.}
    \end{enumerate}
    A {\em coarse space} is a set $X$ equipped with a coarse structure $\mathcal{E}_X$.  The elements of $\mathcal{E}_X $ are referred to as the {\em controlled sets}, and subsets  $B \subseteq X$ for which $B \times B$ is controlled are called the {\em bounded sets}. If $X$ is a coarse space, $Y \subseteq X$, and $E \subseteq X \times X$ is a controlled subset, then by the {\em $E$-neighborhood of $Y$}  we will mean the set $\{ x \in X|\exists y \in Y \mbox{ with } (x,y) \in E \}$. We will denote it by $E[Y]$. Proposition 2.19(a)  of \cite{jR:03} asserts that if $Y$ is a bounded set, then so is $E[Y]$.  
\end{Definition}

 \begin{Example}
     Let $(X,d) $ be a metric space.  Then the metric $d$ defines a coarse structure $\mathcal{E}(d)$ on $X$, which consists of the sets $U \subseteq X \times X$ on which the function $d$ on $X \times X$ is bounded. 
 \end{Example}
 
 \begin{Definition}
     Given two coarse spaces $(X, \mathcal{E}_X)$ and $(Y, \mathcal{E}_Y)$, a map $f:X \rightarrow Y$ is a \textit{coarse map} if 
    \begin{enumerate}
        \item\label{bounded}{whenever $B$ is a bounded subset of $Y$, $f^{-1}(B) $ is bounded in $X$.}
        \item\label{continuous}{for each controlled set $E$ in $X \times X$, the set $(f \times f)(E) $ is a controlled set in $Y$.  }
    \end{enumerate}
    {\em (\ref{bounded})} is a properness condition, and {\em (\ref{continuous})} is a coarse version of continuity. We denote the category of coarse spaces and coarse maps by $\mathfrak{E}$.  
    We say a coarse map $f \colon X \rightarrow Y$ is a {\em coarse equivalence} if there is a coarse map $g \colon Y \rightarrow X$ so that the graphs of $g \circ f$ and $f \circ g$ are controlled sets in $X$ and $Y$ respectively. 
 \end{Definition}
  \begin{Example}If $f \colon X \rightarrow Y$ is a proper continuous map of proper metric spaces (i.e. closed balls are compact), then $f$ is a coarse map of the corresponding coarse spaces.  On the other hand, there are many coarse maps from $X$ to $Y$ which are not continuous.

 \end{Example} 
 \begin{Remark}
     {\em The coarse category defined by J. Roe in \cite{jR:03} has morphisms which are equivalence classes of coarse maps, not the coarse maps themselves.  We need the maps themselves. }
 \end{Remark}

We will need to extend the definition of the term metric to include the notion of infinite distances.  

\begin{Definition}\label{MetSpsss}
An \textit{extended metric space} is a set $X$ and a function
$d \colon X \times X \to [0,\infty) \cup \{\infty\}$
which is reflexive, symmetric, and satisfies the triangle inequality
in the obvious way.
An extended metric space is \textit{proper} if it is a countable disjoint union of
metric spaces $X_i$ where $\im (d \vert X_i \times X_i) \subset [0,\infty)$,
and all closed metric balls in $X$ are compact.
The metric topology on an extended metric space can be
defined as usual.  An extended metric space is a coarse space in an obvious way. 
\end{Definition}

\begin{Proposition}\label{Bornol2}
If $X$ and $Y$ are proper metric spaces with metric functions $d_X$ and $d_Y$,
a map $f \colon X \to Y$ is \textit{uniformly expansive} if there is a real positive function $l$
such that
$d_X (x_1, x_2) \le r$ implies
$d_Y (f(x_1), f(x_2)) \le l(r)$.
 It should be emphasized that the map is not assumed to be continuous.  
The map $f$ is \textit{proper} if $f^{-1} (S)$ is a bounded subset of $X$ for
each bounded subset $S$ of $Y$.
A map $f$ from $X$ to $Y$  is a \textit{coarse map} precisely if it is proper and uniformly expansive.  
Similarly,  $f$ is a \textit{coarse equivalence} if there is a coarse map $g \colon Y \to X$ such
that $f \circ g$ and $g \circ f$ are bounded maps.
\end{Proposition}

We will  denote the full subcategory of $\mathfrak{E}$ on the extended metric spaces by $\mathfrak{E}^m$.  

There is a functor $\iota \colon \underline{Sets}^{proper} \rightarrow \mathfrak{E}^m$, where $\underline{Sets}^{proper}$ denotes the category of sets and proper (i.e. every finite subset has a finite inverse image) maps of sets, given by assigning to a set $X$ the extended metric space whose underlying set of points is $X$ and where any two distinct points have distance $+ \infty$ from each other. We will refer to this metric on $X$ the \textit{discrete metric}.

There is another simple method for constructing extended metric spaces, which we will find useful.

\begin{Definition}\label{LenMet}
	For any set $Z$, any subset $S \subset Z \times Z$, and any function $\phi \colon S \to [0, \infty]$ we denote by $d_{(S,\phi)}$ the infimum of all extended metrics on $Z$ with the property that $d_{(S,\phi)} (x,x') \ge \phi(x,x')$ whenever $(x,x')$ is in $S$.  It is easy to see that this defines a pseudo-metric, i.e. a function $d \colon Z \times Z \rightarrow [0, + \infty]$ satisfying all the defining properties except that it might take the value zero on some pair of distinct points.  However, as long as there is an actual extended metric on $Z$ which is bounded above by $\phi$ on $S$ then we know the pseudo-metric is a metric.
\end{Definition}

\subsection{The cone construction}\label{conedef}

As we have mentioned, the category $\mathfrak{E}^m $ has the property that it ignores bounded parts of the an extended metric space and captures behavior at infinity. Nevertheless, we will need to inject constructions from ordinary homotopy theory into the homology theories we will construct on $\mathfrak{E}^m$. In order to do this, we will need the construction of an infinite cone on a metric space  whose behavior at infinity exactly replicates that of the original space, but with a dimension shift of one.  A different version of our construction was studied in the paper \cite{ePcW:89}, where it was shown that the Pedersen-Weibel construction applied to an infinite Euclidean cone on a finite complex $X$ produced the suspension of the homology of $X$ with coefficients in the algebraic $K$-theory spectrum of a commutative ring $A$.  

\begin{Definition}\label{TX}
 	Given any extended metric space $(X, d)$ and a real function $f$ with $f(r) = 1$ for all $r \le 0$ and $f(r) \ge 1$ for all $r \ge 0$, by the \textit{cone on $X$ associated with $f$} we mean the set $C_f X = \mathbb{R} \times X$ equipped with the pseudometric $d_{(S,\phi)}$ from Definition \ref{LenMet}, where
 \[
 S = \{ ((t,x), (t,x')) \mid t \in  \mathbb{R}, x, x'  \in  X \} \cup 
 \{ ((t,x), (t',x)) \mid t,t' \in  \mathbb{R}, x \in  X \}
 \]
 and where $\phi ((t,x),(t,x')) = f(t) d(x,x')$ and $\phi ((t,x), (t',x)) = \| t' - t \|$. 
This defines an extended metric because the product extended metric on $\mathbb{R} \times X$ is dominated by $\phi$ on $S$.
\end{Definition}
\vspace{-.5cm}
\begin{center}
    \includegraphics[height=4cm]{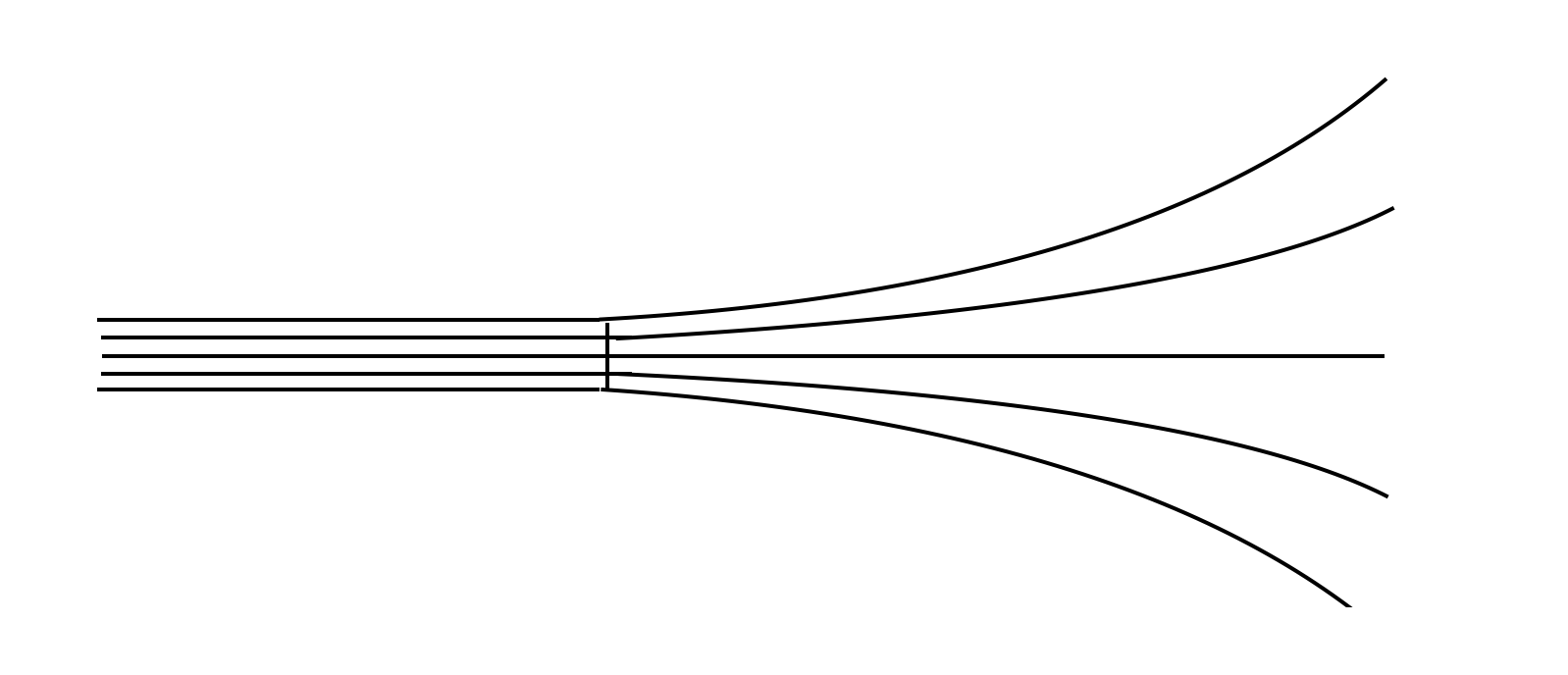}
\end{center}
 
\begin{Notation}\label{NotTX}
In the case $f(r) = 2^r$ for $r \ge 0$, we denote $C_f X$ by simply $CX$.
\end{Notation}

\subsection{Equivariance and the $(-)^{bdd}$ construction}

We will be working with finitely generated groups regarded as metric spaces.  

\begin{Definition}\label{bddefinition}
    For any extended metric space $(X,d)$ and group $\Gamma$, an {\em isometric action} of $\Gamma$ on $(X,d)$ is an action of $\Gamma$ on the underlying set $X$ so that every group element $\gamma$ acts by an isometry of $X$.  Similarly, a {\em coarse action} is an action of $\Gamma$ on $X$ so that every group element acts by a coarse map (and therefore by a coarse equivalence). An action by $\Gamma$ is {\em bounded} if for every $\gamma \in \Gamma$, the function $d(x, \gamma x)$ is bounded above by a (finite) real number $B(\gamma)$ which may depend on $\gamma$.  Of course, bounded actions are necessarily coarse due to the triangle inequality.  
\end{Definition}

In a rough sense, the metric spaces with bounded $\Gamma$-action are regarded as homotopy trivial, because they are to be thought of as close to the action by the identity map. Attached to any extended metric space with isometric action, we would like to construct an analogue to the homotopy orbit space construction, which canonically produces from an isometric $\Gamma$-action on an extended metric space $X$ a new extended metric $d^{bdd}$ on $X$ so that the action of $\Gamma$ on $X$ is bounded in the metric $d^{bdd}$.  In a sense this trivializes the action, and it is exactly this kind of trivialization we will require at a crucial step in our proof.  

Let $Z$ be any extended metric space with a free left $\Gamma$-action by isometries.
We assume that the action is properly discontinuous, that is, for any pair of points $z$ and $z'$,
the infimum over $\gamma \in \Gamma$ of the distances $d(z, \gamma z')$ is attained.

\begin{Definition}\label{OrbitMet}
We define the orbit space metric on $\Gamma \backslash Z$ by
\[
d_{\Gamma \backslash Z} ([z], [z']) = \inf_{\gamma \in \Gamma} d(z, \gamma z').
\]
\end{Definition}

\begin{Lemma}\label{JUQASW}
$d_{\Gamma \backslash Z}$ is a metric on $\Gamma \backslash Z$.
\end{Lemma}

\begin{proof}
The triangle inequality follows from the triangle inequality for $d$.
Symmetry follows from $d (z, \gamma z') = d(\gamma^{-1} z, z') = d(z', \gamma^{-1} z)$.
Finally, $d_{\Gamma \backslash Z} ([z], [z']) = 0$ gives $d (z, \gamma z') = 0$ for some $\gamma \in \Gamma$, so $z = \gamma z'$, and so $[z] = [z']$.
\end{proof}

\begin{Definition}\label{WORD}
The \textit{word-length metric} $d = d_{\Omega}$ on a group
$\Gamma$ with a fixed finite symmetric generating set $\Omega = \Omega^{-1}$ is the metric
induced from the condition that $d (\gamma, \gamma \omega) =1$, whenever $\gamma \in \Gamma$ and $\omega \in \Omega$.  
It is well-known that the pseudometric from Definition \ref{LenMet} makes $\Gamma$ a proper metric space with a free isometric action by $\Gamma$ via left multiplication, as well as a coarse action by multiplication on the right.  
\end{Definition}

Now suppose $X$ is an extended metric space with a left $\Gamma$-action by isometries.

\begin{Definition}\label{Xbdd}
Define
\[
X^{bdd} = X \times_{\Gamma} \Gamma
\]
where
the right-hand copy of $\Gamma$ denotes $\Gamma$ regarded as a metric space with the word-length metric associated to a finite generating set as in Definition \ref{WORD},
the group $\Gamma$ acts by isometries on the metric space $\Gamma$ via left multiplication.
The product $X \times \Gamma$ is given the max metric.
Then $X \times_{\Gamma} \Gamma$ denotes the orbit metric space associated to the diagonal left $\Gamma$-action
on $X \times \Gamma$.
We will denote the orbit metric by $d^{bdd}$.
\end{Definition}

A natural left action of $\Gamma$ on $X^{bdd}$ is given by the formula
$\gamma [x, \gamma'] = [x, \gamma' \gamma^{-1}]$.  Then $\gamma [x,e] = [x, \gamma^{-1}] = [\gamma x, e]$.

Suppose $\Gamma$ is a finitely generated group with a chosen finite generating set which determines the word metric $d_{\Gamma}$ and the corresponding norm $\vert \gamma \vert = d_{\Gamma} (e, \gamma)$.

\begin{Lemma}\label{JUQASW2}
If the action of $\Gamma$ on $X$ is bounded, with the bound function denoted $B(\gamma)$, then there is a real valued function $B_{\ast} \colon [0, \infty) \to [0, \infty)$ such that $\vert \gamma \vert \le s$ implies $B({\gamma}) \le B_{\ast} (s)$. 
\end{Lemma}

\begin{proof}
One simply takes $B_{\ast} (s) = \max \{ B({\gamma}) \ \mathrm{such} \ \mathrm{that} \  \vert \gamma \vert \le s \}$.
\end{proof}

\begin{Proposition}\label{ECEjkiu}
The natural action of $\Gamma$ on $X^{bdd}$ is bounded.
\end{Proposition}

\begin{proof}
We choose $B_{\gamma} = \vert \gamma \vert$.
Now
\begin{align}
d^{bdd} ([x,e], [\gamma x,e]) &= \inf_{\gamma' \in \Gamma} d^{\times} ((x,e), \gamma' (\gamma x,e)) \notag \\
&\le d^{\times} ((x,e), \gamma^{-1} (\gamma x,e)) \notag \\
&= d^{\times} ((x,e), (x,\gamma^{-1})) = d_{\Gamma} (e,\gamma^{-1}) = \vert \gamma^{-1} \vert = \vert \gamma \vert, \notag
\end{align}
where $d^{\times}$ stands for the max metric on the product $X \times \Gamma$.
\end{proof}

\begin{Definition}\label {leftbdd44}
Let $b \colon X \to X^{bdd}$ be the natural bijection given by $b(x) = [x,e]$ in the orbit space $X \times_{\Gamma} \Gamma$.
\end{Definition}

\begin{Proposition}\label{HJDSEO}
The map $b$ is a coarse map.
\end{Proposition}

\begin{proof}
Suppose $d^{bdd} ([x_1,e], [x_2,e]) \le D$, then $d^{\times} ((x_1,e), (\gamma x_2,\gamma)) \le D$ for some $\gamma \in \Gamma$,
so $d(x_1, \gamma x_2) \le D$ and $\vert \gamma \vert \le D$.
Since the left action of $\Gamma$ on $X^{bdd}$ is bounded, there is a function $B_{*}$ guaranteed by Lemma \ref{JUQASW2}.
Now
\[
d (x_1, x_2) \le d(x_1, \gamma x_2) + d (x_2, \gamma x_2)
\le D + B_{\ast} (D).
\]
This verifies that $b$ is proper.  It is clearly distance reducing, therefore uniformly expansive with $l(r)=r$.
\end{proof}

If one thinks of $X^{bdd}$ as the set $X$ with the metric induced from the bijection $b$, the map $b$ becomes the coarse identity map between the metric space $X$ with a left action of $\Gamma$ and the metric space $X^{bdd}$ where the \textit{action is made left-bounded} but is no longer by isometries.  

\begin{Proposition}\label{Gammaright}
Let $X$ be the group $\Gamma$ with the discrete metric.  Then $X^{bdd}$ is a metric space isometric to $\Gamma$ with a word metric.  The natural left action of $\gamma$ on $X^{bdd}$ corresponds to the left-bounded right multiplication action by $\gamma^{-1}$ on $\Gamma$. 
\end{Proposition}

\begin{proof}
	The orbits in $X \times \Gamma$ are distance $1$ apart precisely when there are two representatives $(e,\xi)$ and $(e,\xi')$ with $\xi =  \xi' \omega$ for a generator $\omega$ from a finite symmetric generating set $\Omega$.  So the orbit metric $d([e,\xi], [e,\xi'])$ is given by the length of the shortest word $w$ with $\xi = \xi' w$.  
\end{proof}

We now have a construction we can apply in any instance of a free properly discontinuous action by $\Gamma$ on a metric space $X$.  In some cases $X$ has a very canonical set of metrics, for example in the case of a fundamental group acting cocompactly on the universal cover $X$.  On the other hand, if $\Gamma$ is the fundamental group of a manifold embedded in a Euclidean space, the normal bundle doesn't have such a set of metrics.  We are particularly interested in using left-bounded metrics on regular neighborhoods related to the normal bundles.  In these cases the map induced from $b$ preserves the most relevant $K$-theoretic information.

\subsection{Simplicial objects and geometric realization} \label{simpreal}
We have mentioned the functor $\iota $ which embeds the category of sets and proper set maps as a subcategory of $\mathfrak{E}^m$.  We will refer to the image of $\iota$ as the subcategory of {\em discrete extended metric spaces}, and denote it by $\mathfrak{E}^m_{\delta}$. The functor  $\iota$ extends to a functor from the category of simplicial objects in $\underline{Sets}^{proper}$ to the category $s_{\dotr}\mathfrak{E}^m$ of simplicial objects in $\mathfrak{E}^m$.  For any proper simplicial set $X_{\dotr}$ (i.e. simplicial object in the category $\underline{Sets}^{proper}$), we define an extended metric on the geometric realization $|X_{\dotr}|$ as follows.  For each $n \geq 0$, we assign a metric $d^{{\Delta}[n]}$ on the standard simplex $\Delta [n]$ via its standard inclusion in $\mathbb{R}^{n+1}$, normalized so that the diameter of $\Delta [n]$ is 1.    This assignment then determines an extended metric $d$ on the set $\coprod _n X_n \times \Delta [n]$ by (a) $d((x,w), (x^{\prime}, w^{\prime}) )= + \infty$ if $x \neq x^{\prime}$ and 
(b) $d((x, w), (x, w^{\prime}) = d^{\Delta[n]}(w,w^{\prime})$.  The points of the standard geometric realization are now constructed as equivalence classes under an equivalence relation on $\coprod _n X_n \times \Delta [n]$.  
It is well known that, in the case of a locally finite simplicial complex $X$, the infimum of sums of lengths of disjoint euclidean subpaths in each simplex taken over all paths between points in $|X_{\dotr}|$ is a well-defined proper metric.

\section{Homology Theories on $\mathfrak{E}$}
The key topological objects we will need for our proof can be considered homology theories defined on the category $\mathfrak{E}^m. $  The first construction of this type is the {\em bounded $K$-theory} of E.K.~Pedersen and C.~Weibel \cite{ePcW:85}, \cite{ePcW:89}.  It is a theory which takes as input a proper metric space and a commutative ring, and from it produces a spectrum. It has some of our desired properties, in particular that it is coarse in the sense that it is insensitive to bounded structure in the proper metric space and instead reflects its behavior at infinity.  One problem with it from our point of view  is that it does not take into account non-split exact sequences.  For this reason, we have been forced to develop a $G$-theory version, for which a relatively involved notion of exactness is defined. The $G$-theory also has the coarseness property, and is  the target for a comparison map from the Pedersen-Weibel $K$-theory spectrum which can under certain circumstances be proved to be an equivalence.  These theories are  discussed in \cite{gCbG:11}, \cite{gCbG:18}, and \cite{gCbG:19}. In this section, we introduce these theories, including their fibred and  equivariant versions, as well as the excision properties required for our main result.  We will first develop the absolute $G$-theory and follow that first with the fibred theory and then the equivariant theory.  The key theory for the final result will be the equivariant fibred theory.

\subsection{Bounded categories}\label{PWtheories} Let  $X$  be a coarse space, and let  $A$ be a commutative ring. We form   a category $\mathcal{B}(X,A)$, whose objects are triples $(F,B,\varphi)$ where $F$ is a free $A$-module (not necessarily finitely generated) with a basis $B$, equipped with a labeling function $\varphi \colon B \rightarrow X$, which has the  property that for any bounded $U \subseteq X$, $\varphi ^{-1}(U)$ is finite.  A morphism from $(F_1, B_1, \varphi _1)$ to $(F_2, B_2, \varphi _2)$ is an $A$-linear map $f \colon F_1 \rightarrow F_2$ with the property that there is a controlled set $E \subset X \times X$ so that for every $\beta \in B_1$, $f(\beta) $ is contained in the span of $\varphi _2^{-1}(E[\varphi _1(\beta)])$, where $E[-]$ denotes the $E$-neighborhood of $ \varphi _1 (\beta)$ as defined in Definition \ref{coarsedef}. With this definition of morphism, we have defined a category, due to the fact that $E[\varphi _1 (\beta)]$ is a bounded set, as pointed out in Definition \ref{coarsedef}.   

\begin{Remark}
In the case of an extended metric space, this translates exactly to the category defined by Pedersen and Weibel in \cite{ePcW:85} and \cite{ePcW:89}. We have introduced the additional generality so as to include the case of fibred $K$-theory \cite{gCbG:11}, which we will require.  Most proofs using the additional generality translate verbatim, and we will alert the reader when they do not.
\end{Remark}

It is known from Cardenas--Pedersen \cite{mCeP:97} that the category of bounded chain complexes in  $\mathcal{B}(X,A)$ is a Waldhausen category \cite{fW:83}.  Using the $S_{\dotr}$ construction from \cite{fW:83}, one then constructs a spectrum from it.  As explained in \cite{mCeP:97,gCbG:11}, one can use products with the metric space $\mathbb{R}$ to construct a potentially non-connective spectrum, which we will denote by $K(X,A)$, and prove that it satisfies certain excision properties. The bounded $K$-theory construction  is sufficient to prove that the assembly map is a split injection for all coefficient rings  in many geometric situations, see for example \cite{gC:95}, \cite{gCeP:95}, and  \cite{gCbG:04}.   However, it is clear that the bounded $K$-theory alone will not be enough to prove that the assembly map is an equivalence due to the presence of counterexamples, even for the group $\Gamma = \mathbb{Z}$, when the coefficient ring contains nilpotents \cite{cW:79}. On the other hand, these counterexamples are not known to occur for  $G$-theory, i.e. the $K$-theory of the category of finitely presented modules rather than projective ones. This observation suggests desire for the construction of a theory $G(X,A)$ which bears the same relationship to $K(X,A)$ as $G$-theory of ordinary rings does to $K$-theory. 

  To explain how this is carried out, we remark  that the objects in $\mathcal{B}(X,A) $ can be regarded as cosheaves (see \cite{yK:60}) $\Phi$  on the category of all subsets of $X$, with values in the category of $A$-modules, via the assignment $\varphi(U) = \mbox{span}(\varphi ^{-1}(U))$. 
  
  \begin{Definition}
      Let $X$ be a topological space, and let $\underline{C}$ denote a cocomplete category.  A $\underline{C}$-valued  precosheaf on $X$ is a covariant functor $\Phi$ from the category of open sets in $X$  and inclusions to $\underline{C}$. $\Phi$ is a {\em cosheaf} if for every open covering $\{U_i \}_{i \in I}$  of an open set $U$, the natural map 
\[
\truecolim{i \in I} {\Phi(U_i)} \longrightarrow \Phi (U)
\]
is an isomorphism in $\underline{C}$. 
  \end{Definition}
  \begin{Remark} Assuming $X$ is given the discrete topology, and
      letting $(F,B,\varphi)$ be an object of $\mathcal{B}(X,A)$, the associated cosheaf $\Phi _F$ is given by $\Phi _F(U)  = \mbox{span} (\varphi ^{-1}(U))$. Note that the values of $\Phi$ are finitely generated free $A$-modules. 
  \end{Remark}

  Bounded $G$-theories are built by forming Waldhausen categories of precosheaves satisfying certain restrictions.  We now introduce a number of these. 
  \begin{Definition} \label{splitleanins}
  Let $\Phi$ denote an $A$-module valued precosheaf on $X$, where  $X$ is a coarse space, equipped with the discrete topology.   We say $\Phi$ is an {\em $X$-filtration} of the $A$-module $M = \Phi(X)$ if it factors through the subcategory $\mathcal{C}_M$ of $A$-$mod$ whose objects are all submodules of $M$ and whose only morphisms are inclusions of one submodule in another. Let $\Phi$ denote an $X$-filtration.  
  \begin{enumerate} 

  \item\label{locfinite} { $\Phi$ is {\em locally finite} if $\Phi (U)$ is a finitely generated submodule of $\Phi(X)$ whenever $U$ is a bounded subset of $X$}.
  \item\label{lean}{$\Phi$ is $E$-{\em lean}, or just {\em lean} if $E$ is not relevant to the context, if for  a
 controlled set $E$ we have that
\[
\Phi \left( U \right) \subseteq \sum_{x \in U} \Phi \left( E[x]) \right)
\]
for any subset $U$ of $X$. }
  \item\label{insular}{ $\Phi$ is {$E$-\em insular}, or just {\em insular},  if for  a controlled set 
$E$ we have that 
\[
\Phi(U_1) \cap \Phi(U_2) \subseteq \Phi \big( E[U_1]  \cap E[U_2]  \big)
\]
for any pair of subsets $U_1$, $U_2$ of $X$. }
  \item\label{admissible}{ $\Phi$ is {$E$-\em admissible}, or just {\em admissible} if it is both $E$-lean and $E$-insular with respect to a single controlled set $E$. }
  \end{enumerate}

  \end{Definition}
  \begin{Remark}
      Conditions  {\em(\ref{lean})}  and {\em (\ref{insular})} are coarse or approximate versions of the cosheaf condition. They permit the cosheaf condition to hold up to a bounded difference, in an appropriate sense. 
  \end{Remark}

  Let $\Phi$ denote an $X$-filtration of $A$-modules.  For any submodule $M \subseteq \Phi(X)$, we define a new {\em induced} $X$-filtration $\Phi _M$, which is defined by $\Phi _M(U) = M \cap \Phi(U)$.  By a {\em subconstruction of $\Phi$} we will mean an $X$-filtration of the form $\Phi _{\Phi(U)}$ for some $U \subseteq X$, and we'll denote the subconstruction by $\Phi_U$.  We will also say that a subfiltration  $\Psi \subseteq \Phi$ is a {\em controlled enlargement} or simply an {\em enlargement} of a subconstruction $\Phi _U$ if $\Psi = \Phi_{\Psi (X)}$ and there is a controlled set $E$ so that 
  $\Phi_U \subseteq \Psi \subseteq \Phi_{E[U]}$. 
  
  We will require some conditions on an $X$-filtration which assert that certain properties hold for $\Phi$ and all subconstructions of $\Phi$.  

  \begin{Definition} \label{adapted}  Let $\Phi$ be an $X$-filtration.  
  \begin{enumerate}
      \item{We say $\Phi$ is {\em adapted} if it and all of its subconstructions are admissible. }

        \item{We say $\Phi$ is {\em coarsely adapted} if each of its subconstructions has an enlargement which is admissible.}

      \item{If $Y \subseteq X$ is a subset, we say $\Phi$ is {\em supported near $Y$} if there is a controlled set $E$ so that $\Phi (X) = \Phi(E[Y])$.}
  \end{enumerate}
      
  \end{Definition}

Next, we need to define a notion of control of the morphisms in the category which will define $G$-theory.  

\begin{Definition}\label{controldefinition} Let $\Phi$ and $\Psi$ be two $X$-filtrations, where $X$ is an extended metric space. Then a morphism from $\Phi$ to $\Psi$ is a homomorphism of $A$-modules $f \colon \Phi(X) \rightarrow \Psi (X)$
which has the property that there is a controlled set  $E$ so that for all $U \subseteq X$,  $f(\Phi (U)) \subseteq \Psi (E[U])$. This agrees with the definition in the Pedersen-Weibel construction when applied to the cosheaves defined there. We say a morphism $f \colon \Phi \rightarrow \Psi$ is {\em filtration zero} if $f(\Phi(U)) \subseteq \Psi(U)$ for all $U \subseteq X$. 
\end{Definition}

\begin{Definition}\label{CoAd}
We define the category $\mathcal{G}(X,A)$ to have as its objects $X$-filtrations which are locally finite and coarsely adapted, with the morphisms defined in Definition \ref{controldefinition} above. Let $\mathcal{G}(X,A) _{\leq Y}$ denote the full subcategory of $X$-filtrations supported near $Y$. It is easy to check that the assignments $X \mapsto \mathcal{G}(X,A)$ and $(X,Y) \mapsto \mathcal{G}(X,A) _{\leq Y}$ are functorial for coarse maps in the obvious sense.  
 \end{Definition}
 
We will define exact category structures on $\mathcal{G}(X,A)$ and $\mathcal{G}(X,A)_{\leq Y}$ in section \ref{spectrum}, and consequently  Waldhausen categories and spectra, which will be the  $G$-theory spectra of the coarse space $X$ with coefficients in the ring $A$, and the corresponding $G$-theory for the $X$-filtrations supported near $Y$.

We need two   facts about coarsely adapted objects in $\mathcal{G}(X, A)$ that will be used later in the paper. The first is that the admissibility property is preserved under isomorphisms of $X$-filtrations.  

\begin{Lemma}\label{Propone}
    If $f \colon \Phi_1 \to \Phi_2$ is an isomorphism of $X$-filtrations, and if $M$ is a subobject of $\Phi_1 (X)$ so that the induced filtration $M_1$ is admissible, then the induced filtration $f(M)_2$ of $f(M)$ in $\Phi_2$ is also admissible. 
\end{Lemma}
\begin{proof}
Suppose $M_1$ is $E$-admissible and let $K$ be a control set that controls both $f$ and its inverse.  We can see that $f(M)_2$ is $KEK$-lean from the sequence of inclusions for any given subset $U$:
\begin{multline*}
f(M)_2 (U) \subset f(M_1 (K[U])) \subset \sum_{x \in K[U]} f(M \cap \Phi_1 (E[x])) \subset \\
\sum_{x \in K[U]} f(M) \cap \Phi_2 (KE[x])) \subset 
\sum_{x \in U} f(M) \cap \Phi_2 (KEK[x])) = \sum_{x \in U} f(M)_2 (KEK[x]).
\end{multline*}
A similar estimate 
\begin{multline*}
f(M)_2 (U) \cap f(M)_2 (V) 
\subset 
f(M_1 (K[U])) \cap f(M_1 (K[V])) = \\
f(M \cap \Phi_1 (K[U]) \cap f(M \cap \Phi_1 (K[V]) = \\
f(M) \cap f\Phi_1 (K[U]) \cap f\Phi_1 (K[V]) = \\
f(M) \cap f(\Phi_1 (K[U]) \cap \Phi_1 (K[V])) \mathrm{\ (because \ \mathit{f} \ is \ an \ isomorphism)}
\\
\subset 
f(M) \cap f\Phi_1 (EK[U] \cap EK[V])
\subset \\
f(M) \cap \Phi_2 (KEK[U] \cap KEK[V])
= f(M)_2 (KEK[U] \cap KEK[V])
\end{multline*}
shows that $f(M)_2$ is $KEK$-insular as well.
\end{proof}
The second fact concerns the construction of a full subcategory of $\mathcal{G}(X,A)$. Recall that leanness and insularity of an object $\Phi$ depend on existence of controlled sets $E_{lean}, E_{insular}$ in $X \times Y$. If we let $E_{tot} = E_{lean} \cup E_{insular}$, then $E_{tot}$ is a single controlled set which can be used to verify both conditions.  In other words, $\Phi$ is $E_{tot}$-admissible according to Definition \ref{splitleanins}(\ref{admissible}). 

\begin{Definition} \label{0lean}
    We will be interested in the full subcategory $\mathcal{G}^{\Delta}(X,A)\subseteq \mathcal{G}(X,A)$ on objects which are $\Delta$-lean, where $\Delta$ denotes the diagonal controlled subset. 
\end{Definition}

\begin{Proposition} \label{isomorphism}
    Each object in $\mathcal{G}(X,A)$ is isomorphic to an object of the subcategory $\mathcal{G}^{\Delta}(X,A)$.  
\end{Proposition}

\begin{proof}
For any controlled set $E$ in $X $ and $E$-admissible $X$-filtration $\Phi$, we can define a new $X$-filtration $\Phi ^E$ by setting $\Phi ^E(U) = \Phi (E[U])$. 
We also have  an evident morphism $\Phi \to \Phi ^E$, given at any subset $U$ by the structure map $\Phi(U) \hookrightarrow \Phi (E[U])$ induced by the inclusion $U \hookrightarrow E[U]$.  We claim this is an isomorphism of objects $\Phi \to \Phi^E$.  First, it's clear that it is controlled by $\Delta$ in the forward direction, while the inverse is controlled by $E$. 

It remains to check that $\Phi ^E$ itself is coarsely adapted and is $\Delta$-lean.
The property of $\Phi^E$ that is convenient is the fact that $\Phi^E (x) = \Phi (E[x])$, and so 
\[
\Phi^E (U) \subset \sum_{x \in U} \Phi (E[x]) = \sum_{x \in U} \Phi^E (x)
\]
for any subset $U$.  This means that, in terms of Definition \ref{splitleanins}, the object $\Phi^E$ is $\Delta$-lean.  
 We claim that $\Phi^E$ is also $E^2$-insular.  This follows from 
 \begin{multline*}
 \Phi ^E (U_1) \cap \Phi ^E (U_2) = \Phi (E[U_1]) \cap \Phi (E[U_2]) \\
 \subseteq \Phi (E^2[U_1] \cap E^2[U_2]) \subseteq \Phi ^E
 (E^2[U_1] \cap E^2[U_2]).
 \end{multline*}
Finally, we need to prove that $\Phi ^E$
is coarsely adapted as long as $\Phi$ is coarsely adapted.  The latter means that for every subset $U \subset X$ there is a control set $E_U$ and an admissible subobject $\Psi(U)$ nested between two subconstructions
\[
\Phi_U \subset \Psi(U) \subset \Phi_{E_U [U]}.
\]
In particular, the subconstruction $\Phi_X$ which coincides with $\Phi$ is $E$-admissible.  Let's consider the inclusions of subobjects
\[
\Phi^E_U = \Phi_{E[U]} \subset \Psi(E[U]) \subset 
\Phi_{E_{E[U]} [E[U]]} \subset \Phi^E_{E_{E[U]} [E[U]]}.
\]
As soon as we show that the enlargement $\Psi(E[U])$ is admissible with the respect to the induced filtration from $\Phi^E$, the required control set associated to the subset $U$ to contain an admissible enlargement can be taken to be $E_{E[U]} E$, so $\Phi^E$ is coarsely adapted.

By assumption, $\Psi(E[U])$ is admissible with respect to the filtration induced from $\Phi$.  We already know that $\Phi^E$ is isomorphic to $\Phi$.  The result  follows from Lemma \ref{Propone}
\end{proof}

\begin{Remark}
    An object $\Phi$ is $\Delta$-lean if and only if  $\Phi(U)$ is exactly equal to the sum of all the objects $\Phi(x)$ as $x$ ranges over all points in $U$. Of course, the sum is not necessarily direct.  Nevertheless $\mathcal{G}^{\Delta}(X,A)$ is in a sense closer to  the Pedersen-Weibel category.  
\end{Remark}

Suppose $X$ has a free, properly discontinuous action by a group of isometries $\Gamma$.

\begin{Corollary} \label{eqisomorphism}
     The construction of $\Phi^E$ is invariant under the action by any subgroup $\Gamma_0 \subseteq \Gamma$.  Therefore, each object in $\mathcal{G}(X,A)$ fixed under the action of $\Gamma_0$ is isomorphic to an object of $\mathcal{G}^{\Delta}(X,A)$ fixed by $\Gamma_0$.      
The isomorphism between these fixed objects is an equivariant isomorphism.
\end{Corollary}

\begin{proof}
   This is an immediate consequence of the fact that the construction in the proof of Proposition \ref{isomorphism} is equivariant, since the action is by isometries.  
\end{proof}

The main reason for the introduction of the category $\mathcal{G}^{\Delta}(X,A) $ is to provide an analysis of products (perhaps infinite) of the construction $\mathcal{G}(X,A)$.  This behavior is not very clear when studying the category $\mathcal{G}(X,A) $ directly, but it becomes clear when studying the equivalent subcategory $\mathcal{G}^{\Delta}$. 

\begin{Proposition} \label{infproducts}
    Let $X$ denote a discrete coarse space, so that the diagonal is the maximal controlled set.  Then there is an equivalence of categories 
    $$ \pi \colon \mathcal{G}^{\Delta}(X,A) \longrightarrow \prod_{x \in X} \mathcal{G}(x,A)
    $$
    If $X$ is equipped with an action of a group $\Gamma$, then $\pi$  is a $\Gamma$-equivariant equivalence of categories when $\Gamma$ acts by permuting factors.    
\end{Proposition}
\begin{Proof}
    We consider any $\Delta$-lean  $X$-filtration $\Phi$. For each $x \in X$, we consider the submodule $\Phi (x) \subseteq \Phi (X) $, and claim that the obvious homomorphism 
    $$ j \colon \bigoplus_{x \in X} \Phi (x) \rightarrow \Phi (X)
    $$
    is an isomorphism. We first observe  that $j$ is surjective, since the surjectivity statement is exactly the definition of leanness applied to the set $X$.  To prove that it is injective, any element in the kernel of $j$ is a vector $\{\phi _x\}_{x \in X}$ in the direct sum $\bigoplus _{x \in X} \Phi (x)$. This means that there is a finite subset  $S \subset X$ so that the coordinates $\phi _x$ are all zero unless $x \in S$.  On the other hand, $\Phi$ is an $X$-filtration, so the morphism $\Phi (S) \hookrightarrow \Phi (X)$ is injective, so all the coordinates $\phi _s$ for $s \in S$ are also zero, so the kernel  of $j$ is trivial.  Similarly, consider any subset $U \subseteq X$. Because $\Phi$ is $\Delta$-lean, the homomorphism
    $$ \bigoplus _{x \in U} \Phi (x) \rightarrow \Phi (U)
    $$
    is surjective, and it is also injective since the composite 
    $$\bigoplus _{x \in U}\Phi(x) \rightarrow \Phi (U) \rightarrow \Phi (X) $$ is injective. Consequently, $\Phi$ is isomorphic to the $X$-filtration $\Psi$ defined by $\Psi (U) = \sum _{u \in U} \Phi(u) $. The functor $\pi$ in the statement of the Proposition is given  by sending  $\Phi$ to the vector $  \{ \Phi (x) \}_{x \in X}$.  To see that it is an equivalence, construct a functor $\theta$ in the opposite direction which sends the vector  $\{M_x \}_{x \in X} $ of objects in the category of $A$-modules to the filtration with total module $\bigoplus _{x \in X} M_x$, and assigns to the set $U \in X$ the submodule $\bigoplus _{x \in U} M_x$. It is clear that the composites $\pi \compcirc \theta$
    and $\theta \compcirc \pi$ are both isomorphic to their respective identity functors.  
\end{Proof}




\subsection{Fibred theories}
We will need to define ``fibred" versions of our theories.  These constructions are critical to the proof of our main theorem, and are very much analogous to {\em parametrized homotopy theory } as expounded in \cite{jMjS:06} or \cite{mCiJ:98}. The point is that we will need to develop a context in which the ``cells" are copies of suspensions of the $K$-theory spectrum of $A[\Gamma]$, and the best way to achieve this is through the use of a coarse theory over a base space whose $K$-theory is $K(A[\Gamma ])$. This base space will be the group $\Gamma$ in question regarded as a metric space using the word length metric associated to a finite generating set.   Ultimately we will also require an equivariant version of the  theory, which we will see in Section \ref{equivariantsection}.  

We consider two extended metric spaces $X$ and $Y$.  Of course, we may consider the product $X \times Y$ as an extended metric space,  but it turns out we will require a different coarse structure in which the roles of $X$ and $Y$ are not interchangeable but where $X$ is regarded as a base in a product fibration and $Y$ is regarded as a fiber.

\begin{Definition} \label{fibrecontrol}
Let $X$ be a metric space and $Y$ an extended metric space. We define a coarse structure $\mathcal{E}(X,Y)$ on $X \times Y$ by declaring that $E \subseteq (X \times Y) \times (X \times Y)$ is a controlled set if and only if 

\begin{enumerate}
    \item{ the projection to the 1st and the 3rd coordinates $\pi _{X \times X}(E)$ is a controlled set in $X$,}
    \item{for every bounded set $B \subseteq X$, $\pi _{X \times X}^{-1}(B \times B) \cap E $ is a controlled set  in the product  coarse structure on $X \times Y$.  }
\end{enumerate}

\end{Definition}

\begin{Remark}There is an alternate equivalent characterization in the context of  extended metric spaces,  which says that $E$ is a controlled set
if and only if there is a positive real number $R$, a point $x_0 \in X $, and a monotone function $\theta \colon \mathbb{R}_+ \rightarrow \mathbb{R}_+$ so that $((x,y), (x^{\prime},y^{\prime})) \in E $ if and only if $d_X(x,x^{\prime}) \leq R $ and $d_Y(y,y^{\prime}) \leq \theta (d_X(x_0,x))$. It is easy to show that the statement is independent of choice of $x_0$ as well as of the choice of $x$ over $x^{\prime}$.
\end{Remark}

We now define the fibred category $\mathcal{G}_X(Y, A) $ to be $\mathcal{G}(X \times Y,A)$ as defined in Section  \ref{PWtheories}, where the coarse structure on $X \times Y$ is understood to be the structure  $\mathcal{E}(X,Y)$ defined above. Given a subset $Y_0 \subseteq Y$, we define $\mathcal{G}_X(Y,A)_{\leq Y_0}$ to be the full subcategory of $\mathcal{G}_X(Y,A)$ on $(X \times Y)$-filtrations $\Phi$ for which there exists a controlled subset $E \in \mathcal{E}(X,Y)$ so that $\Phi (E[X \times Y_0])  = \Phi(X \times Y)$. 

We also observe that the following generalization of Proposition \ref{infproducts}.

\begin{Proposition}\label{contproducts}
    Let $Y$ denote any set, regarded as a coarse space with discrete coarse structure.  Then there is a natural equivalence of categories $$\mathcal{G}^{\Delta}_X(Y,A) \cong \prod_{y \in Y}\mathcal{G}_X(y) $$
\end{Proposition}
\begin{Proof}
    The proof is an immediate generalization of the proof of Proposition \ref{infproducts}.  It is critical that the coarse structure is $\mathcal{E}(X,Y)$.  For example, the standard construction of a product  coarse structure on $X \times Y$ would not necessarily satisfy the statement of the theorem.  
\end{Proof}

\subsection{Equivariant  theories}\label{equivariantsection}
When working in equivariant $K$-theory or more generally in equivariant spectra, there is often a naive structure whose group action is not very meaningful.  For instance, if we consider the category $\mathcal{B}(X,A)$ for a metric space $X$ with a free isometric action by a group $\Gamma$, there is an evident $\Gamma$-action on $\mathcal{B}(X,A)$, but the fixed point category of that action consists of only the zero module.  This situation comes up frequently, for instance in the discussion of equivariant stable homotopy or equivariant complex $K$-theory.  The solution is to produce a larger category together with an equivariant inclusion of the original category, but where the larger category has a more meaningful fixed point set.  There is a systematic way of constructing this larger category, originally due to Grothendieck, and discussed in \cite{rT:82} and \cite{mM:16}.

Let $\underline{\Gamma}$ denote the category whose objects are the elements of $\Gamma$ and where there is a  unique morphism between any pair of elements.  The morphisms can be viewed as left multiplication by elements of $\Gamma$, and $\Gamma$ acts on $\underline{\Gamma}$ by right multiplication.  

\begin{Definition}\label{grothendieck}
    Let $\underline{C}$ denote a category with action by a group $\Gamma$.  Then the Grothendieck construction on $\underline{C}$, or the ``equivariant version of $\underline{C}$", is the category $\underline{C}_{\Gamma}$ whose objects are the functors from $\underline{\Gamma}$ to $\underline{C}$, and where the morphisms are the natural transformations of functors.   The right multiplication action of $\Gamma$ on $\underline{\Gamma}$ gives a left action on $\underline{C}_{\Gamma}$.  
\end{Definition}
\begin{Example}
    If the $\Gamma$-action on $\underline{C}$ is trivial, then $\cat{C}_{\Gamma}$ is the {\em category of representations of $\Gamma$ in $\underline{C}$}, i.e. the category whose objects are pairs $(c,\rho)$, where $c$ is an object of $\underline{C}$ and $\rho \colon \Gamma \rightarrow Aut(c)$ is a group homomorphism, and where the morphisms from $(c,\rho)$ to $(c^{\prime}, \rho ^{\prime}) $ are the morphisms from $c$ to $c^{\prime}$ in  $\underline{C}$, with no equivariance requirement.  The $\Gamma$-action on $\underline{C}_{\Gamma}$  is given by $g \cdot (c, \rho) = (c, g \hspace{-.07cm}\cdot \hspace{-.07cm}\rho \hspace{-.07cm}\cdot \hspace{-.07cm} g^{-1})$, and the fixed point category is therefore the category whose objects are pairs $(c,\rho)$, and where the morphisms are  the equivariant morphisms in $\underline{C}$. 
 
\end{Example}

Given a subcategory $\cat{D} \subseteq \cat{C}$ which is closed under the $G$-action, we can form the full subcategory of $\cat{C}_{\Gamma}$ whose objects are the functors from $\underline{\Gamma}$ to $\underline{C}$ with image in $\cat{D}$, and denote it by $\underline{C}_{\Gamma,\cat{D}}$.   It is clearly preserved by the $G$-action on $\cat{C}_{\Gamma}$. 

\begin{Example}
    If $\cat{D}$ consists of the object $e \in \Gamma$, and which therefore is isomorphic to the trivial category with one object and one morphism, then $\cat{C}_{\Gamma, \cat{D}}$ is equivalent to the original category $\cat{C}$ with the given  action. 
\end{Example}
\begin{Example}\label{groupring}
    Let $\cat{C}$ be $\mathcal{B}(\Gamma,A)$, where $\Gamma$ is the word metric space of a finitely generated group $\Gamma$ with respect to a chosen finite generating set $\Omega$.  The metric space is defined using the right multiplication action of $\Omega $ on $\Gamma$, and  $\Gamma$ acts isometrically on itself through left multiplication. Let $\cat{D} \subseteq \cat{C}$ be the category of {\em filtration zero morphisms}, where a filtration zero morphism from $(F_0, B_0, \varphi _0)$ to  $(F_1, B_1, \varphi _1)$ is an $A$-module homomorphism $f$  from $F_0$ to $F_1$ so that for $\beta \in B_0$, we have $f(\beta) \in \mbox{span}(\varphi_1 ^{-1}(\varphi_0(\beta))) $.  Then it is proved in \cite{gC:95} that the fixed point category of the action of $\Gamma$ on $\cat{C}_{\Gamma}$ is equivalent to the category of free  left $A[\Gamma]$-modules. The identification of the category of free $A[\Gamma ]$-modules is the key ingredient in the proofs of the split injectivity of the assembly map appearing in \cite{gC:95}, \cite{gCbG:04}, and \cite{gCeP:95}. 
\end{Example}

A key example of interest is the situation of fibred $G$-theory $\mathcal{G}_X(Y)$,  where $X$ is $\Gamma$ as in Example \ref{groupring} and $Y$ is a metric space with a coarse action by $\Gamma$.  The $\Gamma$-action is diagonal on $X \times Y$.  This produces a $\Gamma$ action on $\mathcal{G}_X(Y)$.  Given a cosheaf $\Phi$ on a space $X$ and a map $f \colon X \rightarrow Y$, we define the {\em direct image cosheaf}  $f_* \Phi$ on $Y$ by the requirement $f_*\Phi (V) = \Phi (f^{-1}(V))$.  The direct image construction is functorial for cosheaves on $X$. 
We define  $\cat{D}(X,Y)  \subseteq \mathcal{G}_X(Y)$ to be the subcategory containing all objects of $\mathcal{G}_X(Y)$ , and for which  the morphisms $f$ satisfy the requirement that
$(\pi _X)_*(f) $ is filtration zero.  

\begin{Definition}
Let  $X$ be an extended metric space (which can be regarded as a coarse space), and let $Y$ be a coarse space.  Suppose that a group $\Gamma$ acts isometrically on $X$, and by coarse maps on $Y$. The actions give a $\Gamma$-action on $\mathcal{G}_X (Y)$, and  we define  $\mathcal{G}_{X, \Gamma}(Y)$ to be the category $\mathcal{G}_X(Y)_{\Gamma, \cat{D}(X,Y)}$.   Let $Y_0 $ be a subset of $Y$, so that for every $\gamma \in \Gamma$, we have $\gamma\cdot Y_0 \subseteq E[Y_0]$ for some controlled set $E$ in $Y \times Y$. Then there is a corresponding construction ${\mathcal{G}}_{X, \Gamma}(Y) _{\leq Y_0}$, defined in the obvious way.  
\end{Definition}

We have the following analogue of Proposition \ref{infproducts}.

\begin{Proposition} \label{eqproducts}
    Let $X$ be an extended metric space, and $Y$ a set regarded as a coarse space with a discrete coarse structure.  Suppose further that  $X$ is equipped with an isometric left action by a group $\Gamma$ and that $Y$ is equipped with an action by $\Gamma$.  Then there are $\Gamma$-equivariant equivalences of categories 
    $$ \mathcal{G}_X^{\Delta}(Y,A) \cong \prod_{y \in Y} \mathcal{G} _X(y,A) \cong \prod _{y \in Y} \mathcal{G}_X(y _0) 
    $$
    and
    $$ \mathcal{G} _{X,\Gamma}^{\Delta}(Y,A) \cong \prod_{y \in Y} \mathcal{G}_{X} (y, A)
    $$
    where the actions on the right hand side of the equivalences by an element $\gamma \in \Gamma$ are given by permuting factors via right multiplication by $\gamma^{-1}$ together with induced functors $\gamma _{\dotr} \colon \mathcal{G}_X(y,A) \rightarrow \mathcal{G}_X(\gamma \cdot y,A) $.
\end{Proposition}
\begin{Proof}
    This is a  straightforward combination of Propositions \ref{infproducts} and \ref{contproducts}.
\end{Proof}

The construction $Y \rightarrow \prod_{y \in Y} \mathcal{G}_X(y,A)$ is naturally a contravariant functor in the set $Y$, since it can be interpreted as the space of functions from the set $\Gamma$ into categories. This functoriality structure is not useful to us, but   there is a covariant functoriality on a modified version of the product construction from the category $\underline{Sets}^p$ of sets and proper maps, which was constructed in Section II of \cite{gC:95}. The modification of the construction that is needed is the replacement of $\prod_{y \in Y}S_y$ by an equivalent but more flexible construction, which we summarize  briefly.  Let $Y$ be any set, and $S$ a spectrum.  We consider the category $\mathcal{S}_Y$  of finite subsets of $Y$, and define $\mathbb{G}_X^{\mathit{lf}} (Y,A)$ to be the homotopy inverse limit over $\mathcal{S}_Y$ of the functor $S \mapsto \mathbb{G}_X(S,A)$. For every set $Y$ (regarded as a discrete coarse space), there is an equivalence of spectra $\mathbb{G}_X^{\mathit{lf}}(Y,A) \rightarrow \mathbb{G}(Y,A) $, induced by a restriction of homotopy inverse limits along the inclusion of subsets of size 1, regarded as a category with only identity morphisms.  It is not natural for proper maps of sets, since $\mathbb{G}_X(Y,A)$ is not covariantly functorial, but $\mathbb{G}_X^{\mathit{lf}}(Y,A)  $ is covariantly functorial in $Y$ for proper maps of sets. 

The theory $\mathbb{G}_{\Gamma}(-,A)$ has an important property when applied to bounded actions, as defined in Definition \ref{bddefinition}. Suppose that $X$ is an extended metric space with a bounded $\Gamma$-action by coarse self-maps.  Let $X^0$ denote the same extended metric space, but with trivial $\Gamma$-action.  Then we have the following. 

\begin{Proposition}\label{trivialization}
    With $X$ and $X^0$ as above, there is a natural $\Gamma$-equivariant equivalence of spectra 
    $$ \eta _X\colon\mathbb{G}_{\Gamma}(X, A) \rightarrow \mathbb{G}_{\Gamma}(X^0,A) 
    $$
\end{Proposition}
\begin{Proof}
    For any extended metric space $X$ with coarse $\Gamma$-action, we define the associated $\Gamma$-action on $\Gamma \times X$ by $\gamma \cdot (g,x) = (\gamma g, \gamma x)$.  We also have an equivariant isomorphism (of $\Gamma$-sets, for now)  
    $\phi _X \colon \Gamma \times X \rightarrow \Gamma \times X^0$, defined by $\phi_X(g,x) = (g,g^{-1}x)$.  In the case where $X$ has a bounded $\Gamma$-action, it is easy to see that $\phi _X$  is in fact an isomorphism of coarse spaces, when we assign both $\Gamma \times X$ and $\Gamma \times X^0$ the coarse structure $\mathcal{E}(\Gamma , X) $ as defined in Definition \ref{fibrecontrol}. Note that the coarse structure is independent of the action, so applies to both $\Gamma \times X$ and $\Gamma \times X^0$. The result follows easily, generalizing the detailed discussion in section 7.2 in \cite{gCbG:19}, particularly Lemmas 7.5 and 7.6.  
\end{Proof}

\subsection{Constructing \textit{G}-theory spectra}\label{spectrum}

There are a number of approaches  that construct spectra from categorical data, one of which is the {\em exact category} approach of Quillen (see \cite{bK:96} or \cite{mS:04}) and another of which is {\em Waldhausen categories} \cite{fW:83}.  We will construct exact categories from the categories $\mathcal{G}(X,A)$, and more generally $\mathcal{G}(X,A)_{\leq Y}$, and follow that construction by a standard construction which assigns to an exact category a Waldhausen category \cite{fW:83}.  One could choose a split exact structure for the two categories, but as we show in Appendix \ref{karoubi} this choice would not permit  us to prove the desired excision results.  

\begin{Definition}
    Let $f \colon F \rightarrow F^{\prime}$ be a morphism in $\mathcal{G}(X,A)$ or $\mathcal{G}(X,A)_{\leq Y}$.  Then we say that $f$ is {\em boundedly bicontrolled} if there exists a controlled set $E$ so that 
    \begin{enumerate}
        \item{$f(F(S)) \subseteq F^{\prime}(E[S])$ for all $S \subseteq X$. }
        \item{$f(F(X))\cap F^{\prime}(S) \subseteq f(F(E[S]))$ for all $S \subseteq X$.}
    \end{enumerate}

We say a map $f \colon F \rightarrow F^{\prime}$, in either $\mathcal{G}(X,A)$ or $\mathcal{G} (X,A)_{\leq Y}$, is an {\em admissible monomorphism} (respectively an {\em admissible epimorphism})   if $f$ is boundedly bicontrolled and $f \colon F(X) \rightarrow F^{\prime}(X)$ is monic (respectively epic). 
\end{Definition}

It is proved in \cite{gCbG:11}, Corollary 2.23,  that ${\mathcal{G}}(X,A)$ and $\mathcal{G}(X,A)_{\leq Y}$ are exact categories with this choice of admissible monics and epics. 

In \cite{fW:83}, F. Waldhausen produces a simplicial set $wS_{\dotr}\cat{E}$ for any exact category $\cat{E}$.  Iteration produces a spectrum 
$$ n \to \Omega \, | \mathop{w \underbrace{S_{\dotr}\hspace{.17cm}.\hspace{.17cm}.\hspace{.17cm}. \hspace{.17cm}S_{\dotr}}}_{\ \ \ n\mathrm{-fold}} \cat{E} \, |
$$ 
which we  call the $K$-theory of $\cat{E}$.  We  define our $G$-theory  spectra using  this construction.
 \begin{Remark}
    This construction produces $\Omega$-spectra, i.e. prespectra $\{X_n \}_{n \geq 0}$ so that the adjoints $X_n \rightarrow \Omega X_{n+1}$ are weak equivalences. Also, if $\cat{E}$ is an exact category with action by a group $\Gamma$, then the construction produces a {\em naive $\Gamma$-spectrum}, i.e. an $\Omega$-spectrum with $\Gamma$ action, so that the adjoints are weak $\Gamma$-equivariant equivalences.  We are not adding any additional structure such as deloopings by one-point compactifications of representations. The fixed point spectra are the spectra associated to the fixed point subcategories, which are exact subcategories.   
 \end{Remark}
 We enumerate the particular exact categories we will be working with.  
\begin{enumerate}
    \item{Given a coarse space $X$, we write $G(X,A)$ for the spectrum constructed from the exact category $\mathcal{G}(X,A)$, with the admissible monomorphisms (respectively epimorphisms) being the boundedly bicontrolled monomorphisms (respectively epimorphisms). For a subset $X_0 \subseteq X$, we also have a spectrum $G(X,A) _{\leq X_0}$ associated to the category $\mathcal{G}(X,A) _{\leq X_0}$. }
    \item{Given coarse spaces $X$ and $Y$, we write $G_X(Y,A)$ for the spectrum constructed from the category $\mathcal{G}_X(Y,A)$ given the exact structure with the admissible monomorphims (respectively epimorphisms) consisting of the boundedly bicontrolled monomorphisms (respectively epimorphisms) with respect to the coarse structure $\mathcal{E}(X,Y)$ from Definition \ref{fibrecontrol}. For a subset $Y_0 \subseteq Y$ we obtain a spectrum $G_{X}(Y)_{\leq Y_0}$ associated to the category $\mathcal{G}_{X}(Y,A)_{\leq Y_0}$. }
    \item{Let $\Gamma$ be a group.  Given an extended metric space $X$ and a coarse space $Y$, an isometric action of $\Gamma$ on $X$, and a coarse action of $\Gamma$ on $Y$, we obtain the category $\mathcal{G}_{X, \Gamma}(Y,A)$ with an action by $\Gamma$. The morphisms in this category are natural transformations $\{ N_{\gamma} \}_{\gamma \in \Gamma}$, and we declare such a natural transformation to be an admissible monomorphism (respectively epimorphism) if each $N_{\gamma}$ is boundedly bicontrolled.  This choice produces a spectrum which we denote by $G_{X, \Gamma}(Y,A)$.  Given a subset $Y_0 \subseteq Y$  so that for each $\gamma \in \Gamma$, there is a controlled set $E_{\gamma}$ so that  $\gamma \cdot Y_0 \subseteq E[Y_0] $, we also obtain a spectrum $G_{X, \Gamma}(Y,A)_{\leq Y_0 }$. Of course, the spectra $G_{X,\Gamma}(Y,A) $ and $G_{X,\Gamma}(Y,A)_{\leq Y_0} $ both are equipped with $\Gamma$-actions.  }
\end{enumerate}

\begin{Notation}
    For our applications, $\Gamma$ will always be a finitely generated group equipped with a word metric defined using right multiplication by elements in a chosen finite generating set.  The base space $X$ will always be the group $\Gamma$ with the isometric action by left multiplication on itself.  To simplify notation, therefore, we will write $\mathbb{G}_{\Gamma}(Y,A)$ for $G_{\Gamma, \Gamma}(Y,A)$ and $\mathfrak{G}_{\Gamma}(Y,A)$ for the category $\mathcal{G}_{\Gamma, \Gamma}(Y,A)$.    
\end{Notation}

The category $\mathfrak{G}_{\Gamma}(Y,A)$, and therefore the spectrum $\mathbb{G}_{\Gamma}(Y,A)$, admits a $\Gamma$-action.  The goal will be to  understand the $\Gamma$-action on $\mathbb{G}_{\Gamma}(Y,A)$. Specifically, we wish to identify the fixed point spectrum $\mathbb{G}_{\Gamma}(Y,A)^{\Gamma}$ as the spectrum associated to a choice of cofibrations and weak equivalences on the   fixed point category $\mathfrak{G}_{\Gamma}(Y,A)^{\Gamma}$. The tools needed to draw these conclusions are Propositions \ref{requirements} and \ref{otherrequirements} of Appendix~\ref{groupactions}. 
\begin{Proposition}
    The set of morphisms $co(\mathfrak{G}_{\Gamma}(Y,A))^{\Gamma} $ endows $\mathfrak{G}_{\Gamma}(Y,A)^{\Gamma}$ with the structure of a category with cofibrations and weak equivalences, where we let the weak equivalences be the minimal choice, i.e. the isomorphisms.  
\end{Proposition}
\begin{Proof} It suffices to show that $(\mathfrak{G}_{\Gamma}(Y,A), co(\mathfrak{G}_{\Gamma}(Y,A)) , *)$ satisfies the hypotheses of Proposition \ref{requirements}. 
The reason for that can be traced back to the basic non-equivariant, non-fibred version of the fact in the form of closure under co-base changes in $\mathcal{B}(X,A)$, checked as part of Theorem 2.13 of \cite{gCbG:11}.  First, we observe that all cofibrations have cokernels, from Proposition 2.6 and part (c) of Proposition 2.18 in \cite{gCbG:11}, and so are admissible monomorphisms in the non-fibred controlled category.  The key is that the fixed point subcategory has the same description of admissible monomorphisms and the construction of cokernels is entirely the same.  There are analogues of the the above facts in the fibred setting as well in Theorem 3.2.5 and Proposition 3.2.7 in \cite{gCbG:18}.  All together this verifies both conditions in Proposition \ref{requirements} for $\mathfrak{G}_{\Gamma}(Y,A)$. 
    \end{Proof}
\begin{Proposition}\label{MM}
    The $K$-theory spectrum of the category with cofibrations and weak equivalences $$(\mathfrak{G}_{\Gamma}(Y,A)^{\Gamma}, co(\mathfrak{G}_{\Gamma}(Y,A))^{\Gamma}, *)  $$ is isomorphic to the fixed point spectrum of the $\Gamma$-action on $\mathbb{G}_{\Gamma}(Y,A)$. 
\end{Proposition}
\begin{Proof} It suffices to show that $(\mathfrak{G}_{\Gamma}(Y,A), co(\mathfrak{G}_{\Gamma}(Y,A)) , *)$ satisfies the hypotheses of Proposition \ref{otherrequirements}. 
In fact, stronger facts hold in this application.  It follows from the proof of the preceding Proposition that the hypotheses of Proposition \ref{otherrequirements} are satisfied for any group action on $\mathfrak{G}_{\Gamma}(Y,A)$ which is induced from bounded actions on the space $Y$.  Since the zero object $*$ is a terminal object in each fixed point subcategory considered, and the cofibrations are admissible monomorphisms, it suffices to observe that cokernels in $\mathfrak{G}_{\Gamma}(Y,A)^{\Gamma'}$ coincide with the cokernels of the same cofibrations in $\mathfrak{G}_{\Gamma}(Y,A)$.  Just as in the proof of Proposition 3.1, this can be tracked from the basic non-equivariant, non-fibred version of the fact through the same sequence of theorems in the literature.  This time the observation is that the same construction of the cokernel in $\mathfrak{G}_{\Gamma}(Y,A)$ gives the cokernel in each of the fixed point subcategory.
\end{Proof}

\begin{Remark}
(a) We will be requiring the spectrum level homotopy fixed point constructions on the spectra constructed via $\mathbb{G}_{\Gamma}(-, A)$. This means that one must ensure that equivariant function spaces  have homotopy invariant meaning.   One could achieve this by, for instance, applying geometric realization or Kan's $Ex^{\infty}$ construction level-wise. However, we  may work directly with Waldhausen's construction due to the results of \cite{gC:95.1}, which ensure that Waldhausen's $S_{\dotr}$ construction is ``quasi-Kan", and therefore homotopy invariant to function space constructions.

(b) An alternative way to think about $\mathfrak{G}_{\Gamma}(Y,A)^{\Gamma}$ as a Waldhausen category is due to Malkiewich and Merling \cite{cMmM:19}.  They call this type of construction the \textit{homotopy fixed points} of a category with a group action and show that this construction always has a good definition of a Waldhausen structure.  Our definition is more basic.  It makes sense already at the level of $\mathcal{G}_{\Gamma}(Y,A)$ which is convenient because it can be readily applied in a variety of categories with fibred control in this paper.  The action by an element of the group in our case preserves the $0$ object and pushouts on the nose.  However, it's easy to see that the construction matches \cite{cMmM:19} in their setting. (Nota bene: before section 2.4 in \cite{cMmM:19} everything is true for a finitely generated group as written, not just finite groups.)
\end{Remark}

\begin{Notation} \label{MM2}
    It is often convenient to use the shorthand notation $\mathcal{G}_{\Gamma} (Y,A)^{h\Gamma}$ for the \textit{homotopy fixed points} category $\mathfrak{G}_{\Gamma}(Y,A)^{\Gamma}$ with the Waldhausen structure as above.  Then Proposition \ref{MM} states there is an equivalence
    \[
    \mathbb{G}_{\Gamma}(Y,A)^{\Gamma} = K(\mathcal{G}_{\Gamma} (Y,A)^{h\Gamma}).
    \]
\end{Notation}

We conclude this section with the following analysis of products. 

\begin{Proposition} \label{prodcriterion}
    Let $X$ be an extended metric space with an isometric action by a group $\Gamma$, and let $Y$ be the  set $\Gamma$ regarded as a coarse space with the discrete coarse structure, and equipped with the left multiplication by $\Gamma$.  Then there is an equivariant equivalence of spectra 
    $$ \mathbb{G}_X (Y,A) \cong F(Y, \mathbb{G}_X(*,A)) 
    $$
    where the function space is equipped with the conjugation action on maps.  Moreover, if $Y$ is instead an extended metric space  with the property that all of its coarse components are bounded in diameter by a uniform bound, then the map $Y \mapsto \pi  _0^{coarse}(Y)$ induces an equivariant equivalence on $\mathbb{G}_X(-, A)$-spectra.  
    \end{Proposition}
    \begin{Proof} 
        The non-equivariant statement is simply that the $K$-theory construction carries infinite products of exact categories to infinite products, as is proved in \cite{gC:95.1} or \cite{dKcW:20}.  To get the equivariant result, we first note that the categories defining $\mathbb{G}$ satisfy the hypotheses of Proposition \ref{otherrequirements}, so that the $G$-theory of the fixed point categories is equivalent to the fixed points of the corresponding spectra.  The result now follows from the fact that the fixed point categories of the $\Gamma$ action on the product categories are themselves product categories, that is the statement that 
        $$(\prod _{Y} \mathcal{G}(y,A))^{\Gamma _0} \cong \prod _{\alpha \in \Gamma / \Gamma _0} \mathcal{G}(g_{\alpha},A)
        $$
        where for each $\alpha$, $g_{\alpha}$ denotes a left coset representative of the coset $\alpha$. We can then apply \cite{dKcW:20} to obtain an equivalence on the fixed point sets of all subgroups $\Gamma _0$, and therefore an equivariant equivalence.  
        \end{Proof} 

        \begin{Remark} The construction $\mathbb{G}_X(Y,A) $ is covariantly functorial in $Y$ on the category of sets and proper maps of sets, while the construction $F(Y, \mathbb{G}_X(*, A))$ is not.  This means that by applying $\mathbb{G}_X(Y,A) $ to any proper simplicial set (regarded as a simplicial object in the coarse category) we obtain a simplicial spectrum, which will be important to us.  We use the equivalence in Proposition \ref{prodcriterion} to perform levelwise equivariant analyses on these simplicial spectra.  
        \end{Remark}

        \begin{Corollary}\label{fixedanalysis}
        There is a natural equivalence of functors 
$$ \mathbb{G}_X(Y,A)^{\Gamma} \cong \mathbb{G}_X(\Gamma \backslash Y,A) 
$$
on the category of free simplicial $\Gamma$-sets with finitely many orbits in each level.  
        \end{Corollary}

We conclude this section with an observation concerning transformation $\eta _X$ from Proposition \ref{trivialization}.
        \begin{Proposition}\label{etacompatibility}
            The transformation $\eta _X$  preserves the localization constructions $(-)_U$, so induces  $\eta _X \colon \mathbb{G}_{\Gamma} (X,A)_{\leq U} \rightarrow \mathbb{G}_{\Gamma}(X^0,A)_{\leq U} $.  The correspondence is natural with respect to inclusions $U \hookrightarrow U^{\prime} $, and so  produces transformations on any diagram of subsets of $X$. 
        \end{Proposition}
        \begin{proof}
            The proof is immediate given the definition of the coarse structure $\mathcal{C}(\Gamma, X)$ on $\Gamma \times X$ and $\Gamma \times X^0$.  
        \end{proof}

\subsection{Relative constructions} 
A critical component of our proof is an excision theorem for the equivariant spectra $\mathbb{G}_{\Gamma}(Y,A)$ and $\mathbb{G}_{\Gamma}(Y,A)_{\leq Y_0}$.  Such a theorem  requires the construction of quotient categories of $\mathcal{G}_{X, \Gamma}(Y,A) $ by subcategories $\mathcal{G}_{X,\Gamma}(Y_0,A)$, where $Y_0 \subseteq Y $ is a $\Gamma$-invariant subset.  In particular, one would want fibration sequences of spectra of the form 
$$ \mathbb{G}_{\Gamma}(Y_0,A) \longrightarrow \mathbb{G}_{\Gamma}(Y,A) \longrightarrow \mathbb{G}_{\Gamma}(Y,Y_0,A)
$$
where $\mathbb{G}_{\Gamma}(Y,Y_0,A) $ is the $K$ theory of a suitable   quotient category $$\mathfrak{G}_{\Gamma}(Y,A)/\mathfrak{G}_{\Gamma}(Y_0,A)$$  There are two obstacles to this plan. The first is that the exact categories $\mathfrak{G}_{\Gamma}(Y,A)$ are not the exact categories associated to abelian categories, and therefore are not amenable to the standard quotient category methods applicable in the case of abelian categories.   It is clear that the  module categories  defined in \cite{gCbG:03} and \cite{gCbG:16}  over a group ring $A[\Gamma]$  do not form an abelian category,  and that fact is reflected in the choice of admissible monomorphisms and epimorphisms in the exact category structure we have defined.  If they had been abelian categories, there would be straightforward choices of quotient categories and corresponding localization sequences as in \cite{dQ:72}. We are therefore forced to resort to a more general notion of quotients of exact categories, namely Schlichting's quotient construction \cite{mS:04}.   The second obstacle is the ``coarse nature" of the homology theories in question, which dictate that we should not deal with the subspace $Y_0$ as it stands but rather with the increasing family of subspaces $E[Y_0]$ as $E$ ranges over the controlled subsets of $Y$. This suggests that we should construct the quotient category $\mathfrak{G}_{\Gamma}(Y,A)/\mathfrak{G}_{\Gamma}(Y,A)_{\leq Y_0} $ instead.  These obstacles are overcome in \cite{gCbG:11}, \cite{gCbG:16}, \cite{gCbG:18}, and \cite{gCbG:19}, to obtain the following result. 
\begin{Theorem}
    Let $Y$ be a coarse space with a coarse action by a group $\Gamma$, and let $Y_0 \subseteq Y$ be a subspace with the property that for every $\gamma \in \Gamma$, there is controlled subset $E$ in $Y \times Y$ so that $\gamma Y_0 \subseteq E[Y_0]$.  Then the subcategory $\mathfrak{G}_{\Gamma}(Y,A)_{\leq Y_0}$ of $\mathfrak{G}_{\Gamma}(Y,A)$ has enough properties that guarantee there is an exact quotient category $\mathfrak{G}_{\Gamma}(Y,A)/\mathfrak{G}_{\Gamma}(Y,A)_{\leq Y_0} $ with the \textit{K}-theory that can be denoted $\mathbb{G}_{\Gamma}(Y,Y_0,A)$ and the corresponding fibration sequence $$ \mathbb{G}_{\Gamma}(Y_0,A) \longrightarrow \mathbb{G}_{\Gamma}(Y,A) \longrightarrow \mathbb{G}_{\Gamma}(Y,Y_0,A).
$$
\end{Theorem}

This sequence is the crucial tool for proving the Excision Theorem \ref{ExcGalltwo} stated in section \ref{ExcThms}.

\subsection{Asymptotic transfer in parametrized \textit{G}-theory} \label{transfer}

First, we discuss the absolute case.  Let $X$ be any locally finite simplicial complex of dimension $\le D$. There is a natural path-length metric $d_X$ on $X$ defined as follows. 

\begin{Definition} \label{metric-simplicial}
Let $\Sigma_X$ denote the set of simplices of $X$, and let $U$ denote the subset $\bigcup_{\sigma \in \Sigma_X} \sigma \times \sigma$ of $X \times X$. 
For any pair $(x, x') \in \sigma \times \sigma$ we define $\phi (x, x')$ to be the distance between the points $\theta(x)$ and $\theta(x')$, where $\theta$ is a chosen linear homeomorphism from $\sigma$ to the standard
Euclidean simplex of the same dimension. We define $d_X$ to be the metric $d_{(U,\phi)}$ as in Definition \ref{LenMet}.  Recall from Notation \ref{NotTX} that $CX$ stands for the cone construction $C_f X$ with $f(r) = 2^r$ for $r \ge 0$. 
\end{Definition}

We claim there  is a discrete subspace $\mathcal{D}_X \subset CX$ such
that the inclusion $\mathcal{D}_X \hookrightarrow CX$ is a coarse equivalence.  It suffices to find a commensurable discrete subspace $\mathcal{D}_X$, which means
there is a real number $\rho$ so that for every point $x \in X$, there is a point $d \in \mathcal{D}_X$
so that $d(x, d) \le \rho$. In order to identify $\mathcal{D}_X$, recall that for a simplex $\sigma$ of dimension
$D$, we have that the diameter of all the simplices in the barycentric subdivision of $\sigma$ is bounded from above by $\frac{D}{D+1} \diam(\sigma)$. We choose $k$ so that $\big( \frac{D}{D+1} \big)^k \le \frac{1}{2}$, and define $V(i) \subset X$ to be the set of vertices in the $i$-th barycentric subdivision $SX$ of $X$.  $V (0)$ is the set of vertices in $X$. The set $\mathcal{D}_X$ will be the set
\[
\bigcup_{i \ge 0} \frac{i}{k} \times V(i) \cup \bigcup_{i \le 0} i \times V(0).
\]
It is clear that $\mathcal{D}_X$ is commensurable with $CX$, and so the inclusion is a coarse equivalence.

We wish to define a simplicial complex structure for $CX$ whose vertex set is $\mathcal{D}_X$.  Let $S(X)$ denote the barycentric subdivision of a simplicial complex $X$. We will define a triangulation $\mathcal{T}$ of $\Delta[n] \times I$, which when restricted to $\Delta[n]\times \{ 0 \}$ is the standard triangulation of $\Delta [n]$, and when restricted to $\Delta [n] \times \{ 1 \}$ is $S(\Delta [n]$).  The vertex set of $\mathcal{T}$ is $\mathcal{V} = \mathcal{V}_0 \coprod \mathcal{V}_1$, where $\mathcal{V}_0 = \{ 0, \ldots , n \}$, and $\mathcal{V}_1$  denotes $\mathcal{P}^0 (\{ 0, \ldots , n \})$,  where $\mathcal{P}^0(X) $ denotes the collection of non-empty subsets of $X$.  A non-empty subset  $S$ of $\mathcal{V}$ spans a simplex in $\mathcal{T}$ if  (a)  $S \cap  \mathcal{V}_1 = \emptyset$ or (b) $S \cap \mathcal{V}_1$ is a non-empty totally ordered collection of subsets of $\{ 0, \ldots , n \} $ and $S \cap \mathcal{V}_0 \subseteq S^{\prime}$ for all $S^{\prime} \in S \cap \mathcal{V}_1$.
This construction is natural for inclusions of sets, and it follows that it extends to the definition of a triangulation of $|X| \times I$ for any simplicial complex $X$, and we denote it by $\mathcal{I}(X)$. 
We have inclusions $i_0 \colon X \to  \mathcal{I} (X)$ and $i_1 \colon \Sdd (X) \to  \mathcal{I} (X)$. We form the disjoint union $\coprod_{i \ge 0} \mathcal{I} ( \Sdd^i (X) )$ and identify the image $i_0 (\Sdd^i(X)) \subset  \mathcal{I} (\Sdd^i (X))$ with the image $i_1 (\Sdd^i(X)) \subset  \mathcal{I} (\Sdd^{i-1} (X))$  for each $i$, to obtain a triangulation of $X \times  [0, \infty)$ whose vertex set is
\[
\mathcal{D}_X \cap (X \times  [0, \infty)) = \bigcup_{i \ge 0} \frac{i}{k} \times V(i)
\]
Note that the interval in the decomposition
$\vert  \mathcal{I} (\Sdd^i(X)) \vert \cong \vert X \vert \times I$
corresponds to the interval $[ \frac{i}{k} , \frac{i+1}{k} ]$ in $[0, \infty)$. 
We will also need to construct
a triangulation of $X \times  (-\infty, 0]$ whose vertex set is
\[
\mathcal{D}_X \cap (X \times  (0, -\infty),0]) = \bigcup_{i \le 0} i \times V(0).
\]
This is much simpler, since there are no increasingly refined subdivisions. Recall that for any simplicial complex $X$, one may obtain a triangulation of $\vert X \vert \times I$ which agrees with the standard one on $\vert X \vert \times 0$ by selecting a total order on the vertex set. That construction proceeds as follows. The vertex set of the simplicial complex is $V_X \times \{0, 1\}$, and a collection of vertices
\[
\{(v_0,0), (v_1,0), \ldots ,(v_s,0),(v_{s+1},1), \ldots, (v_t,1)\}
\]
forms a simplex if and only if there are re-orderings $\sigma$ and $\tau$ of the sets $\{0, \ldots , s\}$
and $\{s + 1,\ldots,t\}$ respectively so that
\[
v_{\sigma(0)} < v_{\sigma(1)} < \dots < v_{\sigma(s)} < v_{\tau(s+1)} < \dots v_{\tau(t)}
\]
where $<$ denotes the chosen ordering on the vertex set of $X$. We denote this simplicial complex by $\mathcal{J}(X)$, and note that we have two inclusions $j_0 \colon  X \to \mathcal{J}(X)$ and $j_1 \colon  X \to \mathcal{J}(X)$. We are now able to form the disjoint union of copies of $\mathcal{J}(X)$,  
$\coprod_{i \ge 0} \mathcal{J}_i (X)$, 
and glue the images $i_1(X) \subset  \mathcal{J}_i (X)$ to the images $i_0(X) \subset  \mathcal{J}_{i+1} (X)$ to obtain a triangulation of $X \times [0, +\infty)$. 
After applying multiplication by $-1$ as a map from $X \times [0, +\infty)$, we obtain a triangulation of $X \times (-\infty,0]$ the restriction of which to the subset $X \times \{0\}$ is clearly identified with the restriction of the restriction of the triangulation of $X \times [0, \infty)$ obtained earlier. The conclusion is that we have constructed a triangulation of all of $CX$. Examination of the construction gives the following result.

\begin{Proposition} \label{HFVH}
	 Given a finite-dimensional locally finite simplicial complex $X$, there is triangulation of $CX$ so that there is a uniform bound on the diameters of the simplices of the triangulation. We will write $\mathcal{C}(X)$ for the simplicial complex, so $\vert \mathcal{C}(X) \vert$ is $CX$. We will also use the associated simplicial set $\mathcal{C}(X)_{\dotr}$.  The natural map $X \to \vert X \vert \times \{0\} \subset \mathcal{C}(X)$ is simplicial with the triangulation of $\vert X \vert$ being the given one corresponding to the simplicial structure $X$. The projection $\mathcal{C}(X) \to \mathbb{R}$ is simplicial with respect to the triangulation as above, with the vertices in $\{ n/k \}$ for the chosen value of $k$.
\end{Proposition}

Consider any locally finite simplicial complex $W$ with a properly discontinuous free $\Gamma$-action, i.e. so that $\Gamma$ acts freely on the set of all simplices in $W$. 
We  assume $W$ is subdivided so that (a)  the orbit space is a simplicial complex which we denote $\Gamma \backslash W$, and (b) so that there is an associated simplicial set $W_{\dotr}$ whose realization is naturally homeomorphic to the simplicial complex $W$ and whose non-degenerate simplices are in bijective correspondence to the simplices of $W$.  
By viewing the sets of simplices of $W_{\dotr}$ as discrete metric spaces, we obtain a simplicial category, yielding a simplicial spectrum with $\Gamma$-action and eventually the equivariant spectrum ${G}_{\Gamma} (W_{\dotr})$. The following result will not be required for the proof of our result, but will add insight into the intuition behind the proof. 

\begin{Proposition} \label{UFSV}
	The fixed point spectrum ${G}_{\Gamma} (W)^{\Gamma}$ is canonically equivalent to  $W_+ \wedge _{\Gamma} \mathbb{G}_{\Gamma}(*)$, and analogously for the simplicial set construction based on $W_{\dotr}$. Here $W_+$ denotes the complex $W $ with a disjoint fixed base point added.  The notation $\wedge _{\Gamma} $ denotes the smash product of $W_+$ with the spectrum $\mathbb{G}_{\Gamma}(*)$, followed by the orbit space construction using the diagonal action on the smash product.  
\end{Proposition}
\begin{Proof}
This is just a reinterpretation of Corollary \ref{fixedanalysis}.
    \end{Proof}

\begin{Definition} \label{HCGCXZ}
    The inclusion $T \colon {G} (\Gamma \backslash W_{\dotr}) \to 
{G}_{\Gamma} (W_{\dotr})$ will be referred to as the {\em asymptotic transfer}  associated to the $\Gamma$-free simplicial set  $W_{\dotr}$.  
	We can apply this construction to $W_{\dotr} = \mathcal{C}(X)_{\dotr}$, with the triangulation as in Proposition \ref{HFVH}, and obtain 
\[
{T} \colon  {G} (\Gamma \backslash \mathcal{C}(X)_{\dotr}) \longrightarrow {G}_{\Gamma} (\mathcal{C}(X)_{\dotr}).
\]	
Passing to parametrized theories over $\Gamma$, we extend ${T}$ to obtain 
\[
\mathbb{G}_{\Gamma} (\Gamma \backslash \mathcal{C}(X)_{\dotr}) \longrightarrow \mathbb{G}_{\Gamma} (\mathcal{C}(X)_{\dotr})
\]	
and the {\em parametrized asymptotic transfer} as the fixed point map
\[
P \colon  \mathbb{G}_{\Gamma} (\Gamma \backslash \mathcal{C}(X)_{\dotr})^{\Gamma} \longrightarrow \mathbb{G}_{\Gamma} (\mathcal{C}(X)_{\dotr})^{\Gamma}.
\]
\end{Definition}

Let $p \colon X \to \Gamma \backslash X$ be the orbit space projection.  
For a subcomplex $A$ of $\Gamma \backslash X$ there is a transfer map 
$T_A \colon {G} (\mathcal{C}(A)_{\dotr}) \to {G}_{\Gamma} (\mathcal{C}(p^{-1} A)_{\dotr})$.
If $\mathcal{U}$ is a finite covering of $\Gamma \backslash X$ by subcomplexes, we define
$p^{\ast} \mathcal{U} = \{ p^{-1} (U) \, \vert \, U \in \mathcal{U} \}$.
We have two homotopy colimits
\[
{G} (\mathcal{C} (\mathcal{U})) =
\hocolim{A \in {\mathcal{U}}} {G} (\mathcal{C}(A)_{\dotr}),
\]
\[
{G}_{\Gamma} (\mathcal{C}(p^{\ast} \mathcal{U})) = 
\hocolim{A \in {\mathcal{U}}} {G}_{\Gamma} (\mathcal{C}(p^{-1} A)_{\dotr})
\]
The maps $T_A$ induce a map
\[
{T}_{\mathcal{U}} \colon 
{G} (\mathcal{C}(\mathcal{U}))
\longrightarrow
{G}_{\Gamma} (\mathcal{C}(p^{\ast} \mathcal{U})).
\]
There is a direct extension of these constructions to the parametrized asymptotic version of the map
\[
P_{\mathcal{U}} \colon 
\mathbb{G}_{\Gamma} (\mathcal{C}(\mathcal{U}))^{\Gamma}
\longrightarrow
\mathbb{G}_{\Gamma} (\mathcal{C} (p^{\ast} \mathcal{U}))^{\Gamma}.
\]
From naturality of $\mathbb{G}_{\Gamma}$ with respect to inclusion maps, we have

\begin{Proposition} \label{FDSWX}
Whenever $\mathcal{U} \subset \mathcal{U}'$ are two finite coverings of $\Gamma \backslash \mathcal{C}(X)$, the square
\[
\xymatrix{
 {G} (\mathcal{C}(\mathcal{U}))
 \ar^-{{T}_{\mathcal{U}}}[r]  \ar_{}[d]
&{G}_{\Gamma} (\mathcal{C}(p^{\ast} \mathcal{U})) 
 \ar_{}[d] \\
{G} (\mathcal{C}(\mathcal{U}'))
 \ar^-{{T}_{\mathcal{U}'}}[r]
&{G}_{\Gamma} (\mathcal{C}(p^{\ast} \mathcal{U}'))
}
\]
commutes.  Similarly, the square
\[
\xymatrix{
 \mathbb{G}_{\Gamma} (\mathcal{C}(\mathcal{U}))^{\Gamma}
 \ar^-{P_{\mathcal{U}}}[r]  \ar_{}[d]
&\mathbb{G}_{\Gamma} (\mathcal{C}(p^{\ast} \mathcal{U}))^{\Gamma}
 \ar_{}[d] \\
\mathbb{G}_{\Gamma} (\mathcal{C}(p^{\ast} \mathcal{U}))^{\Gamma}
 \ar^-{P_{\mathcal{U}'}}[r]
&\mathbb{G}_{\Gamma} (\mathcal{C}(p^{\ast} \mathcal{U}'))^{\Gamma}
}
\]
commutes.  When $\mathcal{U}'$ is the one set covering $\{ \Gamma \backslash \mathcal{C}(X) \}$, the bottom maps in the diagrams are the transfer maps $T$ and $P$, respectively.
\end{Proposition}

The above discussion and Proposition \ref{FDSWX} apply equally well the to the simplicial set constructions based on $W_{\dotr}$.

\subsection{Excision properties} \label{ExcThms}

The goal of this section is to formulate several results about controlled $G$-theory in the generality and the language we use in this paper.  Throughout this section $A$ is a commutative regular Noetherian ring.  The modules are left $A$-modules. 

The results we need are a package of excision theorems analogous to those in controlled $K$-theory and a more recent fibrewise excision theorem that is unique to $G$-theory.  

Controlled $G$-theory is modeled on the Pedersen-Weibel controlled $K$-theory.  
The precise split local structure of the free modules used by Pedersen and Weibel is replaced by filtrations that apriori contain less information. 
However, the necessary control is regained by imposing various conditions on the filtrations such as leanness and insularity, cf. Definition \ref{splitleanins}.  Similarly for the new fibre-wise control situation.  

\begin{Example}
To illustrate the delocalized nature of the new kinds of objects, we give the following example.
    Let $\Phi$ be a lean insular object, for example a Pedersen-Weibel object.  Let $\overline{S}$ denote the complement of $S$ in $X$.  Choose any controlled set $E$.  We define a new filtration $\Phi'$ of the same module by $\Phi' (S) = \Phi (S \setminus E [\overline{S}])$.  This is a new lean insular object.  The identity homomorphism is an isomorphism between $\Phi$ and $\Phi'$ controlled by $E$.  However notice that when we use the standard bounded control structure, where $E$ is the $D$-enlargement of the diagonal for some number $D \ge 0$, $\Phi' (S) = 0$ for all $S$ with diameter smaller than $D$.  
\end{Example}

To recapitulate the point made earlier, a major technical issue in proving excision for controlled $G$-theories is inability to directly localize the objects supported on the whole control space to subspaces.  This is because even if $\Phi$ is lean and insular, a subconstruction $\Phi_U$ may not be lean and insular.  So, in terms of Definition \ref{adapted}, not every lean insular module is adapted. There are two approaches to remedy the situation.  Both have been used in the literature, however one is best suited for the fibred situation that's crucial to the results in this paper.  We will later contrast and relate two full subcategories of lean insular modules: the coarsely adapted modules and the subcategory of filtered modules that are isomorphic to adapted modules.  The latter was the subject of \cite{gCbG:11}.  The former however is the key result we use in this paper, so we start with the details of the theory based on the coarsely adapted modules.

\begin{Proposition}
   A filtered module isomorphic to an adapted module is coarsely adapted.
\end{Proposition}

\begin{proof}
    Let $\Phi$ be the module in question and $f \colon \Phi \to \Theta$ be an isomorphism where $\Theta$ is adapted, $f \Phi (U) \subset \Theta (E[U])$ as well as $f^{-1} \Theta (V) \subset \Phi (E[V])$ for all subsets $U \subset X$, $V \subset X$.  Now given a subset $U$, a desired enlargement of $\Phi_U$ can be defined as $\Psi = f^{-1} \Theta (E[U])$.  In this case, $\Phi_U \subset \Psi \subset \Phi_{E^2 [U]}$, in terms of the notation $E^2 [U]$ for the iterated coarse enlargement $E (E[U])$.
\end{proof}

This fact generalizes Lemma 3.3 from \cite{gCbG:11}.  The use of the isomorphism envelope of adapted modules in \cite{gCbG:11} is through this Lemma and the following two results that we demonstrate in the setting of coarsely adapted modules.  The reader may want to compare this to Theorems 2.22 and 3.10 from \cite{gCbG:11}. 

We know from Proposition 2.18 in \cite{gCbG:11} that lean insular modules form an exact category.  
\begin{Proposition}
    The subcategory of locally finite coarsely adapted $X$-filtered modules is closed under isomorphisms and exact extensions.  Therefore, this subcategory is an exact category.
\end{Proposition}

\begin{proof}
The proof is a translation of that of Theorem 2.22 from \cite{gCbG:11}.  Let
\[
\Phi' \xrightarrow{\ f \ } \Phi \xrightarrow{\ g \ } \Phi''
\]
be an exact sequence of locally finite lean insular modules and let $E$ be a controlled set that verifies both $f$ and $g$ as boundedly bicontrolled maps.  We also assume that both $\Phi'$ and $\Phi''$ are coarsely adapted. We
need to check that $\Phi$ is coarsely adapted.
Consider $U \subset X$, so an enlargement ${\Psi}''_{E[U]}$ of ${\Phi}''_{E[U]}$
is lean and insular. The induced epic
\[
g \colon \Phi(E^2 [U]) \cap g^{-1} {\Psi}''_{E[U]} \longrightarrow
{\Psi}''_{E[U]}
\]
extends to another epic
\[
g' \colon f {\Psi}'_{E^3 [U]} + \Phi(E^2 [U]) \cap g^{-1} {\Psi}''_{E[U]} \longrightarrow
{\Psi}''_{E[U]}
\]
with $\ker (g') = {\Psi}'_{E^3 [U]}$.
Since both ${\Psi}'_{E^3 [U]}$ and ${\Psi}''_{E[U]}$
are lean and insular,
then
\[
\Psi = f {\Psi}'_{E^3 [U]} + \Phi(E^2 [U]) \cap g^{-1} {\Psi}''_{E[U]}
\]
with the standard submodule filtration is lean and insular, and also has the property $\Phi_U \subset \Psi \subset \Phi_{E^4 [U]}$.
\end{proof}

Recall that a filtered module $\Phi$ is supported near a subset $S$ if $\Phi (E[S]) = \Phi (X)$ for some enlargement $E[S]$.  The category of coarsely adapted $X$-filtered modules has the following property, checked by translating the proof of Theorem 3.10 from \cite{gCbG:11}.  

\begin{Proposition}
    Suppose $\Phi (U) \to \Phi$ is the inclusion of a subconstruction, where $\Phi$ is coarsely adapted.  Let $\Psi$ be a submodule satisfying the property $\Phi_U \subset \Psi \subset \Phi_{E [U]}$, with the induced filtration.  Then the quotient module $\Phi/\Psi$ has a natural $X$-filtration and is coarsely adapted.  Moreover, $\Phi/\Psi$ is supported near the complement $X \setminus U$.
\end{Proposition}

Now we can recall the major fibrewise excision result from \cite{gCbG:19}.
The results are for the equivariant theory built from bounded $\Gamma$-actions on the fiber $Y$. 
The objects are modules equipped with $(X \times Y)$-filtrations $\Phi$ so that the induced $X$-filtration $\Phi_X$ is coarsely adapted and $\Phi$ is fibrewise coarsely weakly adapted.  
The paper \cite{gCbG:19} is written in terms of objects which are $(X \times Y)$-filtrations $\Phi$ so that the induced $X$-filtration $\Phi_X$ is just lean and insular and $\Phi$ itself is again fibrewise coarsely weakly adapted. So the difference is in the nature of $\Phi_X$ and mimics precisely the difference between $\mathcal{G}_{X}(\mathrm{pt},A)$ here and all lean insular $X$-filtered modules.  {\em Fibrewise excision} matters are transverse to this difference in $\Phi_X$ since all constructions and proofs are compatible with the additional structure.  This allows to translate the main theorem of \cite{gCbG:19} as follows. 




The following is a combination of Theorem 6.11 from \cite{gCbG:19} and Proposition 4.4.3 from \cite{gCbG:18}.

\begin{Theorem} \label{ExcGalltwo}
Given a single bounded action of $\Gamma$ on $Y$, the natural map 
\[
\delta \colon \hocolim{U \in \mathcal{U}} \mathbb{G}_{\Gamma} (Y, Y',U)^{\Gamma}  \longrightarrow
\mathbb{G}_{\Gamma} (Y,Y')^{\Gamma}
\]
is a weak equivalence.  One consequence of this is a familiar classical restatement: given a subset $V$ of $Y'$, the inclusion of pairs induces a weak equivalence
\[
\mathbb{G}_{\Gamma} (Y,Y')^{\Gamma} \simeq \mathbb{G}_{\Gamma} (Y-V,Y'-V)^{\Gamma}.
\]
\end{Theorem}

With the established $G$-theoretic interpretation of constructions and basic results from bounded $K$-theory in place, we have the key $G$-theoretic analogue of Theorem 2.5 from \cite{gCbG:04}.

Finally, we want to recall that excision properties suffice in \cite{gCbG:11} to show that, in the context of controlled $G$-theory built out of objects isomorphic to strict objects, the $G$-theoretic Novikov conjecture is true for finite asymptotic dimensional groups $\Gamma$ with finite $K(\Gamma,1)$.  That argument now applies directly to the $G$-theory based on locally finite coarsely adapted modules.  

 



\begin{Definition}
    In terms of the Waldhausen structure introduced in Proposition \ref{MM} and Notation \ref{MM2}, for any metric space $X$ with a $\Gamma$-action by isometries we have $\mathcal{G}_{X} (\mathrm{pt}, A)$, the equivariant theory of locally finite coarsely adapted $X$-filtered modules.  The space $X$ can be taken to be the group $\Gamma$ itself with the word metric associated to any chosen finite generating set.  We define the \textit{G}-theory of the group ring $G (A[\Gamma])$ to be the non-connective \textit{K}-theory spectrum $K(\mathcal{G}_{\Gamma} (\mathrm{pt},A)^{\Gamma})$.  
\end{Definition}

This allows to restate Corollary 5.27 of \cite{gCbG:04} for the $G$-theory based on locally finite coarsely adapted modules.

\begin{Theorem}
The controlled assembly map in $G$-theory
\[
\alpha \colon h^{\mathit{lf}} (E\Gamma, G (A)) \longrightarrow 
{G}_{\Gamma} (\mathrm{pt}, A)
\]
is an equivalence.
Therefore, the $G$-theoretic assembly map
\[
A_G = \alpha^{\Gamma} \colon B\Gamma_{+} \wedge G (A) \longrightarrow G
(A[\Gamma])
\]
is a split injection whenever $\Gamma$ is the fundamental group of a finite $K(\Gamma, 1)$ complex and has finite asymptotic dimension.  
\end{Theorem}

In addition, it was argued in \cite{gCbG:16} that $G$-theory of the group ring is close to $K$-theory enough for the Cartan map $\kappa \colon K (A[\Gamma]) \to G (A[\Gamma])$ to be an equivalence under specific conditions.  The only caveat is \cite{gCbG:16} was written for $G$-theory of locally finite lean insular $\Gamma$-filtered modules.  The subcategory of locally finite coarsely adapted $\Gamma$-filtered modules is sandwiched between that older theory and the $K$-theory.  Inspection of \cite{gCbG:16} shows that all results apply verbatim in the context of the present version of $G (A[\Gamma])$. 
The key construction in \cite{gCbG:16} is of a finite resolution by locally finite free equivariant Pedersen-Weibel modules of an equivariant locally finite lean insular $\Gamma$-filtered module.  The same construction gives a resolution of a module that is coarsely adapted.  The result is summarized as follows.

\begin{Theorem}
    The Cartan map $\kappa \colon K (A[\Gamma]) \to G (A[\Gamma])$ induced by the exact inclusion is an equivalence whenever $\Gamma$ is the fundamental group of a finite $K(\Gamma, 1)$ complex and has finite asymptotic dimension and $A$ is a regular Noetherian ring of finite global dimension.
\end{Theorem}

\section{Proof of the Main Theorem} 

As we explained in the introduction, the Main Theorem follows from the more general $G$-theory Isomorphism Theorem, which we poceed to prove in this section.

First we give a brief review of notation and a summary of required space inputs.  

All equivariant theories we need are $G$-theory spectra with fibre-wise control over products $\Gamma \times Y$ for various space entries $Y$.  When $Y$ is a metric space with a bounded action by $\Gamma$, we have the basic theory $\mathbb{G}_{\Gamma}(Y,Y',Z,A)$.  In this general notation, $Y'$ and $Z$ are mutually unrelated, arbitrary metric subsets of $Y$, and $A$ are arbitrary Noetherian ring coefficients. 

\begin{Notation}
    Additional properties of the ring $A$ become important when we apply our conclusions in $G$-theory to the $K$-theory assembly.  But in this part of the paper devoted entirely to $G$-theory, the coefficients act as a ``dummy variable'', so for simplicity we will omit $A$ from the notation.
\end{Notation}


\subsection{Homotopy computability for free simplicial coarse sets with $\Gamma$ action} \label{simplhomot}
Let $\Gamma$ denote a finitely generated group. We say a spectrum $X$ with $\Gamma$-action is {\em homotopy computable} if the natural map $X^{\Gamma} \rightarrow X^{h \Gamma}$ is an equivalence.  For any set $X$, recall that $\iota (X)$ is the object in $\mathfrak{E}$ whose underlying set is $X$ and is equipped with the extended metric which takes the value $+\infty$ for all pairs of distinct points. We will use a simpler notation $X_{\infty}$ for this metric.  $\Gamma_{\infty}$ is equipped with the left multiplication action by $\Gamma$. 
Recall that $\mathbb{G}_{\Gamma}(\Gamma_{\infty})$ is the $G$-theory spectrum  attached to the category with $\Gamma$-action $\mathfrak{G}_{\Gamma}( \Gamma _{\infty},A)$ defined to be the equivariant construction $\cat{C}_{\Gamma} $, where $\cat{C}$ is the  category with $\Gamma$-action $\mathcal{G}_{X}(\Gamma _{\infty})$, with $X = \Gamma$, equipped with word length metric. 

The way we would apply this is the following. 


\begin{Theorem} \label{CorFull}
Suppose $\Gamma$ is a finitely generated group with finite homological dimension, equipped with the action by left multiplication. Let $Y$ denote a locally finite ordered simplicial complex with proper free $\Gamma$-action, and let $Y_{\dotr}$ denote the corresponding simplicial object in the category of free $\Gamma$-sets  and proper $\Gamma$-equivariant maps.  Applying the functor $\iota$ from Section \ref{simpreal} we obtain a simplicial object in $\mathfrak{E}^m$, to which we may apply the functor $\mathbb{G}_{\Gamma}(-)$ to obtain a simplicial  $\Gamma$-spectrum.  
Then the $\Gamma$-spectrum $|\mathbb{G}_{\Gamma}(X_{\dotr})|$ is homotopy computable.  
\end{Theorem}
\begin{Proof} We observe that this result follows from the case of a single free orbit.  This is the case because the realization functor commutes with the functor $(-)^{h \Gamma}$ since the group $\Gamma$  has finite homological dimension, and this in turn means that we can compute the homotopy fixed point spectrum level-wise because the finite inverse limit construction commutes with the colimits defining geometric realization.   The case of a single free orbit follows from Proposition \ref{prodcriterion}, since
$$(\prod_{\gamma \in \Gamma} S)^{\Gamma} \cong (\prod_{\gamma \in \Gamma} S)^{h \Gamma}$$ for any spectrum $S$.   
\end{Proof}

\subsection{Maps derived from coverings}
Let $I$ denote the opposite category of the partially ordered set of non-empty subsets of $\{ 0,1 \}$, and let $I(n)$ denote the $n$-fold product $I^n$.  We will be considering  based space- or spectrum-valued functors from $I(n)$, and their homotopy colimits.  Let $I(n)^0 \subseteq I(n)$ denote the full subcategory on  all objects except the object 
$$  \underbrace{ \{ 0,1 \} \times \mathop{{\cdot \cdot \cdot} } \times \{ 0,1 \}}_{\mbox{\em{n}\hspace{.15cm}} \mbox{factors}} = (\{0,1\}, \ldots , \{ 0,1 \} )
$$
 We recall that given a space valued functor $F$ on $I(n)$, we can form its homotopy colimit as the realization of a simplicial space whose $k$-level space is 
$$  \coprod _{S_0 \supseteq S_1 \supseteq \cdots \supseteq S_k} F(S_0)
$$
with suitably defined face and degeneracy maps (see \cite{aBdK:72}). In particular, it is equipped with a map to the nerve of the category over which the homotopy colimit is taken, which is $I(n)$ in our situation. The nerve of $I(n)$ is a simplicial set which realizes to the $n$-cube $I^n$. 

\begin{Example} \label{mvexample}
    Let $X$ be a space, and for each $1 \leq i \leq n$ let $\{ X_0^i, X_1^i \}$ be a covering of $X$, i.e. $X = X_0^i \cup X_1^i$ for each $i$.  For any non-empty subset $S \subseteq \{ 0,1 \}$, we define $\theta _i (S) = \bigcap _{t \in S}(X^i_t)$.  Then there is a functor  $\Theta$ from $I(n)$ to spaces defined  by $\Theta (S_1, \ldots , S_n) = \bigcap _i \theta _i(S_i) $. Each $\Theta (S_1, \ldots , S_n)$ is equipped with an inclusion into $X$, and those inclusions are compatible with the inclusions among the sets  $\Theta (S_1, \ldots , S_n)$, thereby constructing a map $$
    \tau \colon \hocolim{} \Theta \longrightarrow X
    $$
    This map is an equivalence, by an argument identical to that given in 
    \cite{aBdK:72} to prove that the ``Mayer-Vietoris" blowup of a covering $\{U_1, \ldots, U_n \}$ of a space $X$ maps naturally by an equivalence to $X$.  
   
\end{Example}

For functors from $I(n)$ to based spaces, there is a based version of the homotopy colimit given by $$\hocolim{} F /\hocolim{} *$$  This construction is applied level-wise to give a spectrum level version of the homotopy colimit.  

\begin{Remark}
    {\em We will be constructing spectrum valued functors from coverings in the category $\mathfrak{E}$ by applying the spectrum valued functor $\mathbb{G}_{\Gamma}(-,A)$ to the members of the coverings.}
\end{Remark}

We now consider the category $\mathfrak{I}(n)$ of based space- or spectrum-valued functors on $I(n)$, with natural transformations as the morphisms. We have the functor $H$ which assigns to each functor on $I(n)$ its homotopy colimit.  We have a second functor $\iota$ with assigns to a functor $F$ its value at the object $( \{0,1\}, \{ 0,1 \}, \ldots , \{ 0,1 \}) $

We consider the based homotopy colimit construction for a functor $F$, and modify $F$ by collapsing all the values of $F$ on $I(n)^0$ to the basepoint.  It is clear by inspection that the result is the $n$-fold suspension of $F(\{0,1\}, \ldots , \{ 0,1 \}) $.  The result is a natural transformation $\eta$ from $H$ to $\Sigma ^n \iota$. 

\begin{Proposition} \label{suspension} Let $F$ be a functor from $I(n)$ to based spaces or spectra.  Suppose further that for all objects $x$ in $I(n)^0$, $F(x)$ is contractible.  Then $$\eta(F) \colon \hocolim{I(n)} F \longrightarrow \Sigma ^n F(\{ 0,1 \} , \ldots, \{ 0,1 \})$$ is an equivalence. 
\end{Proposition}
\begin{Proof}
    This is immediate from the based  or spectrum level version of the homotopy invariance of the homotopy colimit construction, proved for instance in \cite{aBdK:72}.
\end{Proof}

We will apply Proposition \ref{suspension} to the following situation.  Suppose that we have two spectrum valued functors $F$ and $G$, two spectra $X$ and $Y$, natural transformations from $F$ and $G$ to $F_X$ and $G_Y$, where $F_X$ (respectively $G_Y$) denotes the constant functor with value $X$ (respectively $Y$). Suppose further that we have a commutative diagram of functors on $I(n)$
$$ \begin{diagram}
    \node{F} \arrow{s} \arrow{e} \node{G} \arrow{s}  \\
    \node{F_X} \arrow{e,t}{\theta} \node{G_Y}
    \end{diagram}
$$
Note that the natural transformation $\theta$ amounts to a map from  $X$ to $Y$.  

\begin{Proposition} \label{verify}
    Let $F$, $G$, $F_X$, and $G_Y$ be as above.  Suppose that the maps on homotopy colimits 
    $$ \hocolim{} F  \rightarrow  \hocolim{}  F_X \cong X \mbox{ and } \hocolim{}  G \rightarrow \hocolim{}  G_Y \cong Y
    $$
    are equivalences.  Suppose further that the values of $F$ and $G$ are contractible for all objects of $I(n)^0$, and that 
    $$\eta \colon F(\{0,1\}, \ldots , \{ 0,1 \}) \rightarrow  G(\{0,1\}, \ldots , \{ 0,1 \})
    $$
    is an equivalence.  
   Then the map $X \rightarrow Y $ encoded by $\theta$ is an equivalence. 
   \end{Proposition}
   \begin{Proof}
       Follows immediately from Proposition \ref{suspension}.
   \end{Proof}

\subsection{Splitting fixed points from homotopy fixed point sets } \label{SPL}

We will need a technical lemma that gives necessary conditions which provides a splitting of fixed points of a group action from a homotopy fixed point set.   

\begin{Proposition}\label{SFP}
    Let $S$ be a spectrum with $\Gamma$-action, where $\Gamma$ is a group.  Suppose there exists an equivariant map $f \colon S \rightarrow E$, where $E$ is another $\Gamma$-spectrum, with the following properties: (a) the natural map $\varepsilon \colon E^{\Gamma} \rightarrow E^{h \Gamma}$ is an equivalence, and (b) there is a splitting $s \colon E^{\Gamma} \rightarrow S^{\Gamma}$ of $f^{\Gamma}$.  Then the natural inclusion $\varepsilon \colon S^{\Gamma} \rightarrow S^{h \Gamma}$ is a split injection of spectra, i.e. there is a map $s^* \colon S^{h \Gamma} \rightarrow S^{\Gamma}$ so that the composite 
    $$ S^{\Gamma} \stackrel{\varepsilon}{\longrightarrow } S^{h \Gamma} \stackrel{s^*}{\longrightarrow} S^{\Gamma}
    $$
    is homotopic to the identity.  
\end{Proposition}
\begin{Proof}
    Consider the diagram
    $$\begin{diagram}\node{S^{\Gamma} } \arrow{s,t}{f^{\Gamma}} \arrow{e,t}{\varepsilon _S} \node{S^{h \Gamma}} \arrow{s,b}{f^{h \Gamma}} \\
    \node{E^{\Gamma}} \arrow{s,t}{s}  \arrow{e,t}{\varepsilon _E} \node{E^{h \Gamma}}  \\
    \node{S^{\Gamma}}
    \end{diagram}
    $$
    The composite $ s \cdot \varepsilon _{E}^{-1} \cdot f^{h \Gamma}$, where $\varepsilon _E^{-1} $ 
    denotes a homotopy inverse to $\varepsilon _{E}$, is the required $s^*$. 
\end{Proof}
\subsection{The proof} \label{proofend}

We plan to apply Proposition \ref{SFP},  and therefore need to supply the ingredients $S$ and $E$.  To construct these $\Gamma$-spectra, we will need to define some spaces attached to the classifying space of the group $\Gamma$, which can by our hypothesis be chosen to be a finite complex $K$. Let $n$ be sufficiently large so that $K$ can be embedded in the interior of an closed $n$-ball   $D^n \subseteq S^n$.  We assume that $D^n $ is triangulated, and that $K$ is embedded as a subcomplex of that triangulation.  We  further suppose  that $N^n$ is a closed regular neighborhood of $K$, that it too is a subcomplex of $D^n$, and that we have a ball $B^n \subseteq N^n$ which is also a subcomplex of $D^n$. We suppose further (after subdivision if necessary) that the triangulation of $D^n$ is an ordered  simplicial complex, and that we can therefore construct the associated simplicial sets, which we denote by $D^n_{\dotr}$, $K_{\dotr}$,  $N_{\dotr}$, and $B_{\dotr}^n$. 
We will regard all of  these as  simplicial objects in the subcategory of discrete coarse spaces.  We also let $A$ be a commutative Noetherian ring.  
For each of these simplicial complexes, we construct the simplicial cone as in Section \ref{conedef}, obtaining metric spaces $CD^n$, $C\partial D^n$, $CK$, $CN$, and $B^n$, as well as the corresponding simplicial sets 
$CD^n_{\dotr}$, $C\partial D^n_{\dotr}$, $CK_{\dotr}$, $CN_{\dotr}$, and $B^n_{\dotr}$.  

Regarding $CD^n_{\dotr}$ as a simplicial set with trivial $\Gamma$-action, we now define 
$$S = \mathbb{G}_{\Gamma}(CD^n_{\dotr}, C \partial D^n_{\dotr};A).
$$
The homotopy type of this simplicial $\Gamma$-spectrum  is 
$$ (\Sigma D^n_{\dotr}, \Sigma \partial D^n_{\dotr}) \wedge \mathbb{G}_{\Gamma}(*, A). 
$$
In particular, its fixed point spectrum is 
$$ (\Sigma D^n_{\dotr}, \Sigma \partial D^n_{\dotr}) \wedge  \mathbb{G}(*,A[\Gamma]) \cong \Sigma ^{n+1}K(A[\Gamma ]).
$$

Next, we will need to define the homotopy computable $\Gamma$-spectrum $E$.  Let $Y$ be the universal cover of $N^n$.  It inherits a triangulation from $N$, so that the covering projection $p \colon Y \rightarrow N^n $ is a simplicial map. The subset $\partial Y = \pi ^{-1}(\partial N^n)$ is a subcomplex, and we have the corresponding simplicial sets $Y_{\dotr}$ and $\partial Y_{\dotr}$ with free $\Gamma$-action.  We now define $E$ to be 
$$  |\mathbb{G}_{\Gamma}CY_{\dotr}, C \partial Y_{\dotr})|
$$
The simplicial complexes $CY$ and $C\partial Y$ are clearly locally finite ordered simplicial complexes, and are both  equipped with a proper free  $\Gamma$-actions, so Theorem \ref{CorFull} shows  that $\mathbb{G}_{\Gamma} (CY_{\dotr})$ and $\mathbb{G}_{\Gamma}(C\partial Y_{\dotr})$ are both homotopy computable.  The fibration sequence of $\Gamma$-spectra 
$$|\mathbb{G}_{\Gamma}(C\partial Y_{\dotr}) |\rightarrow |\mathbb{G}_{\Gamma} (CY_{\dotr})| \rightarrow  |\mathbb{G}_{\Gamma}(CY_{\dotr}, C \partial Y_{\dotr})|=E
$$
 demonstrates that $E$ is homotopy computable.  

 \begin{Remark}
     From Proposition \ref{UFSV}, it follows that $E^{\Gamma}$ is equivariantly equivalent to the smash product of $\mathbb{G}_{\Gamma} (*) \wedge _{\Gamma} (CY_{\dotr}, C\partial Y_{\dotr})$, which by an equivariant $S$-duality (for free complexes with a finite number of $\Gamma$-cells)  argument using \cite{aR:80} is equivalent to a suspension of the homotopy fixed point set $\mathbb{G}_{\Gamma}(*)^{h \Gamma}$. So, we will produce a splitting of the inclusion $\mathbb{G}_{\Gamma}(*)^{\Gamma} \hookrightarrow \mathbb{G}_{\Gamma}(*)^{h \Gamma}$.  This is useful intuition, but is not required for the proof, because of the flexibility that Proposition \ref{SFP} affords us.   
 \end{Remark}

We proceed to define the maps $f$ and $s$.

We assume that $CN$ is a subcomplex of $CD$, so there is a quotient induced by factoring out the complement $\complement CN$, which also has a simplicial subcomplex structure.  This results in a map 
\[
\phi \colon \mathbb{G}_{\Gamma}(CD^n_{\dotr}, C \partial D^n_{\dotr}) \to \mathbb{G}_{\Gamma}(CN_{\dotr}, C \partial N_{\dotr}).
\]

This is an equivariant map as the action on all spaces acting as fibers here is trivial.
The parametrized transfer from Definition \ref{HCGCXZ} is
\[
\mathcal{P} \colon \mathbb{G}_{\Gamma}(CN_{\dotr}, C \partial N_{\dotr} ) \longrightarrow \mathbb{G}_{\Gamma}(CY_{\dotr}, C \partial Y_{\dotr} ).
\]
We define $f$ as the composition
\[
\mathcal{P} \circ \phi \colon \mathbb{G}_{\Gamma}(CD_{\dotr}, C \partial D_{\dotr} ) \longrightarrow \mathbb{G}_{\Gamma}(CY_{\dotr}, C \partial Y_{\dotr} ),
\]
so $f^{\Gamma}$ is ${P} \circ \phi^{\Gamma}$.


We continue to define $s$.  In section \ref{simpreal}, we presented an extended metric on the geometric realization $CY$ and a coarse map that now induces an equivariant natural transformation 
\[
\ell \colon \mathbb{G}_{\Gamma}(CY_{\dotr}, C \partial Y_{\dotr}) \longrightarrow \mathbb{G}_{\Gamma}(CY, C \partial Y).
\]

\begin{Notation}
    We will make frequent use of imposed left-bounded metrics in the proof; it will be easier to use a compact notation $||  \underline{\phantom{W}} ||$ for this construction, in place of $(\underline{\phantom{W}})^{bdd}$.
\end{Notation}

The change to a left-bounded metric $b \colon Y \to || Y ||$  is a proper coarse map giving a proper coarse map $CY \to || CY ||$, according to Proposition \ref{HJDSEO}. This induces 
\[
\alpha \colon \mathbb{G}_{\Gamma}(CY, C \partial Y )^{\Gamma} \longrightarrow \mathbb{G}_{\Gamma}( || CY ||,  || C \partial Y || )^{\Gamma}.
\]

Recall that $\mathbb{G}_{\Gamma}(( || CY ||,  || C \partial Y || ))$ is the $K$-theory of Schlichting's exact quotient category $\mathcal{G}_{\Gamma}(|| CY || )/\mathcal{G}_{\Gamma}(|| CY || )_{<C \partial Y}$.  
From Proposition \ref{etacompatibility} we see that the equivariant natural transformation $\eta \colon \mathbb{G}_{\Gamma}(|| CY ||) \to \mathbb{G}_{\Gamma}(|| CY ||^0)$ respects the filtering subcategories and so gives a natural transformation
\[
\beta \colon \mathbb{G}_{\Gamma}(|| CY ||, || \partial CY || )^{\Gamma} \longrightarrow \mathbb{G}_{\Gamma}(|| CY ||^0, || \partial CY ||^0 )^{\Gamma}.
\]

Let $\overline{B}$ be a homeomorphic lift of $B$ in $|| Y ||$.  We can in fact assume that for a sufficiently small ball $B$ the projection $\overline{B} \to B$ is an isometry.  
Then we have a further quotient $\mathcal{G}_{\Gamma}(|| CY ||^0, || \partial CY ||^0 )/\mathcal{G}_{\Gamma}(|| CY ||^0, || \partial CY ||^0 )_{<\complement C\overline{B}}$ by the filtering subcategory of objects fiber-wise supported near the complement of the closed subset $C\overline{B}$ of $CY$.  The fibred excision theorem \ref{ExcGalltwo} gives an equivariant equivalence of the latter spectrum to $\mathbb{G}_{\Gamma}(C\overline{B}, C \partial \overline{B} )$. 
This quotient map induces a \textit{restriction} map on $K$-theory 
\[
\gamma \colon \mathbb{G}_{\Gamma}(|| CY ||^0, || \partial CY ||^0 )^{\Gamma} \longrightarrow \mathbb{G}_{\Gamma}(C\overline{B}, C \partial \overline{B} )^{\Gamma}.
\]

We define $s = \gamma \circ \beta \circ \alpha$.  

The proof is finished by showing that $s \circ f^{\Gamma}$ is an equivalence.
We do that by employing Proposition \ref{verify}.  This requires us to identify the functors $F$ and $G$ with the property that there is a natural transformation $F \to G$ and with the property that all values of $F$ and $G$ on $I(n)^0$ are contractible, while the map $\eta$ is an equivalence.  

First let us create some geometric features in the topological balls $D$ and $\overline{B}$.  Without loss of generality, we can assume that $D$ and $B$ are concentric closed disks centered at $0$ in the $n$-dimensional real vector space $E$ with the max metric.

We identify the following subsets of $E$,
\begin{align}
&E^i_{0} = \{ (x_1, \ldots, x_{n}) \mid x_i \le 0 \}, \notag \\
&E^i_{1} = \{ (x_1, \ldots, x_{n}) \mid x_i \ge 0 \}, \notag
\end{align}
and apply the constructions from Example \ref{mvexample} to this collection of coverings of $E$.
There are also associated coverings of $D$ given by specifying for all $1 \le i \le n$
\begin{align}
&D^i_{0} = E^i_{0} \cap D, \ D^i_{1} = E^i_{1} \cap D,  \notag \\
&B^i_{0} = B^i_{0} \cap B, \ B^i_{1} = B^i_{1} \cap B. \notag 
\end{align}

The covering of $B$ by the subsets $B_{\ast}^i$ (where $\ast$ is equal either 0 or 1) can be lifted to $\overline{B}$ via the isometry $\overline{B} \to B$ resulting in the corresponding coverings by $\overline{B}^i_{\ast}$ of $\overline{B}$ with isometries $\overline{B}^i_{\ast} \to {B}^i_{\ast}$.  So now we have from Example \ref{mvexample} the two functors on the category  $I(n)$: $\Theta_F$ applied to the subsets $CD_{\ast}^i$ of $CD$ and $\Theta_G$ applied to the subsets $C\overline{B}_{\ast}^i$ of $C\overline{B}$.  

We are ready to define the functors $F$ and $G$ on $I(n)$.  They are given by 
\begin{align}
&F(s) = \mathbb{G}_{\Gamma}(CD_{\dotr},C \partial D_{\dotr},\Theta_F (s))^{\Gamma},  \notag \\
&G(s) = \mathbb{G}_{\Gamma}( C\overline{B}, C \partial \overline{B}, \Theta_G (s))^{\Gamma}. \notag 
\end{align}
We will define a natural transformation from $F$ to $G$ as a sequence of transformations in several stages, each stage using various spaces $X$ and the associated coverings $\{ X^i_0, X^i_1 \}$.  
It's easy to generalize the notation above: given a space $X$, we will be able to describe its boundary $\partial X$ as a subset of $X$, and $\Theta (s)$ are described as subsets of $X$.  It then makes sense to define
$F_X (s)$ as either $\mathbb{G}_{\Gamma}(X, \partial X, \Theta (s))^{\Gamma}$ if $X$ is an extended metric space or $| \mathbb{G}_{\Gamma}(X, \partial X, \Theta (s))^{\Gamma} |$ if $X$ is a proper simplicial set.

So now $F = F_{CD_{\dotr}}$ and $G= G_{C\overline{B}}$.

First, let $X=CN_{\dotr}$. The boundary $\partial N_{\dotr}$ is the usual boundary of the subcomplex $N_{\dotr}$ of $D_{\dotr}$.  The coverings of $N$ are by $N^i_{\ast} = N \cap D^i_{\ast}$.  There are corresponding coverings of $CN$ by $CN^i_{\ast} = C(N \cap D^i_{\ast})$ giving the functor $\Theta_{CN_{\dotr}}$.
Now the general construction gives a functor $F_{TN_{\dotr}}$.  

Consider $\complement N$, the complement of the interior of $N$ in $D$.  
We denote by $\complement (CN)$ the complement of the interior of $CN$ in ${CD}$.  Clearly $\complement (CN)_{\dotr}$ is a subcomplex of $CD_{\dotr}$.  Now for each $s$ the exact quotient functor on the level of fibrewise controlled categories induces a map
\[
\mathbb{G}_{\Gamma}(CD_{\dotr}, C \partial D_{\dotr}, \Theta_{CD} (s))^{\Gamma}
\longrightarrow
\mathbb{G}_{\Gamma}(CD_{\dotr}, \complement (CN)_{\dotr}, \Theta_{CD} (s) )^{\Gamma}.
\]
Using Theorem \ref{ExcGalltwo}, we simplify the target to the equivalent spectrum $F_{TN_{\dotr}}$:
\[
\mathbb{G}_{\Gamma}(CD_{\dotr}, \complement (CN)_{\dotr}, \Theta_{CD} (s) )^{\Gamma}
\simeq
\mathbb{G}_{\Gamma}(CN_{\dotr}, C \partial N_{\dotr}, \Theta_{CN} (s))^{\Gamma}.
\]

\begin{center}
    \includegraphics[height=6cm]{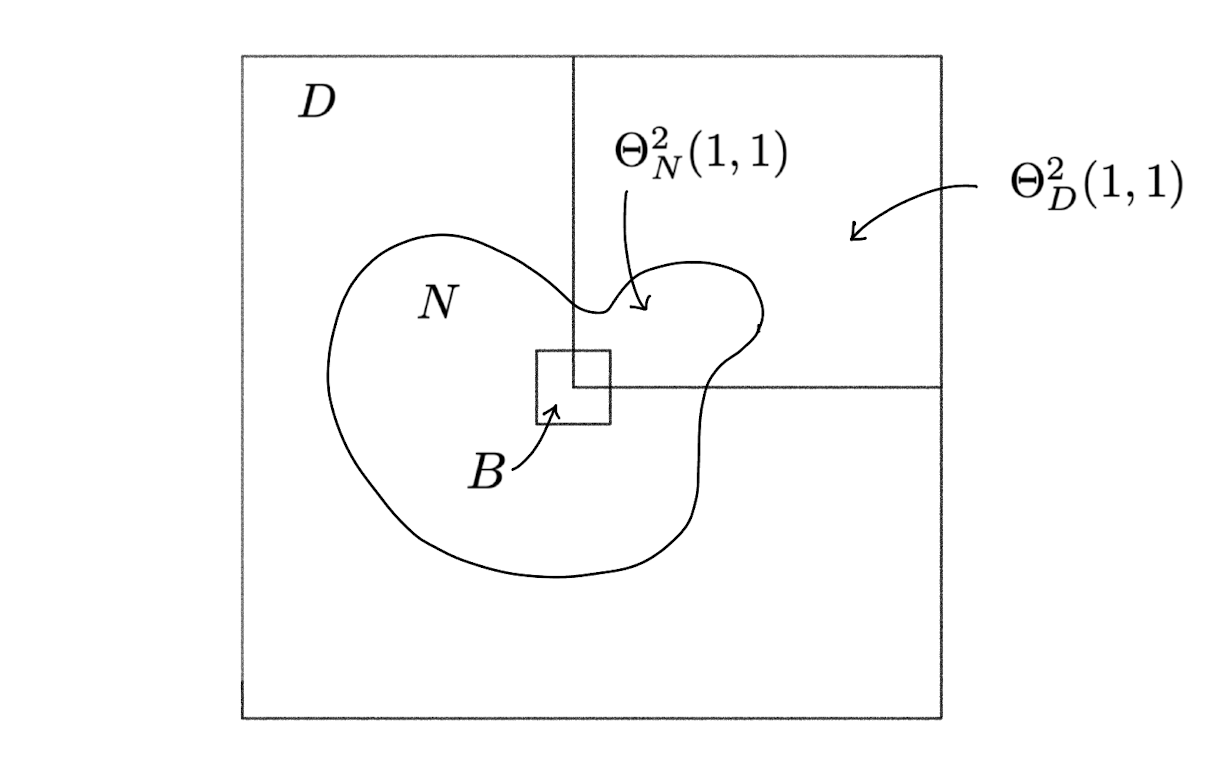}\\
    Figure 1. An illustration of notation for assorted subsets of $D$ in dimension $n=2$.
\end{center}

We have so far constructed a natural transformation $F_{CD_{\dotr}} \to F_{CN_{\dotr}}$.

The transformation $F_{CN_{\dotr}} \to F_{CY_{\dotr}}$ is the parameterized transfer applied to each set in the covering of $N$, using the restriction of the projection map $p \colon CY_{\dotr} \to CN_{\dotr}$.  The sets $CN^i_{\ast}$ lift to full preimages under $p$, denoted $CY^i_{\ast}$.  

The transformation $F_{CY_{\dotr}} \to F_{CY}$ is clear.

The introduction of left-bounded metric in $CY$ is compatible with doing the same in $CY^i_{\ast}$, in the sense that $|| CY || ^i_{\ast} = || CY^i_{\ast} ||$.  This defines a natural transformation $F_{CY} \to F_{|| CY || } \to F_{|| CY ||^0 }$.

Consider the complement $\complement (C\overline{B})$  of the interior of $C\overline{B}$ in $CY$.  
For each $s$ in $I(n)$, the exact quotient functor induces a map
\[
\mathbb{G}_{\Gamma}({|| CY ||^0}, {|| C \partial Y ||^0}, \Theta_{|| CY ||}(s) )^{\Gamma}
\longrightarrow
\mathbb{G}_{\Gamma}( {|| CY ||^0}, {\complement C\overline{B}},  \Theta_{|| CY ||} (s) )^{\Gamma}.
\]
Using Theorem \ref{ExcGalltwo}, we again simplify the target to an equivalent spectrum $F_{T\overline{B}}$ by passing to the equivariant spectrum with respect to the trivial action and excising out the complement of the closure of $C\overline{B}$ in $CY$:
\[
\mathbb{G}_{\Gamma}( {|| CY ||^0}, {\complement C\overline{B}},  \Theta_{|| CY ||} (s)  )^{\Gamma}
\simeq
\mathbb{G}_{\Gamma}(C\overline{B},  C \partial \overline{B}, \Theta_{C\overline{B}} (s) )^{\Gamma}.
\]

We have constructed a sequence of natural transformations
\[
F = F_{CD_{\dotr}} \to F_{CN_{\dotr}} \to F_{CY_{\dotr}} \to F_{CY} \to F_{|| CY ||^0} \to F_{C\overline{B}} = G.
\]
Note, we don't intend to make any claim about the values of intermediate functors on $s$ in $I(n)$.  It remains to verify that only the values of $F$ and $G$ on $I(n)^0$ are contractible.  
This is clear for $F$ which is a locally finite homology theory.

\begin{Lemma}\label{conesuspension}
    $\mathbb{G}_{\Gamma}(C\overline{B},  C \partial \overline{B}, \Theta_{C\overline{B}} (s) )^{\Gamma}$ is contractible for all $s$ in $I(n)^0$.
\end{Lemma}

\begin{proof}
For the purpose of this computation, we can identify $\overline{B}$ with the standard metric ball $D$ centered at $0$ in the $n$-dimensional Euclidean space $E$.  The cone $C(D)$ can be viewed as a metric subspace of $C(E)$. The boundary $C (\partial D)$ consists of points $(t, z)$ where $z$ is an $n$-vector of length $f(t)$.  

The subsets
$E^i_{-1} = \{ (x_1, \ldots, x_{n}) \mid x_i \le 0 \}$, $E^i_{1} = \{ (x_1, \ldots, x_{n}) \mid x_i \ge 0 \}$ give a two set covering of $E$ for each $1 \le i \le n$.
Consider the covering of $E$ by the subsets $\Theta (S_1, \ldots , S_n)$ given in Example \ref{mvexample}, except it will be convenient to use the subscript $-1$ in this proof in place of $0$ to simplify the formulas.
As long as at least one $S_i \ne \{ -1, 1 \}$, we want to claim that $\mathbb{G} (\Theta (S_1, \ldots , S_n))$ is contractible.  This will follow from the Additivity Theorem by showing that the category $\mathcal{G} (\Theta (S_1, \ldots , S_n))$ is flasque.  So choose one specific instance of $S_i = \{ m \}$ where $m$ is either $-1$ or $1$.
We use the shift functor given by $\Phi (m, k) (S) = \Phi (S- m \overline{k})$, where $\overline{k}$ is the vector $(0,...,k,...,0)$ with the natural number $k$ in the $i$-th position.  Then we have an exact endofunctor 
$$T^i_m (\Phi) = \Phi(m,1) \oplus \Phi(m,2) \oplus ... $$ 
of $\mathcal{G} (\Theta (S_1, \ldots , S_n))$.
Clearly, $T^i_m \oplus 1$ is naturally equivalent to $T^i_m$, which shows that $\mathcal{G} (\Theta (S_1, \ldots , S_n))$ is flasque.

The same proof can be used to show that $\mathcal{G} (E)_{< \Theta (S_1, \ldots , S_n)}$ is flasque by applying the construction of $T^i_m$ to objects supported in a controlled neighborhood of $\Theta (S_1, \ldots , S_n)$ rather than $\Theta (S_1, \ldots , S_n)$ itself. 
Further, the same proof can be applied over a cone to show that $\mathcal{G} (CE)_{< C\Theta (S_1, \ldots , S_n)}$ is flasque.
The same proof works for the quotient category $\mathcal{G} (CE, \complement CD)_{< C\Theta (S_1, \ldots , S_n)}$ which has the same objects, but new quotient category morphisms.  The key is that the natural transformation that was an equivalence between $T^i_m \oplus 1 \simeq T^i_m$ is still an equivalence.  

There are associated two set coverings of $D$ given by 
$D^i_{-1} = E^i_{-1} \cap D$ and $D^i_{1} = E^i_{1} \cap D$ for all $1 \le i \le n$.  
Using excision theorem \ref{ExcGalltwo} we interpret each spectrum $\mathbb{G} (CD, C \partial D)_{< \Theta_{CD} (s)}$ as exactly $\mathbb{G} (CE, \complement CD)_{< C\Theta (S_1, \ldots , S_n)}$, so these spectra are contractible as well.  Further, what we've done can be done fibrewise, so $\mathbb{G}_{\Gamma} (CD, C \partial D)_{< \Theta_{CD} (s)}^{\Gamma}$ are flasque as well.
\end{proof} 

Finally, we consider the composition
\[
\eta \colon \mathbb{G}_{\Gamma}(CD_{\dotr} ,  C \partial D_{\dotr}, \Theta_{CD} (s) )^{\Gamma} \longrightarrow \mathbb{G}_{\Gamma}(C\overline{B},  C \partial \overline{B}, \Theta_{C\overline{B}} (s) )^{\Gamma}
\]
for the single value of $s = (\{0,1\}, \ldots , \{ 0,1 \})$.  In this case, we have the center $\overline{0}$ of $\overline{B}$ which projects to the center $0$ of $D$, so the element of the cone $C\overline{B}$ through $\overline{0}$ projects isometrically onto the element of the cylinder $CD$ through $0$ under $p$.  The inverse isometry $C(0) \to C(\overline{0})$ induces the map that coincides with the restriction of the construction in the proof to the category of modules supported on $C(0)$,
\[
\eta \colon \mathbb{G}_{\Gamma}(C(0))^{\Gamma} \longrightarrow \mathbb{G}_{\Gamma}(C(\overline{0}))^{\Gamma}.
\]
This last map is of course an equivalence.
There is in fact a commutative square
$$\begin{diagram}
\node{\mathbb{G}_{\Gamma}(CD_{\dotr},  C \partial D_{\dotr}, \Theta_{CD} (s) )^{\Gamma}} \arrow{s,t}{\pi_{CD}} \arrow{e,t}{\eta} \node{\mathbb{G}_{\Gamma}(C\overline{B},  C \partial \overline{B}, \Theta_{C\overline{B}} (s) )^{\Gamma}} \arrow{s,b}{\pi_{C\overline{B}}} \\
    \node{\mathbb{G}_{\Gamma}(C(0))^{\Gamma}}  \arrow{e,t}{\eta} \node{\mathbb{G}_{\Gamma}(C(\overline{0}))^{\Gamma}}  
    \end{diagram}
    $$
where the vertical maps are induced by the projections $CD \to C(0)$ and $C\overline{B} \to C(\overline{0})$.  The vertical maps are equivalences because all objects in this case are supported near $C(0)$ and $C(\overline{0})$.  This verifies that that the top map $\eta$ is an equivalence and completes the proof.




\appendix

\section{Group actions on Waldhausen categories}\label{groupactions}

We collect some technical results we require concerning the properties of the $K$-theory spectra we work with. Throughout we will deal with Waldhausen categories in which the weak equivalences consist exactly of the isomorphisms, and therefore are only dependent on the underlying category with cofibrations, as defined in \cite{fW:83}. We'll call such Waldhausen categories {\em special} 

 Suppose  that we have an action of a group $\Gamma$ on a special Waldhausen category $\cat{C}$, in the sense that we are given a $\Gamma$-action on the category $\cat{C}$, which carries cofibrations to cofibrations, and which preserves the initial object.  There are then  evident actions of $\Gamma$ on the iterated $wS_{\dotr}$ constructions, and therefore we obtain a $\Gamma$-action on the $K$-theory spectrum of $K\cat{C}$.  We therefore obtain the fixed point spectra $(K\cat{C})^{\Gamma ^{\prime} }$ for any subgroup $\Gamma ^{\prime} \subseteq \Gamma$. The purpose of this appendix is twofold.  We need to (1) identify conditions where the $\Gamma$ fixed points on the spectrum $K \underline{C}$ are equivalent to the $K$-theory spectrum $K( \underline{C}^{\Gamma})$ and (2) identify conditions under which inclusion of  categories with $\Gamma$-action  induces equivalences of equivariant spectra.  Propositions \ref{requirements} and \ref{otherrequirements} provide the conditions for (1), while Propositions \ref{noneqcriterion} and \ref{eqcriterion} will provide the right conditions for (2). 
 

\begin{Proposition}\label{requirements}  Let $(\cat{C}, co(\cat{C}),*)$ be a category with cofibrations, and suppose that a group $\Gamma$ acts on $\cat{C}$ in such a way that it preserves $co(\cat{C})$ and fixes the initial object $*$. Let $\cat{C}^{\Gamma}$ denote the fixed point subcategory, and let $co(\cat{C})^{\Gamma}$ denote the fixed points of the $\Gamma$-action on $co(\cat{C})$.       Suppose that the following two conditions hold.  
\begin{enumerate}
    \item{Any diagram in $\cat{C}^{\Gamma}$  of the form $C \longleftarrow A \stackrel{i}{\longhookrightarrow} B$, where $i$ is in $co(\cat{C})^{\Gamma}$, admits a pushout in $\cat{C}^{\Gamma}$.  }
    \item{For any diagram $C \longleftarrow A \stackrel{i}{\longhookrightarrow} B$ in $\cat{C}^{\Gamma}$, with $i \in co(\cat{C})^{\Gamma}$, the natural map from $C $ to the pushout of the diagram, which we denote by $C \cup _A B$, is also a cofibration.  }
\end{enumerate}
Then the triple $(\cat{C}^{\Gamma}, co(\cat{C})^{\Gamma},*)$ is a category with cofibrations.  We will define $co(\cat{C}^{\Gamma}) = co(\cat{C})^{\Gamma}$. 
\end{Proposition}
\begin{Proof}
    This is simply a statement of the definition of categories with cofibrations in \cite{fW:83}.  Conditions $cof1$ and $cof2$ are essentially trivial, and (1) and (2) above are  statements of the requirements in $cof3$. 
\end{Proof}

While we clearly need this statement, what we really need is a statement that identifies the $K$-theory spectrum  $K(\cat{C}^{\Gamma}, co(\cat{C}^{\Gamma}), *)$ with the fixed point spectrum $K(\cat{C}, co(\cat{C}), *)^{\Gamma}  $.  The requirements of Proposition \ref{requirements} are not sufficient to guarantee this.  The reason is that 
the definition of the $wS_{\dotr}$ construction requires that  certain diagrams of the form 
\begin{equation}\label{pushout}
\begin{diagram} \node{A} \arrow{e,t,J}{i} \arrow{s} \node{B} \arrow{s} \\
\node{*} \arrow{e} \node{C}
\end{diagram} 
\end{equation}
are pushouts. In  $wS_{\dotr}(\cat{C})^{\Gamma}$, we then find diagrams of the form (\ref{pushout}) which are fixed under $\Gamma$. Such diagrams are diagrams in $\cat{C}^{\Gamma}$   which are pushouts {\em in $\cat{C}$}. The requirement that a diagram is a pushout in $\cat{C}$ is more stringent than the requirement that it be a pushout in $\cat{C}^{\Gamma}$, so that what we have is an inclusion 
$$ K(\cat{C}, co(\cat{C}), *) ^{\Gamma} \hookrightarrow K(\cat{C}^{\Gamma}, co(\cat{C}^{\Gamma}), *)
$$ which need not be surjective.  The following result gives the condition that is required.  
\begin{Proposition}\label{otherrequirements}
    Let $\cat{C}$ be a category with cofibrations, and suppose that a group $\Gamma$ acts on $\cat{C}$, preserving the set  $co(\cat{C})$.  Suppose further that the following two hypotheses hold. 
\begin{itemize}
\item{The hypotheses of Proposition \ref{requirements} hold for the $\Gamma ^{\prime}$-action on $\cat{C}$ for all subgroups $\Gamma ^{\prime} \subseteq \Gamma$. }
    \item{For any subgroup $\Gamma ^{\prime} \subseteq \Gamma$ and any pushout diagram $\mathcal{P}$ of the form 
    $$
\begin{diagram} \node{A} \arrow{e,t,J}{i} \arrow{s} \node{B} \arrow{s} \\
\node{*} \arrow{e} \node{C}
\end{diagram} 
$$
in $\cat{C}^{\Gamma ^{\prime}}$, with  $i \in co(\cat{C})^{\Gamma ^{\prime}}$, $\mathcal{P}$ is also a pushout diagram in $\cat{C}$.  }
\end{itemize}
Then for every $\Gamma ^{\prime} \subseteq \Gamma $, the triple $(\cat{C}^{\Gamma ^{\prime}}, co(\cat{C})^{\Gamma ^{\prime}},*)$ is a category with cofibrations.  Moreover,  if for each $\Gamma ^{\prime}$, we set $w(\cat{C}^{\Gamma^{\prime}}) $ equal to the isomorphisms in
$\cat{C}^{\Gamma ^{\prime}}$, the inclusion 
$$ I \colon K(\cat{C}, co(\cat{C}), *)^{\Gamma ^{\prime}} \hookrightarrow K(\cat{C}^{\Gamma ^{\prime}}, co(\cat{C}^{{\Gamma} ^{\prime}}), *)
$$
is an isomorphism of spectra.  

\end{Proposition}
\begin{Proof}
    We provide the proof for the first delooping $wS_{\dotr}$.  The iterated situation for the higher deloopings works identically. Let $\cat{C}$ denote a category with cofibrations.  Recall that $wS_k\cat{C}$ is the nerve of a category whose objects are diagrams $\{ A_{ij} \}_{0 \leq i \leq j \leq k}$, with $A_{ij} \in \cat{C}$,  equipped with with morphisms $f_{ij}^{i^{\prime}j^{\prime}}$ in $\cat{C}$ from $A_{ij} $ to $A_{i^{\prime}j^{\prime}}$ whenever $i \leq i^{\prime}$ and $j \leq j^{\prime}$, satisfying the following two conditions. 
    \begin{enumerate}
    \item{The morphisms $f_{ij}^{ij^{\prime}} $ are all cofibrations. }
    \item{The diagrams $$ \begin{diagram}
    \node{A_{ij}}\arrow{s,t}{f_{ij}^{jj}} \arrow{e,t}{f_{ij}^{ij^{\prime}}} \node{A_{ij^{\prime}} } \arrow{s,b}{f_{ij^{\prime}}^{jj^{\prime}}} \\
    \node{A_{jj} \cong *} \arrow{e,t}{f_{jj}^{jj^{\prime}}} \node{A_{j j^{\prime}}}
        \end{diagram}
    $$
    are all pushouts in $\cat{C}$.  }
    \end{enumerate}
  Let $\cat{C}$ denote the category with cofibrations from the statement of the theorem, and $wS_{\dotr}(\cat{C})$ the first delooping in the delooping defined in \cite{fW:83}, with the assumption throughout that the weak equivalences are the isomorphisms.    Since $\Gamma$ preserves $co(\cat{C})$ and certainly preserves the isomorphisms and the initial object, it acts on $wS_{\dotr}(\cat{C})$, and we may consider the fixed point space $wS_{\dotr}(\cat{C})^{\Gamma}$, and more generally $wS_{\dotr}(\cat{C})^{\Gamma'}$.  On the other hand, we may consider $wS_{\dotr}(\cat{C}^{\Gamma'})$, where $\cat{C}^{\Gamma'}$ is equipped with the structure of a category with cofibrations and weak equivalences as in the discussion preceding the statement of the proposition. For each $k$ we have an isomorphism of categories $S_{k}(\cat{C})^{\Gamma'} \to S_{k}(\cat{C}^{\Gamma'})$, which gives an isomorphism of simplicial categories $S_{.}(\cat{C})^{\Gamma'} \to S_{.}(\cat{C}^{\Gamma'})$.  We conclude that $I$ in the statement of the theorem is an isomorphism of spectra.
\end{Proof}

Finally, we need to introduce a finality criterion guaranteeing an equivariant equivalence of spectra. We'll  prove it in two steps.  

\begin{Proposition}\label{noneqcriterion}
    Let $\underline{C}$ denote a  category with cofibrations, and let $\underline{C}_0 \subseteq \underline{C}$ denote a full subcategory. Suppose moreover that every object of $\underline{C}$ is isomorphic to an  object of $\underline{C}_0$.   The intersection $co(\underline{C}) \cap Mor(\underline{C}_0)$, which we denote by $co(\underline{C}_0)$,   turns $\underline{C}_0$  into a category with cofibrations, and the the inclusion 
    $$ (\underline{C}_0, co(\underline{C}_0), *) \hookrightarrow (\underline{C}, co (\underline{C}), *)
    $$
    induces an equivalence of $K$-theory spectra.  
\end{Proposition}
\begin{Proof}
    The choice  for any object $x$ of $\underline{C}$ of an isomorphic object $x_0$ of $\underline{C}_0$ provides a functor $\rho \colon \underline{C} \rightarrow \underline{C}_0$, with $\rho (x) = x_0$, and a natural isomorphism from $Id_{\underline{C}}$ to $\rho$.  This means that $N_{\dotr}(i:\underline{C}_0 \rightarrow \underline{C})$ is an equivalence.  It is also clear that a morphism $f$ in $\underline{C}$ is a cofibration if and only if $\rho(f)$ is a cofibration.  It easily follows that the associated level-wise maps on the simplicial categories defining the $K$-theory spectra of the two categories are also equivalences.  
\end{Proof}

Here is the equivariant version of this result. 

\begin{Proposition} \label{eqcriterion} 
Let $\underline{C}$ denote a category with cofibrations acted on by a group $\Gamma$, so that $\Gamma $ preserves the cofibrations. Suppose further that $\underline{C}$ satisfies all the requirements for Propositions \ref{requirements} and \ref{otherrequirements}. Let $\underline{C}_0$ denote a $\Gamma$-invariant  subcategory, which is given the structure of a category with cofibrations by requiring that the  cofibrations in $\underline{C}_0$ are precisely the morphisms in $\underline{C}_0$ which are cofibrations when regarded as morphisms in $\underline{C}$. Suppose further that for every subgroup $\Gamma _0 \subseteq \Gamma$ and every object $c \in \underline{C}^{\Gamma _0}$, there is an object $x \in \underline{C}_0^{\Gamma _0} $ and an  isomorphism $f:c \rightarrow x$ in $\underline{C}^{\Gamma _0} $.  Then the map of spectra $K\underline{C}_0 \rightarrow K \underline{C}$ is an equivariant equivalence of $\Gamma$-spectra. We also remark that if these conditions hold for the Grothendieck (or equivariant) version of $\underline{C}$ from Definition \ref{grothendieck}  also satisfies these hypotheses.  
\end{Proposition}
\begin{Proof}
    This is immediate from Proposition \ref{otherrequirements} and the application of Proposition \ref{noneqcriterion} to the various fixed point subcategories $\underline{C}^{\Gamma _0}$ as $\Gamma _0$ varies over all subgroups of $\Gamma$.  
\end{Proof}

\section{Schlichting's quotients}

We recall Schlichting's quotient construction from \cite{mS:04}. We first define the terms we'll need.  

\begin{Definition} Let $\mathcal{U}$ be a small  exact category, and let $\mathcal{A} \subseteq \mathcal{U}$ be a full subcategory closed under isomorphisms in $\mathcal{U}$.  

\begin{enumerate}
\item{We say  $\mathcal{A} $ is {\em extension closed} if for any exact sequence $$* \rightarrow A \rightarrow B \rightarrow C \rightarrow *$$ in $\mathcal{U}$, in  which $A$ and $C$ are both objects in $\mathcal{A}$, we also have that $B$ is an object of $\mathcal{A}$.  When this condition holds, $\mathcal{A}$ is an exact category in its own right, with the exact sequences consisting of the sequences in $\mathcal{A}$ which are exact in $\mathcal{U}$}
\item{We say $\mathcal{A}$ is {\em right filtering} if 
\begin{itemize}
\item{A is closed under taking admissible subobjects and admissible quotients in U}
\item{Every map $U \rightarrow A$ from an object $U$ of $\mathcal{U}$ to an object $A$ of $\mathcal{A}$ factors through an object $B$ of
$A$ such that the arrow $U \rightarrow  B$ is an admissible epimorphism in $\mathcal{U}$.  }
\end{itemize}}
\item{We say $\mathcal{A}$ is {\em right $s$-filtering} if it is right filtering and if   for every  admissible monomorphism $A \rightarrow U$
from an object $A$ of $\mathcal{A}$  to an object $U$ of $\mathcal{U}$, there is an admissible epimorphism $U \rightarrow B$, with $B$ an object of $\mathcal{A}$, so that the composition $A \rightarrow U \rightarrow B$ is an admissible monomorphism.   }
\item{A morphism in $\mathcal{U}$ is said to be a weak isomorphism if it is a finite composite of morphisms of two types:
\begin{enumerate}
\item{ An admissible monomorphism $i$ which fits into an exact sequence $U_0 \stackrel{i}{\rightarrow} U_1 \rightarrow U_2$, where $U_2$ is an object of $\mathcal{A}$.}
\item{An admissible epimorphism $p$ which fits into an exact sequence $U_0 \rightarrow U_1 \stackrel{p}{\rightarrow} U_2$ where $U_0$ is an object of $\mathcal{A}$. }
\end{enumerate}

The set of weak isomorphisms will be denoted by $\Sigma$.}
\end{enumerate}

\end{Definition}

There is a standard construction which inverts the the elements of $\Sigma$, and which we denote by $\mathcal{U}[\Sigma ^{-1}]$.  Schlichting in \cite{mS:04} proves the following technical result.

\begin{Proposition} If $\mathcal{A}$ is right filtering in $\mathcal{U}$, then the morphisms $\Sigma$ admit a calculus of right fractions (see  \cite{pGmZ:67}).  If $\mathcal{A}$ is right $s$-filtering in $\mathcal{U}$, then the category $\mathcal{U}[\Sigma^{-1}]$ is itself an exact category, with the set of exact sequences consisting of sequences $U_0 \rightarrow U_1 \rightarrow U_2$ which are isomorphic in $\mathcal{U}/\mathcal{A}$ to exact sequences in $\mathcal{U}$. We will denote $\mathcal{U}[\Sigma ^{-1}]$ by $\mathcal{U}/\mathcal{A}$. 
\end{Proposition}

\section{Absence of equivariant fibred Karoubi filtrations} \label{karoubi}

It is known that the $K$-theory assembly map is not an equivalence in general, even for torsion-free groups, without an additional constraint imposed on the coefficient ring.  This is related to the Nil phenomena.  In the course of the argument the spectrum $K(\mathcal{B}_{\Gamma} ( ( CY )^{bdd} )^{\Gamma})$, which is constructed entirely in terms of geometric modules, can be used to replace the $G$-theory constructions of the same type.  This would allow to eschew the introduction of fibred $G$-theory.  An expert may expect this $K$-theory spectrum to be computable using the standard localization and excision techniques based on Karoubi filtrations \cite{mK:70,mCeP:97,dK:15}.  If that was the case, the argument could be finished without passing to $G$-theory and so without requiring constraints on the coefficient ring $A$ or even the group $\Gamma$, therefore ``proving too much''. However, such computation is not possible.   Here we offer a precise explanation.

The Karoubi filtration technique exploits direct sum decompositions of objects of the controlled category. A concise description of Karoubi filtrations is due to Kasprowski \cite{dK:15}.

\begin{Definition} \label{KF}
Let $\mathcal{U}$ be a full additive subcategory of an additive category $\mathcal{A}$. We say that $\mathcal{A}$ is \emph{Karoubi filtered by} $\mathcal{U}$ if for each object $A \in \mathcal{A}$, each $U \in \mathcal{U}$ and any pair of morphisms $f \colon A \to U$, $g \colon U \to A$, there is a decomposition $A \cong E \oplus D$ with $E$ an object of $\mathcal{U}$ and factorizations
\[
\xymatrix{A\ar[rr]^f\ar[rd]_p&&U\\&E\ar[ru]&
}
\quad \text{and}\quad 
\xymatrix{
U\ar[rr]^g\ar[rd]&&A\\&E\ar[ru]_i&
}
\]
where $p$ is the projection onto a direct summand, and $i$ is the inclusion of a direct summand. 
\end{Definition}

Let us restrict to the simplest case of the infinite cyclic group $\Gamma = C$ acting by addition on itself equipped with the standard word metric associated to the generator $1$, denoted $\mathbb{Z}$.  And let us use a shorthand notation $\mathcal{B}_C (\mathbb{Z})^C$ for the relevant fibred category of free $C$-modules.  We show an object of this category has no required direct summands.  Since the coefficient ring $A$ is treated by these methods as a dummy variable, we are free to choose $A$ which contains nilpotent elements.  Let $k$ be a field, then $A = k[x]/(x^2)$.  In this case, where the normal bundle is trivial, we are free to discard the inert extra dimensions and treat the fibre, which should be the cone $CY$, as simply $\mathbb{R}$, and even further as a commensurable subset of integers $\mathbb{Z}$ with the usual metric as above.

We start with a self-contained description of $\mathcal{B}_{C,C} (\mathbb{Z})^C$.

\begin{Definition} \label{WGpre2}
First, we recall the additive category $\mathcal{B}_{C,C} (\mathbb{Z})$ whose
objects are functors
$\theta \colon \underline{E\Gamma} \rightarrow \mathcal{B}_{C} (\mathbb{Z})$
with the additional property that the morphisms $\theta(f)$ are bounded by $0$ but only as homomorphisms between $A$-modules parametrized over the metric space $C$.
The \textit{homotopy fixed points} $\mathcal{B}_C (\mathbb{Z})^{hC}$ is the fixed point category of the action of $C$ on $\mathcal{B}_{C,C} (\mathbb{Z})$ induced from the diagonal action on $C \times \mathbb{Z}$ (see Notation \ref{MM2}). 
Explicitly, cf. Definition 1.3 \cite{gCbG:19}, it is a category with objects which are sets of data $( F,\{ \psi_{\gamma} \} )$ where 
\begin{itemize}
\item $F$ is an object of $\mathcal{B}_{C} (\mathbb{Z})$,
\item $\psi_{\gamma}$ is an isomorphism $F \to \gamma F$ in $\mathcal{B}_{C} (\mathbb{Z})$,
\item $\psi_{\gamma}$ has filtration $0$ but only when viewed as a morphism in $\mathcal{B} (C, \mathrm{Mod}_A)$,
\item $\psi_e = \mathrm{id}$,
\item $\psi_{\gamma_1 \gamma_2} = \gamma_1 \psi_{\gamma_2} \circ \psi_{\gamma_1}$ for all $\gamma_1$, $\gamma_2$ in $C$.
\end{itemize}
The morphisms $( F, \{ \psi_{\gamma} \} ) \to ( F', \{ \psi'_{\gamma} \} )$
can be given by specifying a morphism $\phi_0 \colon F \to F'$ in $\mathcal{B}_{C} (\mathbb{Z})$.
The additive structure on $\mathcal{B}_C (\mathbb{Z})^{hC}$ is induced from that on $\mathcal{B}_{C} (\mathbb{Z})$.
This means in particular that the operation $\oplus$ has the property ${\gamma}(F \oplus G) ={\gamma}F \oplus {\gamma}G$. 
\end{Definition} 

In the present situation, it is possible to exploit the group structure in the fiber $\mathbb{Z}$ in order to simplify the discussion.  
The object we want to construct $(G, \{ \psi_n \})$ will correspond to a particular representation of the cyclic group $C$ in the group of automorphisms of $G$ in $\mathcal{B}_{\mathbb{Z}} (\mathbb{Z})$.
The object $G$ itself can be viewed as a group ring $A[\mathbb{Z} \times \mathbb{Z}]$ equipped with the evident $A$-basis $s^i \otimes t^j$. Note that once we specify an automorphism $\phi$ of $A[\mathbb{Z} \times \mathbb{Z}]$ in $\mathcal{B}_{\mathbb{Z}} (\mathbb{Z})$ which commutes with the diagonal action of $\mathbb{Z}$ and has filtration 0 over the first factor, the maps $\psi_n \colon G \to nG$ can be given as $\phi^n$.
This means that we have to define the action of $\phi$ on $A[\mathbb{Z} \times \mathbb{Z}]$ as an $B$-module, where $B=A[\mathbb{Z}]$ regarded as the diagonal subring. It suffices to specify the action on a $B$-basis. Such a basis is given by the elements $b_i= s^i \otimes 1 = s^i \otimes t^0$.  

\begin{Definition} \label{HWOC}
 The automorphism $\phi$ is defined by the formula
\[
\phi (b_i) = b_{i} + x \xi b_{i-1}.
\]
\end{Definition}

Suppose $F$ is an equivariant direct summand of $G$ which is fibrewise bounded from above.  This means there is an integer $K$ so that any element of $F$ has coefficient 0 for every $A$-basis element $b_i$ with $i \le  K$. Over $R$ this means all coefficients attached to any position $(i,j)$ with $j \ge i +K$ are $0$.  Let $K_0$ be the smallest integer with this property.  This means we can find an element $m$ which has a non-zero coefficient at the $A$-generator $s^{i_0} \otimes t^{i_0 + K_0 -1}$.

\begin{Proposition} \label{HSPDJ}
Given an object $M$ of $\mathcal{B}_X (Y)$, we can think of $M$ as a based free $A$-module with support of each basis element in a single point of $X \times Y$.  Let $M_k$ be the $k$-span of the given basis of $M$, so $M \cong A \otimes_k M_k$.  Given any direct summand $M_0$ of $M_k$, we may construct the direct summand $A \otimes_k M_0$ of $M$.  Every direct summand of $M$ is of the form $A \otimes_k M_0$ for some $M_0$, and $M_0$ is defined uniquely.
\end{Proposition}

\begin{proof}
	We note that $\overline{M} = k \otimes_A M$ is an object of $\mathcal{B}_X (Y,k)$.  Summands of $M$ or of $M$ correspond to idempotent operators on $M$ or, respectively, $\overline{M}$.  Since the ideal $(x)$ is nilpotent, and $A/(x) = k$, idempotent operators on $\overline{M}$ lift uniquely to idempotent operators on $M$, cf. Lemma 10.31.5 and Lemma 15.11.6 \cite{stacks}.   This means that there is a bijection between summands of $M$ and summands of $\overline{M}$, which is isomorphic to $M_k$.  For any summand $M_0$ of $M_k$, we have $M_0 \cong k \otimes_A (A \otimes _k M_0)$, so the correspondence $M_0 \mapsto A \otimes_k M_0$ is a bijection between the summands of $M_k$ and the summands of $M$. 
\end{proof}

Therefore we can assume the coefficient is a unit, or further that it is 1.  Now $\phi(m)$ will have a linear $x$-term in the coefficient at $s^{i_0} \otimes t^{i_0 + K_0}$ so it is not 0, which contradicts the properties of $F$.  We conclude that there are no equivariant summands of $G$ which are fibrewise bounded from above.  These are precisely the summands one requires for a controlled  excision computation of $K(\mathcal{B}_C (\mathbb{Z})^{hC})$.

\section{Comparison with the Farrell-Jones conjecture}

The Farrell-Jones isomorphism conjectures were introduced by F.T.~Farrell and L.E.~Jones  in \cite{tFlJ:93} as an attempt to model the $K$-theory of the group ring in terms of group homology up to the contribution from the virtually cyclic subgroups.
They verified the conjectures for a variety of fundamental groups of manifolds, most notably the nonpositively curved Riemannian manifolds.
A more flexible formulation of the conjectures in \cite{jDwL:98} is in terms of orbit categories for families of subgroups.

A \textit{family of subgroups} $\mathcal{F}$ of a group $\Gamma$ is a nonempty family closed under conjugation and under passage to subgroups. 
Examples of such families are the collections of all finite subgroups $\mathit{Fin}$,
all virtually cyclic subgroups $\mathit{VCyc}$,
and the two extreme cases of the only trivial subgroup $1$ and the family of all subgroups denoted simply by $\Gamma$.
The orbit category of the group $\Gamma$ is the category $\Or_\Gamma$ of all coset spaces of $\Gamma$ regarded 
as $\Gamma$-spaces, so the morphisms between $\Gamma/H$ and $\Gamma/K$ are the $\Gamma$-maps $\Gamma/H \to \Gamma/K$.
Given a family $\mathcal{F}$ of subgroups of $\Gamma$, one defines the full 
subcategory $\Or_\Gamma \mathcal{F}$ on the coset spaces $\Gamma/F$ where $F$ belongs to the family 
$\mathcal{F}$.
So in particular $\Or_\Gamma \Gamma = \Or_\Gamma$.

The canonical maps $K(R) \to K(R\Gamma)$ induce the assembly map 
\[
A(1) \colon \hocolim{\Or_\Gamma 1} K(R) \longrightarrow K(R\Gamma).
\]
One can check that this is precisely the Loday assembly.
For a general family $\mathcal{F}$, there is a way to define a functor 
$K_R \colon \Or_\Gamma \to \textit{spectra}$ so that
for all subgroups $F < \Gamma$ there is a homotopy equivalence 
$K_R (\Gamma/F) \simeq K(RF)$.
Using this functor and the canonical maps $K(RF) \to K(R\Gamma)$,
one obtains the Farrell-Jones assembly map 
\[
A(\mathcal{F}) \colon \hocolim{\Or_\Gamma \mathcal{F}} K_R \longrightarrow K(R\Gamma).
\]
The Farrell-Jones conjecture relative to the family $\mathcal{F}$ is the statement that the assembly map $A(\mathcal{F})$ is an equivalence.

The following diagram illustrates the relationship between various maps.
\[
\xymatrix{
 \hocolim{\Or_\Gamma 1} G(R)  \ar[rr]^-{A_G (1)}
&&G(R\Gamma) 
\\
 \hocolim{\Or_\Gamma 1} K(R) \ar[rr]^-{A(1)} \ar[d]^-{\gamma}
 \ar[u]^-{\hocolimprep \kappa}
&& K(R\Gamma) \ar[u]_-{\kappa_\Gamma}
\\
\hocolim{\Or_\Gamma \mathit{VCyc}} K_R 
\ar@/_1.3pc/[rru]_-<<<<<{\ A(\mathit{VCyc})}  
}
\]
The assembly map relevant to the study of topological rigidity of closed aspherical manifolds is $A(1)$.  Even more specifically, the assembly map of genuine geometric interest is 
\[
A(1) \colon \hocolim{\Or_\Gamma 1} K(\mathbb{Z}) \longrightarrow K(\mathbb{Z}\Gamma),
\]
where $\Gamma$ is the fundamental group of a closed aspherical manifold, and is therefore torsion-free.
Now we observe the following facts.  The Cartan map $\kappa \colon K(R) \to G(R)$ is an equivalence for the regular Noetherian ring $\mathbb{Z}$.  The map $\gamma$ is also an equivalence whenever $R$ is regular Noetherian and $\Gamma$ is torsion-free, which is a much more subtle fact.  Finally, for a very large class of groups $\Gamma$ with strict finite decomposition complexity, the Cartan map $\kappa_\Gamma$ on the right is an equivalence.  So, in this geometrically important and very general situation, the two extreme maps $A_G (1)$ and $A (\mathit{VCyc})$ coincide up to homotopy.
In both lines of research the fact is that studying these maps is technically and conceptually easier than $A (1)$ itself.

The two parts of the diagram above and below $A(1)$ are two different approaches to analyzing the map.  The Farrell-Jones conjecture is an attempt to modify the domain of $A(1)$ with the inherited problems of analyzing and possibly adjusting the domain.  The G-theoretic approach is essentially an attempt to modify the target of $A(1)$ with the inherited need to compute the Cartan map.

It should be interesting to compare these two approaches with the goal of possibly combining the strengths of each method.  The strength of the work on the Farrell-Jones conjecture, even the latest general advances 
\cite{aBwLhR:08a, aBwL:12, aBwLhRhR:14,aB:19}, is in its geometric nature.  The isomorphism results on the Farrell-Jones assembly come with various geometric conditions on the $\Gamma$-spaces but are true for any choice of a ring $R$.  So in order to relate these results to $A (1)$, one simply needs to assume that the coefficient ring $R$ is regular Noetherian, which is a very weak assumption.
On the other hand, the contributions in this direction still use (up to some induction procedures) the flows on $\Gamma$-spaces, and therefore some possibly immanent version of nonpositive curvature.
This geometric nature restricts the range of accessible classes of groups.
In a different way but in much the same fashion, the transition from $A (1)$ to the more algebraic map $A_G (1)$, which is an equivalence for any $\Gamma$ with finite $B\Gamma$ and any Noetherian $R$, requires constraints on both $R$ and $\Gamma$.  However, in this approach the algebraic constraints on $R$ are more strict (finite homological dimension) but the geometric constraints are considerably more lax and abstract.  For example, most likely all groups of strict finite decomposition complexity can be included in the accessible class.
It would be great to find some middle ground where the two approaches through the G-theory assembly and through the Farrell-Jones assembly can be made to interact and hopefully enrich each other.

One intriguing step in this direction can be seen in the work of 
Corti\~{n}as/Cirone \cite{gCeC:14}.
A corollary to their main theorem states that if the Farrell-Jones conjecture is true for a given group and coefficients in any commutative smooth $\mathbb{Q}$-algebra $R$ then it is also true for the same group and any commutative $\mathbb{Q}$-algebra.
Since every smooth $\mathbb{Q}$-algebra $R$ is Noetherian regular of finite homological dimension, the Cartan map $\kappa_\Gamma$ in the diagram is an equivalence for a geometric group $\Gamma$ of finite decomposition complexity.
The fact that the assembly map $A_G (1)$ is an equivalence gives the Farrell-Jones in the special smooth cases, therefore $A(\mathit{VCyc})$ is also an equivalence for these groups and coefficients in any commutative $\mathbb{Q}$-algebra.



\begin{thebibliography}{99}


\bibitem{aB:19}
{A.~Bartels},
{\it K-theory and actions on Euclidean retracts},
Proceedings of the ICM 2018 (2019), 1041--1062. 

\bibitem{aBwLhR:08a}
A.~Bartels, W.~L\"{u}ck, and H.~Reich, 
\textit{The K-theoretic Farrell-Jones Conjecture for hyperbolic groups},
Invent. Math. (2008), 29--70.

\bibitem{aBwL:12}
A.~Bartels and W.~L\"{u}ck
\textit{The Borel Conjecture for hyperbolic and CAT(0)-groups},
Ann. Math. 175 (2012), 631--689.

\bibitem{aBwLhRhR:14}
A.~Bartels, W.~L\"{u}ck, H.~Reich, and H.~R\"{u}ping,
\textit{K- and L-theory of group rings over $GL_n (\mathbb{Z})$}, 
Publ. Math. Inst. Hautes \'{E}tudes Sci.  \textbf{119} (2014), 97--125.



\bibitem{aBdK:72} A.K. Bousfield and D.M. Kan, {\em  Homotopy limits, completions and localizations}
Lecture Notes in Math., Vol. 304
Springer-Verlag, Berlin-New York, 1972, v+348 pp.

\bibitem{mCeP:97}
{M. Cardenas and E.K. Pedersen},
{\it On the Karoubi filtration of a category},
$K$-theory, {\bf 12} (1997), 165--191.

\bibitem{gC:95}
{G. Carlsson},
{\it Bounded $K$-theory and the assembly map
in algebraic $K$-theory},
in {\it Novikov conjectures, index theory and rigidity},
{\it Vol. 2}
(S.C. Ferry, A. Ranicki, and J. Rosenberg, eds.),
Cambridge U. Press (1995), 5--127.

\bibitem{gC:95.1}\bysame, {\it On the algebraic K-theory of infinite product categories},
$K$-Theory v. 9, n. 4 (July, 1995): 305--322

\bibitem{gCbG:03}
{G.~Carlsson and B.~Goldfarb}, \textit{On homological coherence of
discrete groups}, J.~Algebra {\bf 276} (2004), 502--514.

\bibitem{gCbG:04}
\bysame, \textit{The integral K-theoretic Novikov conjecture for
groups with finite asymptotic dimension}, Invent. Math. {\bf 157} (2004), 405--418.

\bibitem{gCbG:11}
\bysame,
\textit{Controlled algebraic $G$-theory, I}, J. Homotopy Relat. Struct. \textbf{6} (2011), 119--159.

\bibitem{gCbG:16}
\bysame, \textit{On modules over infinite group rings}, Int. J. Algebra Comput.\,\textbf{26} (2016), 1--16.

\bibitem{gCbG:18}
\bysame, \textit{Bounded $G$-theory with fibred control}, J. Pure Appl. Algebra \textbf{223} (2019), 5360--5395. 

\bibitem{gCbG:19}
\bysame, \textit{Excision in equivariant fibred \textit{G}-theory}, Annals of K-theory \textbf{5} (2020), 721--756.

\bibitem{gCeP:95}  G. Carlsson and E.K. Pedersen, \textit{Controlled algebra and the Novikov conjectures for $K$- and $L$-theory} Topology \textbf{4}, 3, (1995), 731--758.

\bibitem{mCiJ:98} M. Crabb and I. James, \textit{Fibrewise homotopy theory}
Springer Monogr. Math.
Springer-Verlag London, Ltd., London, 1998, viii+341 pp.



\bibitem{gCeC:14}
G.~Corti\~{n}as and E.R.~Cirone,
\textit{Singular coefficients in the K-theoretic Farrell-Jones conjecture}, Algebr. Geom. Topol. \textbf{16} (2016), 129--147.

\bibitem{jDwL:98}
J.F. Davis and W. L\"{u}ck, 
\textit{Spaces over a category and assembly maps in isomorphism conjectures in K-and L-theory}, 
K-Theory \textbf{15} (1998), 201--252.

\bibitem{tFlJ:93}
F.T. Farrell and L. Jones,
\textit{Isomorphism conjectures in algebraic \textit{K}-theory},
Journal of Amer. Math. Soc. \textbf{6} (1993), 249--297.


\bibitem{pGmZ:67}  
P. Gabriel and M. Zisman, {\it Calculus of fractions and homotopy theory}
Ergeb. Math. Grenzgeb., Band 35[Results in Mathematics and Related Areas],Springer Verlag, 1967, 168 pp. 


\bibitem{bG:13} {B. Goldfarb}, \textit{Weak coherence of groups and finite decomposition complexity}, preprint (2013). \texttt{arXiv:1307.5345v2}

\bibitem{mK:70}
M. Karoubi, \textit{Foncteur d\'{e}riv\'{e}es et K-th\'{e}orie}, in S\'{e}minaire Heidelberg-Saarbr\"{u}cken-Strasbourg sur la $K$-th\'{e}orie (1967/68), Lecture Notes in Mathematics, vol. 136, Springer, Berlin, 1970, 107--186.

\bibitem{dK:15}
D. Kasprowski, \textit{On the K-theory of groups with finite decomposition complexity},
Proc. Lond. Math. Soc. \textbf{110} (2015), 565--592.

\bibitem{dKcW:20}
D. Kasprowski and C. Winges, 
\textit{Shortening binary complexes and commutativity
of k-theory with infinite products},
Trans. Amer. Math. Soc. Series B \textbf{7} (2020), 1--23.

\bibitem{yK:60}
Y. Kawada, 
{\em Cosheaves}, Proc. Japan Acad.36(1960), 81–85.

\bibitem{bK:96} B. Keller, \textit{ Derived categories and their uses}, in: Handbook of Algebra, Vol. 1, North-Holland, Amsterdam, 1996, pp. 671–701.




\bibitem{wL:22}
W. L\"{u}ck, {\it Isomorphism conjectures in K- and L-theory}, book preprint, 2022.

\bibitem{cMmM:19}
C. Malkiewich and M. Merling,
\textit{Equivariant \textit{A}-theory},  Documenta Mathematica \textbf{24} (2019), 815--855.

\bibitem{jMjS:06} J.P.May and J. Sigurdsson,  \textit{Parametrized homotopy theory}
Math. Surveys Monogr., 132
American Mathematical Society, Providence, RI, 2006, x+441 pp.

\bibitem{mM:16}
{M. Merling},
\textit{Equivariant algebraic K-theory of G-rings},
Math. Z.
\textbf{285} (2017), 1205--1248.

\bibitem{ePcW:85}
{E.K. Pedersen and C. Weibel},
{\it A nonconnective delooping of algebraic $K$-theory},
in {\it Algebraic and geometric topology}
(A. Ranicki, N. Levitt, and F. Quinn, eds.),
Lecture Notes in Math. {\bf 1126},
Springer-Verlag (1985), 166--181.

\bibitem{ePcW:89}
\bysame,
\textit{K-theory homology of spaces}, 
in \textit{Algebraic Topology} (A. Dold and B. Eckmann, eds.), Lecture Notes in Math. {\bf 1370}, Springer-Verlag (1989), 346--361.

\bibitem{dQ:72} D. Quillen,  \textit{Higher algebraic  K -theory. I}
Lecture Notes in Math., Vol. 341
Springer-Verlag, Berlin-New York, 1973, pp. 85–147.


\bibitem{aR:80} A.A. Ranicki, \textit{Algebraic theory of surgery II. Applications to topology.} Proc. Lond. Math. Soc.
40, 193–283 (1980)


\bibitem{aR:16}
A. Ranicki, 
\textit{The birth of the Borel conjecture}, a historical note on the personal web site \texttt{https://webhomes.maths.ed.ac.uk/\string~v1ranick//surgery/borel.pdf}, 2016.

\bibitem{hRmV:19}
H. Reich and M. Varisco, 
\textit{Algebraic K-theory, assembly maps, controlled algebra, and trace methods},
in \textit{Space--Time--Matter. Analytic and Geometric Structures}, De Gruyter, 2018, 1--50.

\bibitem{jR:03}
{J. Roe},
{\it Lectures on coarse geometry}, University Lecture Series, vol.~31, American Mathematical Society, 2003.


\bibitem{mS:04} M. Schlichting, \textit{Delooping the K-theory of exact categories.} Topology 43 (2004), no.5, 1089–1103.


\bibitem{rT:82}
 R.W. Thomason, \textit{First quadrant spectral sequences in algebraic  K -theory via homotopy colimits}
Comm. Algebra 10 (1982), no. 15, 1589–1668


\bibitem{stacks}
Various authors, \textit{The Stacks project}, \texttt{stacks.math.columbia.edu}, 2018.


\bibitem{fW:83}
{F.~Waldhausen},
\textit{Algebraic \textit{K}-theory of spaces}, in Proc. of 1983 Rutgers Conf. on Algebraic Topology, Springer Lecture Notes in Math. \textbf{1126} (1985), 318--419.
\bibitem{cW:79} C.A. Weibel, \textit{Nilpotence and $K$-theory}, Journal of Algebra {\bf 61}, (1979), 298-307. 




\bibitem{sW:22}
S. Weinberger, {\it Variations on a Theme of Borel}, Cambridge U. Press, 2023.
\end{thebibliography}
\end{document}